\documentclass[12pt]{article}

\usepackage{amsmath,enumerate,amsbsy,amsfonts,amssymb,mathabx,amscd,graphicx,algorithm,amstext,euscript,multirow}
\usepackage[round,authoryear]{natbib}
\usepackage{authblk}
\usepackage{mathrsfs}
\usepackage{appendix}
\usepackage{chngcntr}
\usepackage{epsfig}
\newtheorem{definition}{Definition}
\newtheorem{theorem}{Theorem}

\newtheorem{lemma}{Lemma}

\newtheorem{proposition}{Proposition}
\newtheorem{condition}{Condition}
\newtheorem{remark}{Remark}

\newcommand{\beps}{\boldsymbol \epsilon}

\newcommand{\bSigma}{\boldsymbol \Sigma}

\newcommand{\bGamma}{\boldsymbol \Gamma}
\newcommand{\bTheta}{\boldsymbol \Theta}

\newcommand{\bDelta}{\boldsymbol \Delta}

\newcommand{\bpsi}{\boldsymbol \psi}
\newcommand{\bPsi}{\boldsymbol \Psi}
\newcommand{\bphi}{\boldsymbol \phi}
\newcommand{\bPhi}{\boldsymbol \Phi}

\newcommand{\bxi}{\boldsymbol \xi}
\newcommand{\btheta}{\boldsymbol \theta}

\newcommand{\bbeta}{\boldsymbol \beta}
\newcommand{\bnu}{\boldsymbol \nu}

\newcommand{\bvar}{\boldsymbol \varepsilon}
\newcommand{\bvarphi}{\boldsymbol \varphi}
\newcommand{\bvarPhi}{\boldsymbol \varPhi}

\newcommand{\bvarepsilon}{\boldsymbol \varepsilon}

\newcommand{\bzeta}{\boldsymbol \zeta}
\newcommand{\bPi}{\boldsymbol \Pi}

\newcommand{\bbf}{{\boldsymbol f}}
\newcommand{\ba}{{\mathbf a}}
\newcommand{\be}{{\mathbf e}}
\newcommand{\bx}{{\mathbf x}}
\newcommand{\by}{{\mathbf y}}

\newcommand{\br}{{\mathbf r}}

\newcommand{\bg}{{\mathbf g}}

\newcommand{\bD}{{\bf D}}
\newcommand{\bA}{{\bf A}}
\newcommand{\bB}{{\bf B}}
\newcommand{\bC}{{\bf C}}
\newcommand{\bE}{{\bf E}}
\newcommand{\bI}{{\bf I}}
\newcommand{\bK}{{\bf K}}
\newcommand{\bM}{{\bf M}}

\newcommand{\bX}{{\bf X}}
\newcommand{\bY}{{\bf Y}}
\newcommand{\bZ}{{\bf Z}}
\newcommand{\bR}{{\bf R}}
\newcommand{\bU}{{\bf U}}
\newcommand{\bV}{{\bf V}}
\newcommand{\bQ}{{\bf Q}}

\newcommand{\bW}{{\bf W}}

\newcommand{\cM}{{\cal M}}

\newcommand{\cU}{{\cal U}}
\newcommand{\cV}{{\cal V}}
\newcommand{\cS}{{\cal S}}
\newcommand{\cZ}{{\cal Z}}

\newcommand{\eZ}{\mathbb{Z}}
\newcommand{\eR}{\mathbb{R}}
\newcommand{\cH}{\mathbb{H}}

\newcommand{\tE}{\text{E}}
\newcommand{\tF}{\text{F}}

\newcommand{\cov}{\text{Cov}}
\newcommand{\var}{\text{Var}}

\newcommand{\LS}{\text{LS}}
\newcommand{\OLS}{\text{OLS-O}}

\newcommand{\eulB}{\EuScript B}

\newcommand{\llangle}{\langle \langle}
\newcommand{\rrangle}{\rangle \rangle}

\newcommand{\T}{\scriptscriptstyle T}

\makeatletter
\DeclareRobustCommand\widecheck[1]{{\mathpalette\@widecheck{#1}}}
\def\@widecheck#1#2{%
	\setbox\z@\hbox{\m@th$#1#2$}%
	\setbox\tw@\hbox{\m@th$#1%
		\widehat{%
			\vrule\@width\z@\@height\ht\z@
			\vrule\@height\z@\@width\wd\z@}$}%
	\dp\tw@-\ht\z@
	\@tempdima\ht\z@ \advance\@tempdima2\ht\tw@ \divide\@tempdima\thr@@
	\setbox\tw@\hbox{%
		\raise\@tempdima\hbox{\scalebox{1}[-1]{\lower\@tempdima\box
				\tw@}}}%
	{\ooalign{\box\tw@ \cr \box\z@}}}
\makeatother
\RequirePackage[colorlinks,citecolor=blue,urlcolor=blue]{hyperref}

\begin{document}

\title{\bf \large Finite Sample Theory for High-Dimensional Functional/Scalar Time Series with Applications}


\author[1]{Qin Fang}
\author[2]{Shaojun Guo}
\author[1]{Xinghao Qiao}
\affil[1]{Department of Statistics, London School of Economics and Political Science, U.K.}
\affil[2]{Institute of Statistics and Big Data, Renmin University of China, P.R. China}
\setcounter{Maxaffil}{0}

\renewcommand\Affilfont{\itshape\small}
\date{\vspace{-5ex}}
\maketitle

\begin{abstract}
Statistical analysis of high-dimensional functional times series arises in various applications. Under this scenario, in addition to the intrinsic infinite-dimensionality of functional data, the number of functional variables can grow with the number of serially dependent observations. In this paper, we focus on the theoretical analysis of relevant estimated cross-(auto)covariance terms between two multivariate functional time series or a mixture of multivariate functional and scalar time series beyond the Gaussianity assumption. We introduce a new perspective on dependence by proposing functional cross-spectral stability measure to characterize the effect of dependence on these estimated cross terms, which are essential in the estimates for additive functional linear regressions. With the proposed functional cross-spectral stability measure, we develop useful concentration inequalities for estimated cross-(auto)covariance matrix functions to accommodate more general sub-Gaussian functional linear processes and, furthermore, establish finite sample theory for relevant estimated terms under a commonly adopted functional principal component analysis framework. Using our derived non-asymptotic results, we investigate the convergence properties of the regularized estimates for two additive functional linear regression applications under sparsity assumptions including functional linear lagged regression and partially functional linear regression in the context of high-dimensional functional/scalar time series.
\end{abstract}


\section{Introduction}
\label{sec.intro}
Functional time series have received a great deal of attention in the last decade in order to provide methodology for functional data objects that are observed sequentially over time. Despite progress being made in this area, existing literature has focused on the statistical analysis of a single or small number of random functions.
The increasing availability of large dataset with multiple functional features corresponds to the data structure of 
\begin{equation*}
\label{mfts}
\bX_t(u)=\big(X_{t1}(u), \dots, X_{tp}(u)\big)^{\T}, ~~t=1,\dots, n,u \in \cU,
\end{equation*}
with covariance matrix function $\bSigma_0^X(u,v)=\cov\{\bX_t(u), \bX_t(v)\},$
where, under the high-dimensional and dependent setting, the number of functional variables ($p$) can be comparable to, or even larger than, the number of serially dependent observations ($n$), posing new challenges to existing work. 

Examples of high-dimensional functional time series include daily electricity consumption curves \cite[]{cho2013} for a large collection of households, half-hourly measured PM10 curves \cite[]{aue2015} over a large number of sites and cumulative intraday return curves \cite[]{horvath2014} for hundreds of stocks. These applications require developing learning techniques to handle such type of data. One large class considers imposing various functional sparsity assumptions  on the model parameter space, e.g. {\it vector functional autoregressions} (VFAR) \cite[]{guo2019} and, under a special independent setting, functional graphical models \cite[]{qiao2019} and functional additive regressions \cite[]{fan2014,fan2015,kong2016,luo2017,xue2020}, where the corresponding regularized estimates are proposed.

Within the high-dimensional time series framework, it is essential to 
establish necessary concentration inequalities for dependent data and assess how the presence of serial dependence affect non-asymptotic error bounds. See relevant concentration results for Gaussian process \cite[]{basu2015a},
linear process or linear spatio-temporal model with more general noise distributions \cite[]{sun2018,shu2019} and heavy tailed time series \cite[]{wong2019}. Compared with theoretical analysis of scalar time series, the added technical challenges that arise to handle functional time series involve developing non-asymptotic results for dependent processes within an abstract Hilbert space and characterizing the effect of serial dependence in $\{\bX_t(\cdot)\}$ with infinite, summable and decaying eigenvalues of $\bSigma_0^X.$

Theoretical investigation of high-dimensional functional time series is rather incomplete. \cite{guo2019} proposed a functional stability measure for Gaussian functional time series by controlling the functional Rayleigh quotients of spectral density matrix functions relative to $\bSigma_0^X$ and hence can precisely capture the effect of small eigenvalues. Moreover, they relied on it to establish concentration bounds on sample (auto)covariance matrix function of $\bX_t(\cdot),$ serving as a fundamental tool to provide theoretical guarantees for the proposed three-step procedure and the regularized VFAR estimate, in a high dimensional regime. However, their proposed stability measure only facilitates finite sample theory to accommodate Gaussian functional time series and is not sufficient to evaluate the effect of serial dependence on the estimated cross-(auto)covariance terms in a non-asymptotic way, which plays a crucial role in the theoretical analysis of a wide class of additive functional linear regressions under the high-dimensional regime when the serial dependence exists.


To illustrate, we consider two important examples of additive functional linear regressions in the context of high-dimensional functional/scalar time series.
The first example considers the high-dimensional extension of functional linear lagged regression \cite[]{hormann2015b} in the additive form:
\begin{equation} \label{model_fully}
    Y_t(v) = \sum_{h=0}^L \sum_{j=1}^p \int_{\cU}X_{(t-h)j}(u)\beta_{hj}(u,v)du + \epsilon_t(v),~~t=L+1,\dots,n, (u,v) \in \cU \times \cV,
\end{equation}
where $p$-dimensional functional covariates $\{\bX_t(\cdot)\}$ and functional errors $\{\epsilon_t(\cdot)\}$ are generated from independent, centered, stationary functional processes, and $\{\beta_{hj}(\cdot,\cdot):h=0, \dots, L, j=1, \dots, p\}$ are sparse functional coefficients to be estimated.
Under an independent setting without lagged functional covariates, model~(\ref{model_fully}) reduces to the additive function-on-function linear regression \cite[]{luo2017}.

The second example studies partially functional linear regression \cite[]{kong2016} consisting of a mixture of $p$-dimensional functional time series $\{\bX_t(\cdot)\}$ and $d$-dimensional scalar time series $\bZ_t=(Z_{t1}, \dots, Z_{td})^{\T}$ for $t=1, \dots, n,$ both of which are independent of errors $\{\epsilon_t\},$ as follows:
\begin{equation} \label{model_scalar}
    Y_t = \sum_{j=1}^p \int_{\cU}X_{tj}(u)\beta_{j}(u)du + \sum_{k=1}^d Z_{tk}\gamma_{k} + \epsilon_t,  ~~t=1,\dots,n, u \in \cU,
\end{equation}
where $\{\beta_j(\cdot): j=1\dots,p\}$ are sparse functional coefficients and $\{\gamma_k: k=1, \dots,d\}$ are sparse scalar coefficients. Whereas \cite{kong2016} focused on an independent scenario and treated $p$ as fixed, we allow both $p$ and $d$ to be diverging with $n$ under a more general dependence structure. See also special cases of model~(\ref{model_scalar}) without functional covariates or scalar covariates in \cite{basu2015a,Wu2016} or \cite{fan2015,xue2020}, respectively.

In addition to existing non-asymptotic results in \cite{guo2019}, the central challenge to provide theoretical supports for the regularized estimates for models~(\ref{model_fully}) and (\ref{model_scalar}) is: (i) to characterize how the underlying dependence structure affects the non-asymptotic error bounds on those essential estimated cross-(auto)covariance terms, e.g. estimated cross-covariance functions between $\bX_t(\cdot)$ and $Y_{t+h}(\cdot)$ (or $\epsilon_{t+h}(\cdot)$) for $h=0, \dots, L$ in model~(\ref{model_fully}) and estimates of $\cov(\bX_t(\cdot),\bZ_t),$  $\cov(\bZ_t,\epsilon_t)$ and $\cov(\bX_t(\cdot),\epsilon_t)$ in model~(\ref{model_scalar}); (ii) to develop useful non-asymptotic results beyond Gaussian functional/scalar time series.

To address such challenges, the main contribution of our paper is threefold.
\begin{itemize}
    \item First, we propose a novel functional cross-spectral stability measure between $\{\bX_t(\cdot)\}$ and $d$-dimensional functional (or scalar) time series, i.e. $\{\bY_t(\cdot)=(Y_{t1}(\cdot), \dots, Y_{td}(\cdot))^{\T}\},$ defined on $\cV$ or $\{\bZ_t\},$ based on their cross-spectral density properties. Compared with the direct functional extension of the cross-stability measure in \cite{basu2015a}, our functional cross-spectral stability measure can more precisely capture the small eigenvalues effect to handle truly infinite-dimensional functional objects. It also facilitates the development of non-asymptotic results for $\widehat\bSigma^{X,Y}_h$ and $\widehat \bSigma^{X,Z}_h,$ which respectively  are estimates of cross-(auto)covariance terms, $\bSigma^{X,Y}_h(u,v)=\cov(\bX_t(u),\bY_{t+h}(v))$ and $\bSigma^{X,Z}_h=\cov(\bX_t(u),\bZ_{t+h})$ for all integer $h.$ Moreover, it provides insights into how  $\widehat\bSigma^{X,Y}_h$ and $\widehat\bSigma^{X,Z}_h$ are affected by the presence of serial dependence.
    
    \item Second, we establish finite sample theory in a non-asymptotic way for relevant estimated (cross)-(auto)covariance terms beyond Gaussian functional (or scalar) time series to accommodate more general multivariate functional linear processes with sub-Gaussian functional errors. Our finite sample results and adopted techniques are general, and can be applied broadly to provide theoretical guarantees for regularized estimates of other high-dimensional functional time series models, e.g. the autocovariance-based estimates of sparse functional linear regressions \cite[]{chang2020} and the functional factor model \cite[]{guo2019}.
    
    \item Third, due to the infinite dimensionality of the functional covariates, dimension reduction is necessary in the estimation. One common approach is {\it functional principal component analysis} (FPCA). We hence establish useful deviation bounds on relevant estimated terms under a FPCA framework. To illustrate using models~(\ref{model_fully}) and (\ref{model_scalar}), we implement FPCA-based three-step procedures to estimate unknown parameters under sparsity constraints. With derived non-asymptotic results, we verify functional analogs of routinely used restricted eigenvalue and deviation conditions in the lasso literature \cite[]{loh2012,basu2015a} and, furthermore, investigate the convergence properties of regularized estimates under a high-dimensional and serially dependent setting.
\end{itemize}

{\bf Literature review}. Our work lies in the intersection of two strands of literature: functional time series and high-dimensional time series. In the context of functional time series, many standard univariate or low-dimensional time series methods have been recently adapted to the functional domain with theoretical properties explored from a standard asymptotic perspective, see, e.g., \cite{Bbosq1,bathia2010,hormann2010,panaretos2013,aue2015,hormann2015b,pham2018,li2019} and reference therein. In the context of high-dimensional time series, some lower-dimensional structural assumptions are often incorporated on the model parameter space and different regularized estimation procedures have been developed for the respective learning tasks including, e.g.,
high-dimensional sparse linear regression \cite[]{basu2015a,Wu2016,han2020} and high-dimensional sparse vector autoregression \cite[]{Guo2016,Lin2017,Gao2019,Ghosh2019,Zhou2019,wong2019,Lin2020}.

{\bf Outline}. The remainder of the paper is organized as follows. In Section~\ref{sec_non-asymptotic results}, we propose cross-stability measures under functional and mixed-process scenarios, define sub-Gaussian functional linear processes and rely on them to present finite sample theory for estimated (cross-)terms used in subsequent analyses. In Section~\ref{sec_fully}, we consider sparse high-dimensional functional linear lagged model in (\ref{model_fully}), develop the penalized least squares estimation procedure and apply our derived non-asymptotic results to provide theoretical guarantees for the estimates. Section~\ref{sec_partial} is devoted to the modelling, regularized estimation and application of established deviation bounds on the theoretical analysis of sparse high-dimensional partially functional linear model in (\ref{model_scalar}). Finally, we examine the finite-sample performance of the proposed methods for both models (\ref{model_fully}) and (\ref{model_scalar}) through simulation studies in Section~\ref{sec.sim}.
All technical proofs are relegated to the appendix.


{\bf Notation}. Let $\mathbb{Z}$ and $\mathbb{R}$ denote the sets of integers and real numbers, respectively. For $x, y \in \mathbb{R},$ we use $x \vee y = \max(x,y).$ 
For two positive sequences $\{a_n\}$ and $\{b_n\},$ we write $a_n \lesssim b_n$ or $a_n=O(b_n)$ or $b_n \gtrsim a_n$ if there exists a positive constant $c$ independent of $n$ such that $a_n/b_n \leq c .$ We write $a_n \asymp b_n$ if $a_n \lesssim b_n$ and $a_n \gtrsim b_n.$ For a vector $\bx \in \mathbb{R}^p,$ we denote its $\ell_1$, $\ell_2$  and maximum norms  by $\|\bx\|_1 = \sum_{j = 1}^p |x_j|,$  $\|\bx\| = (\sum_{j = 1}^p |x_j|^2)^{1/2}$ and $||\bx||_{\max}=\max_{j}|x_{j}|,$ respectively. For a matrix $\bB \in {\eR}^{p \times q},$ we denote its Frobenius norm by $||\bB||_{\tF}=\big(\sum_{i,j} \text{B}_{ij}^2\big)^{1/2}.$  Let  $L_2(\cU)$ be a Hilbert space of square integrable functions on a compact interval $\cU.$ For $f,g \in L_2(\cU),$ we denote the inner product by $\langle f,g \rangle=\int_{\cU}f(u)g(u)du$ for $f,g \in L_2(\cU)$ with the norm $\|\cdot\|=\langle \cdot,\cdot \rangle^{1/2}.$ For a Hilbert space $\mathbb{H} \subseteq L_2(\cU),$ we denote the $p$-fold Cartesian product by $\mathbb{H}^p = \mathbb{H} \times \cdots \times \mathbb{H}$ and the tensor product $\mathbb{S} = \mathbb{H} \otimes \mathbb{H}.$ For $\bbf=(f_1, \ldots, f_p)^{\T}$ and $\bg=(g_1, \dots,g_p)^{\T}$ in $\cH^p,$ we denote the inner product by
$\langle\bbf,\bg\rangle=\sum_{i=1}^p\langle f_i, g_i\rangle$ with induced norm of $\bbf$ by $\|\bbf\| = \langle \bbf, \bbf \rangle^{1/2},$ $\ell_1$ norm by $\|\bbf\|_1 = \sum_{i=1}^p \|f_i\|,$ and $\ell_0$ norm by $\|\bbf\|_0 = \sum_{i=1}^p I(\|f_i\| \neq 0),$ where $I(\cdot)$ is the indicator function.
For an integral matrix operator $\bK: \mathbb{H}^p \rightarrow \mathbb{H}^q$ induced from the kernel matrix function $\bK=(K_{ij})_{q \times p}$ with each $K_{ij} \in \mathbb{S}$ through
$\bK(\bbf)(u)=\Big(\sum_{j=1}^p\langle K_{1j}(u,\cdot),f_{j}(\cdot)\rangle, \dots, \sum_{j=1}^p\langle K_{qj}(u,\cdot),f_{j}(\cdot)\rangle\Big)^{\T} \in \mathbb{H}^q,$
for any given $\bbf \in \mathbb{H}^p.$ To simplify notation, we will use $\bK$ to denote both the kernel function and the operator. When $p=q=1,$ $\bK$ degenerates to $K$ and we denote its Hilbert--Schmidt norm by $\|K\|_{\cS} =\big(\int\int K(u,v)^2dudv\big)^{1/2}.$ For general $\bK,$ we define functional versions of 
Frobenius, elementwise $\ell_{\infty}$, matrix $\ell_{1}$ and matrix $\ell_{\infty}$ norms by 
$\|\bK\|_{\tF} = \big(\sum_{i,j} \|K_{ij}\|_{\cS}^2\big)^{1/2},$ $\|\bK\|_{\max} =  \max_{i,j} \|K_{ij}\|_{\cS},$ $\|\bK\|_{1} = \max_{j} \sum_{i} \|K_{ij}\|_{\cS}$ and $\|\bK\|_{\infty} = \max_{i} \sum_{j} \|K_{ij}\|_{\cS},$ respectively. 

\section{Finite sample theory} \label{sec_non-asymptotic results}
In this section, we first review functional stability measure and propose functional cross-spectral stability measure. We then introduce the definitions of sub-Gaussian process and multivariate functional linear process. Finally, we rely on our proposed stability measures to develop finite sample theory for useful estimated terms to accommodate sub-Gaussian functional linear processes.

\subsection{Functional stability measure}
Consider a $p$-dimensional vector of weakly stationary functional time series $\{\bX_t(\cdot)\}_{t\in\mathbb{Z}}$ defined on $\cU,$ with mean zero and $p \times p$ autocovariance matrix functions, 
\begin{equation} \nonumber
    \bSigma_h^X(u,v) = \cov\{\bX_t(u),\bX_{t+h}(v)\} = \{\Sigma_{h,jk}^X(u,v)\}_{1\leq j,k\leq p},  ~~ t,h \in \mathbb{Z},(u,v)\in \mathcal{U}^2.
\end{equation}


These autocovariance matrix functions (or operators) encode the second-order dynamical properties of $\{\bX_t(\cdot)\}$ and typically serve as the main focus of functional time series analysis. From a frequency domain analysis prospective, spectral density matrix function (or operator) aggregates autocovariance information at different lag orders $h \in \mathbb Z$ at a frequency $\theta\in[-\pi,\pi]$ as
$$\bbf_\theta^X= \frac{1}{2\pi}\sum_{h\in \mathbb{Z}}\bSigma_h^X\mathrm{exp}(-ih\theta).$$ According to \cite{guo2019}, the functional stability measure of $\{\bX_t(\cdot)\}$ is defined based on the functional Rayleigh quotients of $\bbf_\theta^X$ relative to $\bSigma_0^X,$ 
\begin{equation} \label{eq_M_X}
\cM^X = 2\pi \mathop{\text{ess sup}}\limits_{\theta \in [-\pi,\pi],\bPhi\in \mathbb{H}_0^p}\frac{\langle\bPhi,\bbf_\theta^X(\bPhi)\rangle}{\langle\bPhi,\bSigma_0^X(\bPhi)\rangle},
\end{equation}
where $\mathbb{H}_0^p = \{\bPhi \in \mathbb{H}^p:\langle\bPhi,\bSigma_0^X(\bPhi)\rangle \in (0, \infty)\}.$ To handle truly infinite-dimensional objects $\{\bX_t(\cdot)\}$ with infinite, summable and decaying eigenvalues of $\bSigma_{0}^X,$ such stability measure $\cM^X$ can more precisely capture the effect of small eigenvalues of $\bSigma_0^X$ on the numerator in (\ref{eq_M_X}).




We next impose a condition on $\cM^X$ and introduce the functional stability measure of subprocesses of $\{\bX_t(\cdot)\},$ which will be used in our subsequent analysis.
\begin{condition} \label{con_stability}
(i) The spectral density matrix operator $\bbf^X_\theta, \theta \in [-\pi,\pi]$ exists;
(ii) $\cM^X<\infty.$
\end{condition}


For any $k$-dimensional subset $J\subseteq\{1,\dots,p\}$ with its cardinality $|J| \leq k,$ the functional stability measure of $\big\{\big(X_{tj}(\cdot)\big):j \in J\big\}_{t \in \eZ}$ is defined by
\begin{equation} \label{def_sub_fsm}
\cM_k^X = 2\pi \cdot \underset{\theta\in [-\pi, \pi],\|\bPhi\|_0 \le k,\bPhi \in \cH_0^p}{\text{ess} \sup} \frac{\langle \bPhi, \bbf_{\theta}^X(\bPhi)\rangle}{\langle \bPhi, \bSigma_0^X(\bPhi)\rangle}, ~k=1,\dots,p.
\end{equation} Under Condition~\ref{con_stability}, we have $\cM_k^X \leq \cM^X <\infty.$ 

\subsection{Functional cross-spectral stability measure}
Consider $\{\bX_t(\cdot)\}$ and $\{\bY_t(\cdot)\},$ where $\{\bY_t(\cdot)\}_{t\in\mathbb{Z}}$ is a $d$-dimensional vector of centered and weakly stationary functional time series, defined on $\cV,$ with lag-$h$ autocovariance matrix function given by
$$\bSigma_h^Y(u,v) = \cov\{\bY_t(u),\bY_{t+h}(v)\} = \{\Sigma_{h,jk}^Y(u,v)\}_{1\leq j,k\leq d},~~t,h \in \mathbb{Z},(u,v)\in  \cV^2.$$

To characterize the effect of dependence on the cross-covariance between two sequences of joint stationary multivariate functional time series, we can correspondingly define the cross-spectral density matrix function (or operator) and functional cross-spectral stability measure. The proposed cross-spectral stability measure plays a crucial role in the non-asymptotic analysis of relevant estimated cross terms, e.g., estimated cross-(auto)covariance matrix functions in Section~\ref{sec_cross cov}.
\begin{definition}
The cross-spectral density matrix function  between $\{\bX_t(\cdot)\}_{t\in\mathbb{Z}}$ and $\{\bY_t(\cdot)\}_{t\in\mathbb{Z}}$ is defined by
\begin{equation}\nonumber
\bbf_\theta^{X,Y} = \frac{1}{2\pi}\sum_{h\in \mathbb{Z}}\bSigma_h^{X,Y}\mathrm{exp}(-ih\theta),~~~ \theta\in[-\pi,\pi],
\end{equation}
where $\bSigma_h^{X,Y}(u,v) = \cov\{\bX_t(u),\bY_{t+h}(v)\} = \{\Sigma_{h,jk}^{X,Y}(u,v)\}_{1\leq j \leq p, 1 \leq k \leq d},t,h \in \mathbb{Z},(u,v)\in \cU \times \cV.$  
\end{definition}
\begin{condition} \label{cond_cfsm}
For $\{\bX_t(\cdot)\}_{t\in\mathbb{Z}}$ and $\{\bY_t(\cdot)\}_{t\in\mathbb{Z}},$
$\bbf^{X,Y}_\theta, \theta \in [-\pi,\pi]$ exists and the functional cross-spectral stability measure defined in (\ref{eq_M_cross}) is finite, i.e.
\begin{equation} \label{eq_M_cross}
\cM^{X,Y} = 2\pi \mathop{\text{ess sup}}\limits_{\theta \in [-\pi,\pi],\bPhi_1\in \mathbb{H}_0^p,\bPhi_2\in \mathbb{H}_0^d}\frac{\left|\langle\bPhi_1,\bbf_\theta^{X,Y}(\bPhi_2)\rangle\right|}{\sqrt{\langle\bPhi_1,\bSigma_0^{X}(\bPhi_1)\rangle}\sqrt{\langle\bPhi_2,\bSigma_0^{Y}(\bPhi_2)\rangle}}<\infty,
\end{equation}
where $\mathbb{H}_0^p = \{\bPhi \in \mathbb{H}^p:\langle\bPhi,\bSigma_0^X(\bPhi)\rangle \in (0, \infty)\}$ and $\mathbb{H}_0^d = \{\bPhi \in \mathbb{H}^d:\langle\bPhi,\bSigma_0^Y(\bPhi)\rangle \in (0, \infty)\}.$
\end{condition}

\begin{remark}
\begin{enumerate} [(a)]
    \item  If $\{\bX_t(\cdot)\}$ are independent of $\{\bY_t(\cdot)\},$ then $\mathcal{M}^{X,Y} = 0.$ Moreover, in the special case that $\{\bX_t(\cdot)\}$ and $\{\bY_t(\cdot)\}$ are identical, $\cM^{X,Y}$ degenerates to $\cM^X$ in (\ref{eq_M_X}). 

    \item Under the non-functional setting where $\bX_t \in {\mathbb R}^p$ and $\bY_t \in {\mathbb R}^d,$ \cite{basu2015a} introduced an upper bound condition for their proposed cross-spectral stability measure with $p=d,$ i.e. 
\begin{equation} \label{eq_basuM}
    \widetilde{\mathcal{M}}^{X,Y}=  \mathop{\text{ess sup}}\limits_{\theta \in [-\pi,\pi],\bnu\in \widetilde{\mathbb{R}}_0^d }\sqrt{\frac{\bnu^{\T} \{\bbf^{X,Y}_\theta\}^* \bbf^{X,Y}_\theta\bnu}{\bnu^{\T}\bnu}} <\infty,
\end{equation}
where $\widetilde{\mathbb{R}}_0^d = \{\bnu \in \mathbb{R}^d: \bnu^{\T}\bnu \in (0, \infty)\}$ and $*$ denotes the conjugate. This measure relates the cross-stability condition to the largest singular value of the cross-spectral density matrix $\bbf_\theta^{X,Y}.$ 
On the other hand, the non-functional analog of (\ref{eq_M_cross}) is equivalent to 
\begin{equation} \nonumber
    \mathop{\text{ess sup}}\limits_{\theta \in [-\pi,\pi],\bnu_1\in \widetilde{\mathbb{R}}_0^p,\bnu_2\in \widetilde{\mathbb{R}}_0^d}\frac{\left|\bnu_1^{\T} \bbf_\theta^{X,Y}\bnu_2\right|}{\sqrt{\bnu_1^{\T}\bnu_1}\sqrt{\bnu_2^{\T}\bnu_2}} <\infty,
\end{equation}
whose upper bound is $\widetilde{\mathcal{M}}^{X,Y}$ as justified in Lemma~\ref{lemma_compare} in Appendix~\ref{ap.lemma.secA3}. This demonstrates that, compared with (\ref{eq_basuM}), our proposed cross-stability measure 
corresponds to a milder condition. 

\item  For two truly infinite-dimensional functional objects, one limitation of the functional analog of $\widetilde{\cM}^{X,Y}$ is that it only controls the largest singular value of  $\bbf_\theta^{X,Y}.$ By contrast, our proposed $\cM^{X,Y}$ can more precisely characterize the effect of singular values of  $\bbf_\theta^{X,Y}$ relative to small eigenvalues of $\bSigma_0^X$ and $\bSigma_0^Y.$ Furthermore, it facilitates the development of finite sample theory for normalized versions of relevant estimated cross terms, where the normalization is formed by the corresponding eigenvalues in the denominator of $\cM^{X,Y}$. See Sections~\ref{sec_cross cov} and \ref{sec_fpca} for details.
\item We can generalize  (\ref{eq_M_cross}) to measure the serial and cross dependence structure between a mixture of multivariate functional and scalar time series. Specifically, consider $\{\bX_t(\cdot)\}_{t \in \mathbb{Z}} $ and $d$-dimensional vector time series $\{\bZ_t\}_{t \in \mathbb{Z}}$ with autocovariance matrices $\bSigma_h^Z$ for $h \in \mathbb Z.$ We can similarly define $\bbf_{\theta}^{X,Z}=\frac{1}{2\pi}\sum_{h \in \mathbb Z} \bSigma_h^{X,Z}\exp(-ih\theta)$ with $\bSigma_h^{X,Z}(\cdot)=\cov(\bX_t(\cdot), \bZ_{t+h}).$ 
According to (\ref{eq_M_cross}), the mixed cross-spectral stability measure of $\{\bX_t(\cdot)\}$ and $\{\bZ_t\}$ can be defined by
\begin{equation} \label{def_m_cfsm}
\mathcal{M}^{X,Z} = 2\pi \mathop{\text{ess sup}}\limits_{\theta \in [-\pi,\pi],\bPhi\in \mathbb{H}_0^p,\bnu\in \mathbb{R}_0^d}\frac{\left|\langle \bPhi, \bbf_\theta^{X,Z}\bnu\rangle\right|}{\sqrt{\langle\bPhi,\bSigma_0^{X}(\bPhi)\rangle}\sqrt{\bnu^{\T}\bSigma_0^{Z}\bnu}}
\end{equation}
and the non-functional stability measure of $\{\bZ_t\}$ reduces to

\begin{equation} \label{def_v_fsm}
\cM^Z = 2 \pi \cdot \underset{\theta\in [-\pi, \pi],\bnu \in \mathbb{R}_0^d}{\text{ess}\sup} \frac{\bnu^{\T} \bbf_{\theta}^Z\bnu}{\bnu^{\T}\bSigma_0^Z \bnu},
\end{equation}
where $\mathbb{R}_0^d = \{\bnu \in \mathbb{R}^d: \bnu^{\T}\bSigma_0^Z\bnu \in (0, \infty)\}.$
The proposed stability measures in (\ref{def_m_cfsm}) and (\ref{def_v_fsm}) play an essential role in the convergence analysis of the regularized estimates for 
model~(\ref{model_scalar}). 
See Section~\ref{sec_partial} for details.
\end{enumerate}
\end{remark}

For any $k_1$-dimensional subset $J$ of $\{1, \dots, p\}$ and $k_2$-dimensional subset $K$ of $\{1, \dots, d\},$ we can accordingly define the functional cross-stability measure of two subprocesses.
\begin{definition}
\label{def_cfsm}
Consider subprocesses 
$\left\{(X_{tj}(\cdot)):j \in J\right\}_{t \in \eZ}$ for $J \subseteq \{1,\dots,p\}$ with $|J| \leq k_1$ $(k_1=1,\dots,p)$ and 
$\left\{(Y_{tk}(\cdot)):k \in K\right\}_{t \in \eZ}$ for $K \subseteq \{1,\dots,d\}$ with $|K| \leq k_2$ $(k_2=1,\dots,d),$ their functional cross-spectral stability measure is defined by
\begin{equation} \label{def_sub_cfsm}
  \mathcal{M}_{k_1,k_2}^{X,Y} = 2\pi \underset{\|\bPhi_1\|_0 \leq k_1,\|\bPhi_2\|_0 \leq k_2}{\underset{\theta \in [-\pi,\pi],\bPhi_1\in \mathbb{H}_0^p,\bPhi_2\in \mathbb{H}_0^d}{\text{ess sup}}}\frac{\left|\langle\bPhi_1,\bbf_\theta^{X,Y}(\bPhi_2)\rangle\right|}{\sqrt{\langle\bPhi_1,\bSigma_0^{X}(\bPhi_1)\rangle}\sqrt{\langle\bPhi_2,\bSigma_0^{Y}(\bPhi_2)\rangle}}.
\end{equation}
\end{definition}
Under Condition~\ref{cond_cfsm}, it is easy to verify that,
$$
\mathcal{M}_{k_1,k_2}^{X,Y} \leq \mathcal{M}_{k_1^{'},k_2^{'}}^{X,Y} \leq \mathcal{M}^{X,Y} < \infty~~ \text{for }k_1 \leq k_1^{'} \text{ and } k_2 \leq k_2^{'}.
$$ 
According to (\ref{def_sub_fsm}), (\ref{def_m_cfsm}), (\ref{def_v_fsm}) and (\ref{def_sub_cfsm}), we can similarly define $\cM_{k_1,k_2}^{X,Z}$ and $\cM_{k_2}^Z$ for $k_1=1, \dots, p$ and $k_2=1, \dots, d,$ which will be used in our subsequent analysis.

\subsection{Sub-Gaussian functional linear process} 
\label{sec_subG}
Before presenting relevant non-asymptotic results beyond Gaussian functional time series, we introduce the definitions of sub-Gaussian process and multivariate functional linear process in this section.


Provided that our non-asymptotic analysis is based on the infinite-dimensional analog of Hanson--Wright inequality \cite[]{rudelson2013}
for sub-Gaussian random variables taking values within a Hilbert space, we first define sub-Gaussian process as follows.
\begin{definition}\label{def_subGaussian}
Let $X_t(\cdot)$ be a mean zero random variable in $\mathbb{H}$ and $\Sigma_0: \mathbb{H} \to \mathbb{H}$ be a covariance operator. Then $X_t(\cdot)$ is a sub-Gaussian process if there exists an $\alpha \geq 0$ such that for all $ x \in \mathbb{H},$
\begin{equation}
    \mathbb{E}\{e^{\langle x, X \rangle}\} \leq e^{\alpha^2\langle x,\Sigma_0(x)\rangle/2}.
\end{equation}
\end{definition}

The proof of Hanson--Wright inequality for serially dependent random functions relies on the fact that uncorrelated Gaussian random functions are also independent, which does not apply for non-Gaussian random functions. However, we show that, for a larger class of non-Gaussian functional time series, it is possible to develop finite sample theory for useful estimated terms in Sections~\ref{sec_cross cov} and \ref{sec_fpca}. We focus on multivariate functional linear processes with sub-Gaussian errors, namely sub-Gaussian functional linear processes:
\begin{equation} \label{eq_subG linear process}
    \bX_t(\cdot) = \sum\limits_{l=0}^{\infty}\bA_l(\bvarepsilon_{t-l}),~~ t \in \mathbb{Z},
\end{equation}
where $\bA_l=(A_{l,jk})_{p\times p}$ with each $A_{l,jk} \in {\mathbb S}$ and $\bvarepsilon_{t}(\cdot) = (\varepsilon_{t1}(\cdot),\dots,\varepsilon_{tp}(\cdot))^{\T} \in \cH^p.$ $\{\bvarepsilon_t(\cdot)\}_{t\in \mathbb Z}$ 
denotes a sequence of $p$-dimensional vector of random functions, whose components are independent sub-Gaussian processes satisfying Definition~\ref{def_subGaussian}.  It is worth noting that (\ref{eq_subG linear process}) not only extends the functional linear processes \cite[]{Bbosq1} to the multivariate setting but also can be seen as a generalization of $p$-dimensional linear processes \cite[]{Yao2019} to the functional domain.

Denote the polynomial $\mathcal{B}(z)(u,v) = \sum_{l = 0}^\infty \bA_l(u,v) z^l$ for $u,v \in {\cal U}.$
Under (\ref{eq_subG linear process}), we can derive the spectral density matrix function as
	\begin{equation}
	\label{f.X}
	\bbf^X_\theta(u,v) = \frac{1}{2\pi}  \int \int \mathcal{B}\left(e^{-i \theta}\right)(u,u') \bSigma^{\varepsilon}_0(u',v') \mathcal{B}\left(e^{-i \theta}\right)^*(v,v')du'dv'
		\end{equation}
	and the covariance matrix function as
	\begin{equation}
	 \label{cov.X}
	\bSigma_0^{X}(u,v) = \sum_{l = 0}^\infty \int \int \bA_l(u,u') \bSigma^{\varepsilon}_0(u',v')\bA_{l}^*(v,v')du'dv'.
	\end{equation}
Then we can express the functional stability measure ${\cal M}^X$ in (\ref{eq_M_X}) based on (\ref{f.X}) and (\ref{cov.X}).
The cross-spectral stability measure ${\cal M}^{X,Y}$ in (\ref{eq_M_cross}) or ${\cal M}^{X,Z}$ in (\ref{def_m_cfsm}) can be expressed in a similar fashion.


\begin{condition} \label{con_sub_coefficient}
The coefficient functions satisfy
$\sum_{l=0}^{\infty}\|\bA_l\|_{\infty} = O(1).$ 
\end{condition}

\begin{condition} \label{con_e}
(i) The marginal-covariance functions of $\{\bvarepsilon_t(\cdot)\},$ $\Sigma_{0,jj}^{\varepsilon}(u,v)$'s, are continuous on $\cU^2$ and uniformly bounded over $j \in \{1,\dots,p\};$ 
(ii) $\omega_0^{\varepsilon} = \max_{j}\int_\cU \Sigma_{0,jj}^{\varepsilon}(u,u)du =O(1) .$
\end{condition}

Condition~\ref{con_sub_coefficient} ensures functional analog of standard condition of elementwise absolute summability of
moving average coefficients for multivariate linear processes \cite[]{hamilton1994} under Hilbert--Schmidt norm.
It also guarantees the stationarity of $\{\bX_t(\cdot)\}$ and, furthermore together with  
Condition~\ref{con_e}, implies that  $\omega_0^{X} = \max_{j}\int_\cU \Sigma_{0,jj}^{X}(u,u)du = O(1),$ both of which are essential in our subsequent analysis. See Lemma~\ref{lm_X_traceclass} in Appendix~\ref{ap.lemma.secA3} for details. In general, we can relax Conditions~\ref{con_sub_coefficient} and \ref{con_e} by allowing $\sum_{l=0}^{\infty}\|\bA_l\|_{\infty}$ and $\omega_0^{\varepsilon}$ to grow at very slow rates as $p$ increases, then our subsequent non-asymptotic bounds will depend on $\omega_0^{X},$ or, more precisely, these two terms, which complicate the presentation of theoretical results.

\subsection{Concentration bounds on sample (cross-)(auto)covariance matrix function}
\label{sec_cross cov}
We construct estimated (auto)covariance of $\{\bX_t(\cdot)\}_{t=1}^n$  by
\begin{equation} \nonumber
    \widehat{\bSigma}_h^X(u,v) = \frac{1}{n-h}\sum_{t = 1}^{n-h}\bX_t(u)\bX_{t+h}(v)^{\T},~~h = 0,1,\dots, ~(u,v) \in \mathcal{U}^2,
\end{equation}
and estimated cross-(auto)covariance matrix functions between $\{\bX_t(\cdot)\}$ and $\{\bY_t(\cdot)\}$ by
\begin{equation} \nonumber
    \widehat{\bSigma}_h^{X,Y}(u,v) = \frac{1}{n-h}\sum_{t = 1}^{n-h}\bX_t(u)\bY_{t+h}(v)^{\T},~~ h = 0,1,\dots, ~(u,v) \in \mathcal{U}\times \cV.
\end{equation}
\begin{theorem}\label{th_XY}
Suppose that Conditions~\ref{con_stability}--\ref{con_e} hold for sub-Gaussian functional linear processes, $\{\bX_t(\cdot)\},$ $\{\bY_t(\cdot)\}$ and $h$ is fixed. Then for any given vectors $\bPhi_1\in \mathbb{H}_0^p$ and $\bPhi_2\in \mathbb{H}_0^d $ with $\|\bPhi_1\|_0 \leq k_1, \|\bPhi_2\|_0 \leq k_2$ $(k_1 =1,\dots,p, k_2 =1,\dots,d),$ there exists some constants $c, c_1, c_2 > 0$ such that for any $\eta > 0,$
\begin{equation} \label{eq_th_X}
     P\left\{\left|\frac{\langle\bPhi_1,(\widehat{\bSigma}_{0}^X - \bSigma_{0}^X)(\bPhi_1)\rangle}{\langle\bPhi_1,\bSigma_{0}^X(\bPhi_1)\rangle}\ \right| > \cM_{k_1}^X\eta\right\} \leq 2 \exp \left\{ -cn \min \left(\eta^2,\eta\right)\right\},
\end{equation}
and
\begin{equation} \label{eq_th_XY}
\begin{split}
    P\left\{\left|\frac{\langle\bPhi_1,(\widehat{\bSigma}_h^{X,Y}-\bSigma_h^{X,Y})(\bPhi_2)\rangle}{\langle\bPhi_1,\bSigma_0^X(\bPhi_1)\rangle +\langle\bPhi_2,\bSigma_0^Y(\bPhi_2)\rangle}\right| >  \left(\mathcal{M}_{k_1}^{X} +\mathcal{M}_{k_2}^{Y} +\mathcal{M}_{k_1,k_2}^{X,Y}\right)\eta \right\} \\
\leq c_1\exp\{-c_2n\min(\eta^2,\eta)\}.
\end{split}
\end{equation}
\end{theorem}
\begin{remark}
(\ref{eq_th_X}) extends the concentration inequality for normalized quadratic form of $\widehat \bSigma_0^X$ 
in Theorem~1 of \cite{guo2019} under the Gaussianity assumption to accommodate a larger class of sub-Gaussian functional linear processes and serves as a starting point to establish further useful non-asymptotic results, e.g. those listed in Theorems~1--4 and Proposition~1 of \cite{guo2019}, so we present some results used in our subsequent analysis in Appendix~\ref{ap.subgaussian}.
The concentration inequality in (\ref{eq_th_XY}) illustrates that the tail for normalized bilinear form of $\widehat{\bSigma}_h^{X,Y}-\bSigma_h^{X,Y}$ behaves in a sub-Gaussian or sub-exponential way depending on which term in the tail bound is dominant. It is also crucial in deriving subsequent concentration results, e.g. with suitable choices of $\bPhi_1$ and $\bPhi_2,$ it facilitates the elementwise concentration bounds on $\widehat\bSigma_h^{X,Y}$ in the following theorem.

\end{remark}

\begin{theorem} \label{th_Dev}
Suppose that conditions in Theorem~\ref{th_XY} hold. Then there exists some constants $c_1,c_3>0$ such that for any $\eta>0$ and each $j = 1,\dots,p,$ $k = 1,\dots,d,$
\begin{equation}\label{eq_Dev_j}
P\left\{\| \widehat{\Sigma}_{h,jk}^{X,Y}-\Sigma_{h,jk}^{X,Y}\|_\cS > (\omega_0^X + \omega_0^Y)\cM_{X,Y}\eta\right\} \leq c_1\exp\left\{-c_3n\min(\eta^2,\eta)\right\},
\end{equation}
where $\cM_{X,Y} = \mathcal{M}_1^{X} +\mathcal{M}_1^{Y}+\mathcal{M}_{1,1}^{X,Y}, \omega_0^{X} = \max_{ j }\int_\cU \Sigma_{0,jj}^{X}(u,u)du $ and $\omega_0^{Y} = \max_{ k }\int_\cU \Sigma_{0,kk}^{Y}(u,u)du.$ In particular,  
there exists some constant $c_4>0$ such that, for sample size $n \gtrsim \log (pd),$ with probability greater than $1-c_1(pd)^{-c_4},$ the estimate $\widehat{\bSigma}_h^{X,Y}$ satisfies the bound 
\begin{equation}\label{eq_Dev_max}
\| \widehat{\bSigma}_h^{X,Y} -\bSigma_h^{X,Y}\|_{\max} \lesssim \cM_{X,Y}\sqrt{\frac{\log (pd)}{n}}.
\end{equation}
\end{theorem}

\begin{remark}
In the deviation bounds established above, the effects of dependence are commonly captured by the sum of marginal-spectral and cross-spectral stability measures, $\cM_{X,Y}=\mathcal{M}_1^{X} +\mathcal{M}_1^{Y}+\mathcal{M}_{1,1}^{X,Y},$ with larger values yielding a slower elementwise $\ell_{\infty}$ rate in (\ref{eq_Dev_max}).
Under a mixed-process scenario consisting of $\{\bX_t(\cdot)\}$ and $d$-dimensional time series $\{\bZ_t\}$ belonging to multivariate linear processes with sub-Gaussian errors \cite[]{sun2018}, namely sub-Gaussian linear processes, it is easy to extend (\ref{eq_Dev_max}) as
\begin{equation} \label{eq_Dev_max_XZ}
    \max_{1 \leq j\leq p,1 \leq k\leq d}\|\widehat{\Sigma}_{h,jk}^{X,Z} -\Sigma_{h,jk}^{X,Z}\|\lesssim \cM_{X,Z}\sqrt{\frac{\log (pd)}{n}},
\end{equation}
where $\cM_{X,Z} = \cM_1^X+\cM_1^Z+\cM_{1,1}^{X,Z}.$
(\ref{eq_Dev_max_XZ}) can be justified in the  proof of Proposition~\ref{pr_FPCscores_XZ} in Appendix~\ref{ap.pro.secA2}.
\end{remark}

\subsection{Rates in elementwise $\ell_{\infty}$ norm under a FPCA framework} 
\label{sec_fpca}
For each $j=1, \dots, p,$ suppose that $X_{1j}(\cdot), \dots, X_{nj}(\cdot)$ are $n$ serially dependent observations of $X_j(\cdot).$ The Karhunen-Lo\`eve theorem \cite[]{Bbosq1} serving as the theoretical basis of FPCA allows us to represent each functional observation in the form of 
$X_{tj}(\cdot)=\sum_{l=1}^{\infty}\zeta_{tjl}\psi_{jl}(\cdot).$ Here the coefficients $\zeta_{tjl} = \langle X_{tj}, \psi_{jl} \rangle,$ namely FPC scores, are uncorrelated random variables with mean zero and $\cov(\zeta_{tjl}, \zeta_{tjl'}) = \omega_{jl}^X I(l=l').$ In this formulation, $\{(\omega_{jl}^X, \psi_{jl})\}_{l=1}^{\infty}$ are eigenpairs satisfying $\langle\Sigma_{0,jj}^X(u,\cdot), \psi_{jl}(\cdot)\rangle=\omega_{jl}^X\psi_{jl}(u).$ Similarly, for each $k = 1,\dots, d,$ we represent $Y_{tk}(\cdot) = \sum_{m=1}^\infty\xi_{tkm}\phi_{km}(\cdot)$ with eigenpairs $\{(\omega_{km}^Y, \phi_{km})\}_{m=1}^{\infty}.$


To estimate relevant terms under a FPCA framework, for each $j,$ we perform an eigenanalysis on $\widehat{\Sigma}_{0,jj}^X(u,v) = n^{-1}\sum_{t = 1}^n X_{tj}(u)X_{tj}(v),$ i.e. $\langle \widehat\Sigma_{0,jj}^X(u,\cdot), \widehat\psi_{jl}(\cdot)\rangle  =\widehat \omega_{jl}^X\widehat\psi_{jl}(u),$ where $\{(\widehat\omega_{jl}^X, \widehat \psi_{jl})\}_{l=1}^{\infty}$ denote the estimated eigenpairs. The corresponding estimated FPC scores are given by $\widehat\zeta_{tjl} = \langle X_{tj}, \widehat\psi_{jl}\rangle.$ Furthermore, relevant estimated terms for $\{Y_{tk}(\cdot)\},$ i.e. $\widehat\omega_{km}^Y, \widehat\phi_{km}(\cdot),\\ \widehat \xi_{tkm},$ can be obtained in the same manner. 


Before presenting relevant deviation bounds in elementwise $\ell_{\infty}$ norm, which are essential under high-dimensional regime, $(\log p \vee \log d)/n \rightarrow 0,$ we impose the following lower bound condition on the eigengaps. 

\begin{condition}\label{con_eigen}
For each $j=1, \dots, p$ and $k=1, \dots,d,$ $\omega_{j1}^X > \omega_{j2}^X > \cdots >0$  and $\omega_{k1}^Y > \omega_{k2}^Y > \cdots >0.$ There exist some positive constants $c_0$ and $\alpha_1, \alpha_2 > 1$ such that $\omega_{jl}^X - \omega_{j(l+1)}^X \geq c_0 l^{-\alpha_1-1}$ for $l = 1,\dots,\infty$ and $\omega_{km}^Y - \omega_{k(m+1)}^Y \geq c_0 m^{-\alpha_2-1}$ for $m = 1,\dots,\infty.$
\end{condition}

Condition~\ref{con_eigen} implies the lower bounds on eigenvalues, i.e. $\omega_{jl}^X \geq c_0\alpha_1^{-1}l^{-\alpha_1}$ and $\omega_{km}^Y \geq c_0\alpha_2^{-1}m^{-\alpha_2}.$ 
See also \cite{hall2007} and \cite{kong2016} for similar conditions.

In practice, the infinite series in the Karhunen-Lo\`eve expansions of $X_{tj}(\cdot)$ and $Y_{tm}(\cdot)$ are truncated at $M_1$ and $M_2,$ chosen data-adaptively, which transforms the infinite-dimensional learning task into the modelling of multivariate time series. Given sub-Gaussian functional linear process $\{\bX_t(\cdot)\},$ to aid convergence analysis under high-dimensional scaling, we establish elementwise concentration inequalities and, furthermore, elementwise $\ell_{\infty}$ error bounds 
on relevant estimated terms, i.e. estimated eigenpairs and sample (auto)covariance between estimated FPC scores. These results are of the same forms as those under the Gaussianity assumption \citep{guo2019}, so we only present them in Lemmas~\ref{lm_guo2019_eigen} and \ref{lm_guo2019_score} in Appendix~\ref{ap.subgaussian}. 

In the following, we focus on sample cross-(auto)covariance between estimated FPC scores, $\widehat{\sigma}_{h,jklm}^{X,Y} = (n-h)^{-1}\sum_{t=1}^{n-h}\widehat{\zeta}_{tjl}\widehat{\xi}_{(t+h)km},$ and establish a normalized deviation bound in elementwise $\ell_{\infty}$ norm on how $\widehat{\sigma}_{h,jklm}^{X,Y}$ concentrates around $\sigma_{h,jklm}^{X,Y} = \cov(\zeta_{tjl},\xi_{(t+h)km}).$ 


\begin{theorem} \label{th_FPCscores_XY}
	Suppose that Conditions~\ref{con_stability}--\ref{con_eigen} hold for sub-Gaussian functional linear processes, $\{\bX_t(\cdot)\},$ $\{\bY_t(\cdot)\},$ and $h$ is fixed.
	Let $M_1$ and $M_2$ be positive integers possibly depending on $(n,p,d).$ 
	 If $n \gtrsim \log(pdM_1M_2)(M_1^{4\alpha_1+2} \vee M_2^{4\alpha_2+2})\cM_{X,Y}^2,$
	 then  
there exist some positive constants $c_5$ and $c_6$ such that, with probability greater than 
$ 1- c_5(pdM_1M_2)^{-c_6},$ the estimates $\{\widehat \sigma_{h,jklm}^{X,Y}\}$ satisfy
\begin{equation}
\label{eq_FPCscores_XY}
\underset{\underset{1 \le l \le M_1, 1 \le m \le M_2}{1\le j\le p, 1\le k\le d}}{\max}~ \frac{\left|\widehat \sigma_{h,jklm}^{X,Y} - \sigma_{h,jklm}^{X,Y}\right|}{(l^{\alpha_1+1} \vee m^{\alpha_2+1}){\sqrt{\omega_{jl}^X \omega_{km}^Y}}}\lesssim \mathcal{M}_{X,Y}\sqrt{\frac{\log(pdM_1M_2)}{n}}.
\end{equation}
\end{theorem}


In the special case that $\{\bX_t(\cdot)\}$ and $\{\bY_t(\cdot)\}$ are identical, (\ref{eq_FPCscores_XY}) degenerates to the deviation bound on $\widehat\sigma_{h,jklm}^{X}$ under the Gaussianity assumption \cite[]{guo2019}. We next consider a mixed process scenario consisting of $\{\bX_t(\cdot)\}$ and $\{\bZ_t\}$ and establish a normalized deviation bound in elementwise $\ell_{\infty}$ norm on sample cross-(auto)covariance between estimated FPC scores of $\{X_{tj}(\cdot)\}$ and $Z_{(t+h)k}.$ Define $\widehat\varrho_{h,jkl}^{X,Z} =(n-h)^{-1}\sum_{t=1}^{n-h}\widehat{\zeta}_{tjl}Z_{(t+h)k} $ and $\varrho_{h,jkl}^{X,Z} =\cov({\zeta}_{tjl},Z_{(t+h)k}).$ We are ready to extend (\ref{eq_FPCscores_XY}) to the following mixed-process scenario.


\begin{proposition} \label{pr_FPCscores_XZ}
Suppose that Conditions~\ref{con_stability}--\ref{con_eigen} hold for sub-Gaussian functional linear process $\{\bX_t(\cdot)\},$ $\{\bZ_t\}$ follows sub-Gaussian linear process and $h$ is fixed.
Let $M_1$ be a positive integer possibly depending on $(n,p,d).$ 
If sample size $ n\gtrsim \log(pdM_1)M_1^{3\alpha_1+2}\cM_{X,Z}^2,$ then
there exist some constants $c_7, c_8>0$ such that, with probability greater than 
$ 1- c_7(pdM_1)^{-c_8},$ the estimates $\{\widehat \varrho_{h,jkl}^{X,Z}\}$ satisfy
\begin{equation}
\label{eq_FPCscores_XY_partial}
\underset{\underset{1 \le l \le M_1}{1\le j\le p, 1\le k\le d}}{\max} ~\frac{\left|\widehat \varrho_{h,jkl}^{X,Z} - \varrho_{h,jkl}^{X,Z}\right|}{l^{\alpha_1+1}\sqrt{\omega_{jl}^X}}\lesssim \mathcal{M}_{X,Z}\sqrt{\frac{\log(pdM_1)}{n}}.
\end{equation}
\end{proposition}

We next consider $\{\epsilon_t(\cdot)\}_{t=1}^n$, defined on $\cV,$ which can be seen as the error term in model~(\ref{model_fully}) being independent of $\{\bX_t(\cdot)\}.$ Define $\Sigma_{h,j}^{X,\epsilon}(u,v) = \cov\{X_{tj}(u),\epsilon_{t+h}(v)\}$ and its estimate $\widehat\Sigma_{h,j}^{X,\epsilon}(u,v) = (n-h)^{-1}\sum_{t=1}^{n-h}X_{tj}(u)\epsilon_{t+h}(v).$ To provide theoretical analysis of the estimates for model~(\ref{model_fully}), the FPCA-based representation in Appendix~\ref{sec_sm_fully_matrix} suggests to investigate the consistency properties of the estimated cross terms, i.e. $\widehat{\sigma}_{h,jlm}^{X,\epsilon} = \langle \widehat \psi_{jl},\langle \widehat{\Sigma}_{h,j}^{X,\epsilon}, \widehat\phi_{m}\rrangle$ or 
$\widehat{\sigma}_{h,jlm}^{X,Y} = (n-h)^{-1}\sum_{t=1}^{n-h}\widehat{\zeta}_{tjl}\widehat{\xi}_{(t+h)m}=\langle \widehat \psi_{jl},\langle \widehat{\Sigma}_{h,j}^{X,Y}, \widehat\phi_{m}\rrangle.$ As $\{\bX_{t-h}(\cdot): h=0, \dots, L\}$ and $\{\epsilon_t(\cdot)\}$ are independent and can together determine the response $\{Y_t(\cdot)\}$ via (\ref{model_fully}), 
it is more sensible to study the former term, i.e. how $\widehat \sigma_{h,jlm}^{X, \epsilon}$ deviates from $\sigma_{h,jlm}^{X, \epsilon} = 0$ in the following proposition.


\begin{proposition}
 \label{pr_FPCscores_XE}
	Suppose that Conditions~\ref{con_stability}--\ref{con_eigen} hold for sub-Gaussian functional linear processes $\{\bX_t(\cdot)\},$ $\{\epsilon_t(\cdot)\}$ and $h$ is fixed.
	Let $M_1, M_2$ be positive integers possibly depending on $(n,p).$ 
	If $n \gtrsim \log(pM_1M_2)(M_1^{4\alpha_1+4} \vee M_2^{4\alpha_2+4})(\cM_1^X+\cM^Y)^2,$ then 
there exist some constants $c_9, c_{10}>0$ such that, with probability greater than 
$ 1- c_9(pM_1M_2)^{-c_{10}},$ the estimates $\{\widehat \sigma_{h,jlm}^{X,\epsilon}\}$ satisfy
\begin{equation}
\label{eq_FPCscores_XE}
\underset{\underset{1 \le l \le M_1, 1 \le m \le M_2}{1\le j\le p}}{\max} \frac{\left|\widehat \sigma_{h,jlm}^{X,\epsilon} \right|}{(l^{\alpha_1} \vee m^{\alpha_2}){\sqrt{\omega_{jl}^X \omega_{m}^Y}}}\lesssim (\cM_1^X + \cM^\epsilon)\sqrt{\frac{\log(pM_1M_2)}{n}}.
\end{equation}
\end{proposition}

Finally, we consider a mixed-process scenario in model~(\ref{model_scalar}), where $\{\epsilon_t\}_{t=1}^n$ are scalar errors, independent of both $\{\bX_t(\cdot)\}$ and $\{\bZ_t\}.$ In addition to Proposition~\ref{pr_FPCscores_XZ} above, the following proposition demonstrates how $\widehat\varrho_{h,jl}^{X,\epsilon} =(n-h)^{-1}\sum_{t=1}^{n-h}\widehat{\zeta}_{tjl}\epsilon_{t+h} $ converges to $\varrho_{h,jl}^{X,\epsilon} =\cov({\zeta}_{tjl},\epsilon_{t+h})=0.$ 

\begin{proposition}
 \label{pr_FPCscores_XE_partial}
 Suppose that Conditions~\ref{con_stability}--\ref{con_eigen} hold for sub-Gaussian functional linear process $\{\bX_t(\cdot)\},$ $\{\epsilon_t\}$ is sub-Gaussian linear process and $h$ is fixed.
	Let $M_1$ be positive integer possibly depending on $(n,p).$ If $n \gtrsim\log(pM_1)M_1^{3\alpha_1+2} (\cM_1^X)^2$, then there exist some constants $c_{11}, c_{12}>0$ such that, with probability greater than 
$ 1- c_{11}(pM_1)^{-c_{12}},$ the estimates $\{\widehat \varrho_{h,jl}^{X,\epsilon}\}$ satisfy
\begin{equation}
\label{eq_FPCscores_XE_partial}
\underset{1\le j\le p,1 \le l \le M_1}{\max} \frac{\left|\widehat \varrho_{h,jl}^{X,\epsilon} \right|}{\sqrt{\omega_{jl}^X}}\lesssim (\cM_1^X + \cM^\epsilon)\sqrt{\frac{\log(pM_1)}{n}}.
\end{equation}
\end{proposition}

\begin{remark}
Benefiting from the independence assumption between $\{\bX_t(\cdot)\}$ and $\{\epsilon_t(\cdot)\},$ 
Proposition~\ref{pr_FPCscores_XE} leads to a faster rate of convergence in (\ref{eq_FPCscores_XE}) compared with (\ref{eq_FPCscores_XY}) with $d=1.$ Proposition~\ref{pr_FPCscores_XE} also plays a crucial rule in the proof of Proposition~\ref{pr_fof_max.error} to demonstrate that, with high probability, model~(\ref{model_fully}) satisfies  
the routinely used deviation condition. 
Analogously, taking an advantage of the independence assumption between $\{\bX_t(\cdot)\}$ and $\{\epsilon_t\},$  Proposition~\ref{pr_FPCscores_XE_partial} results in a faster rate in (\ref{eq_FPCscores_XE_partial}) than that in (\ref{eq_FPCscores_XY_partial}) with $d=1.$ In the proof of Proposition~\ref{pr_partial_max.error_X}, we will apply Proposition~\ref{pr_FPCscores_XE_partial} to verify that, with high probability, model~(\ref{model_scalar}) satisfies the corresponding deviation condition.
\end{remark}

\section{High-dimensional functional linear lagged regression} 
\label{sec_fully}
In this section, we first develop a three-step procedure to estimate sparse functional coefficients in model~(\ref{model_fully}) and then apply our derived finite sample results in Section~\ref{sec_fpca} to investigate the convergence properties of the estimates under high-dimensional scaling.


\subsection{Estimation procedure}
\label{sec.m1.est}
Consider functional linear lagged regression model in (\ref{model_fully}),
where $\{\beta_{hj} \in \mathbb{S}: h = 0, \dots, L, j = 1,\dots,p \}$ are unknown functional coefficients and $\{\epsilon_t(\cdot)\}_{t=1}^n$ are mean-zero errors from sub-Gaussian functional linear process, independent of $\{\bX_t(\cdot)\}_{t=1}^n$ from sub-Gaussian functional linear process.
Given observed data $\{Y_t, \bX_t\}_{t=1}^n,$ our goal is to estimate a vector of functional coefficients, $\bbeta = (\beta_{01},\dots,\beta_{0p},\dots, \beta_{L1},\dots,\beta_{Lp})^{\T}$ with each $\beta_{hj}\in \mathbb{S}.$ To assure a feasible solution under a high-dimensional regime, we impose a sparsity assumption on $\bbeta.$ To be specific, we assume that $\bbeta$ is functional $s$-sparse with support set $S = \big\{ (h,j) \in \{0, \dots, L\}\times \{1,\dots,p\}:\|\beta_{hj}\|_\mathcal{S}\neq 0\big\}$ and its cardinality $|S| = s,$ much smaller than the dimensionality, $p(L+1).$

Due to the infinite dimensional nature of functional data, we approximate each $X_{tj}(\cdot)$ and $Y_{t}(\cdot)$ under the Karhunen-Lo\`eve expansion truncated at $q_{1j}$ and $q_2,$ respectively, i.e. 
$$X_{tj}(\cdot)\!\approx\!\sum_{l=1}^{q_{1j}}\zeta_{tjl}\psi_{jl}(\cdot) \!=\! \bzeta_{tj}^{\T}\bpsi_j(\cdot), ~~~ Y_t(\cdot) \!\approx\!\sum_{m=1}^{q_2}\xi_{tm}\phi_m(\cdot) \!=\! \bxi_t^{\T}\bphi(\cdot),$$ where $\bzeta_{tj} \!= \!(\zeta_{tj1},\dots,\zeta_{tjq_{1j}})^{\T},$ $\bpsi_j(\cdot) \!=\! (\psi_{j1}(\cdot),\dots,\psi_{jq_{1j}}(\cdot))^{\T},$ $\bxi_t \!= \!(\xi_{t1},\dots,\xi_{tq_2})^{\T}$ and $\bphi(\cdot)\! = \!(\phi_1(\cdot),\dots,\phi_{q_2}(\cdot))^{\T}.$ The truncation levels $q_{1j}$ and $q_{2}$ are carefully chosen so as to provide reasonable approximations to each $X_{tj}(\cdot)$ and $Y_t(\cdot).$ See \cite{kong2016} for the selection of the truncated dimension in practice.


According to Appendix~\ref{sec_sm_fully_matrix}, we can represent model~(\ref{model_fully}) in the following matrix form
\begin{equation} \label{model_fully_matrix}
    \bU = \sum_{h=0}^L\sum_{j=1}^p\bV_{hj}\bPsi_{hj} + \bR + \bE, 
\end{equation}
where $\bPsi_{hj}=\int_{\cV}\int_{\cU}\bpsi_j(u)\beta_{hj}(u,v)\bphi(v)^{\T} dudv\in\mathbb{R}^{q_{1j}\times q_2},$ $\bU\in \mathbb{R}^{(n-L)\times q_2}$ with its row vectors given 
by $\bxi_{L+1},\dots,\bxi_{n}$ and $\bV_{hj}\in \mathbb{R}^{(n-L)\times q_{1j}}$ with its row vectors given by $\bzeta_{(L+1-h)j},\dots,\bzeta_{(n-h)j}.$  Note $\bR$ and $\bE$ are $(n-L)\times q_2$ matrices whose row vectors are formed by truncation errors $\{\br_t\in\mathbb{R}^{q_2}: t=L+1, \dots, n\}$ and random errors $\{\beps_t\in\mathbb{R}^{q_2}: t=L+1, \dots, n\}$ respectively.


We develop the following three-step estimation procedure.

First, we perform FPCA on $\{X_{tj}(\cdot)\}_{t=1}^n$ for each $j=1,\dots,p$ and $\{Y_t(\cdot)\}_{t=1}^n,$ thus obtaining estimated FPC scores and eigenfunctions, i.e. $\widehat\zeta_{tjl}, \widehat\psi_{jl}(\cdot)$ for $l \geq 1$ and $\widehat\xi_{tm}, \widehat\phi_{tm}(\cdot)$ for $m \geq 1,$ respectively.

Second, it is worth noting that the problem of recovering functional sparsity structure in $\bbeta$ is equivalent to estimating the block sparsity pattern in $\{\bPsi_{hj}: h = 0, \dots, L, j = 1, \dots,p\}.$ Specifically, if $\beta_{hj}(\cdot,\cdot)$ is zero, all entries in $\bPsi_{hj}$ will be zero. This motivates us to incorporate a standardized group lasso penalty \cite[]{simon2012} by minimizing the following penalized regression criterion over $\{\bPsi_{hj}: h = 0, \dots, L, j = 1, \dots,p\}$:
\begin{equation}
\label{lasso}
\frac{1}{2}\|\widehat{\bU}-\sum_{h=0}^L\sum_{j=1}^p\widehat{\bV}_{hj}\bPsi_{hj}\|_\tF^2 + \lambda_n\sum_{h=0}^L\sum_{j=1}^p\|\widehat{\bV}_{hj}\bPsi_{hj}\|_\tF,
\end{equation}
where $\widehat{\bU}$ and $\widehat{\bV}_{hj}$ are the estimates of $\bU$ and $\bV_{hj},$ respectively,  and $\lambda_n$ is a non-negative regularization parameter. Let $\{\widehat\bPsi_{hj}\}$ be the minimizer of (\ref{lasso}). 

Finally, we estimate functional coefficients by  $$\widehat{\beta}_{hj}(u,v) = \widehat{\bpsi}_j(u)^{\T}\widehat{\bPsi}_{hj}\widehat{\bphi}(v), ~~(u,v) \in \cU \times \cV, h = 0,\dots, L, j = 1, \dots, p.$$

\subsection{Theoretical properties}
\label{sec.thm.m1}
We begin with some notation that will be used in this section. For a block matrix $\bB = (\bB_{jk})_{1 \leq j \leq p_1, 1 \leq k \leq p_2} \in \eR^{p_1q_1 \times p_2 q_2}$ with the $(j,k)$-th block $\bB_{jk} \in {\eR}^{q_1 \times q_2},$  we define its $(q_1,q_2)$-block versions of elementwise $\ell_{\infty}$ and matrix $\ell_1$ norms by $\|\bB\|_{\max}^{(q_1, q_2)} = \max_{j,k} \|\bB_{jk}\|_{\tF}$ and
$
\|\bB\|_{1}^{(q_1,q_2)} = {\max}_{k}\sum_{j} \|\bB_{jk}\|_{\tF},
$
respectively.
To simplify notation, we will assume the same $q_{1j}$ across $j=1,\dots, p,$ but our theoretical results extend naturally to the more general setting where $q_{1j}$'s are different.

Let $\widehat{\bZ} = (\widehat{\bV}_{01},\dots,\widehat{\bV}_{0p}, \dots, \widehat{\bV}_{L1},\dots,\widehat{\bV}_{Lp})\in \mathbb{R}^{(n-L)\times (L+1)pq_1},$ $\bPsi = (\bPsi_{01}^{\T},\dots,\bPsi_{0p}^{\T},\dots, \bPsi_{L1}^{\T},\\ \dots,\bPsi_{Lp}^{\T})^{\T}\in\mathbb{R}^{(L+1)pq_1\times q_2}$ and $\widehat{\bD} =  \text{diag}(\widehat{\bD}_{01},\dots,\widehat{\bD}_{0p}, \dots, \widehat{\bD}_{L1},\dots,\widehat{\bD}_{Lp})\in \mathbb{R}^{(L+1)pq_1\times (L+1)pq_1}$ with $\widehat{\bD}_{hj} = \{(n-L)^{-1}\widehat{\bV}_{hj}^{\T}\widehat{\bV}_{hj}\}^{1/2}\in\mathbb{R}^{q_1\times q_1}$ for $h = 0, \dots, L$ and $j=1,\dots,p.$ Then minimizing (\ref{lasso}) over $\{\bPsi_{hj}\}$ is equivalent to the following optimization task:
\begin{equation}
\label{target}
 \widehat \bB = \underset{\bB \in \mathbb{R}^{(L+1)pq_1\times q_2}}{ \text{arg min}} \left\{ \frac{1}{2(n-L)}\|\widehat{\bU}-\widehat{\bZ}\widehat{\bD}^{-1}\bB\|_\tF^2 + \lambda_n\|\bB\|^{(q_1,q_2)}_1 \right\}.
\end{equation}
Then we have $\widehat{\bPsi} = \widehat{\bD}^{-1}\widehat{\bB}$ with its $\{(h+1)j\}$-th row block given by $\widehat{\bPsi}_{hj}.$ 


Before our convergence analysis, we present the following regularity conditions.

\begin{condition} \label{con_fof_beta_a}
For each $(h,j)\in S,$ $\beta_{hj}(u,v) = \sum_{l,m=1}^\infty a_{hjlm}\psi_{jl}(u)\phi_{m}(v)$ and there exist some positive constants $\kappa > (\alpha_1 \vee \alpha_2) /2 + 1$ and $\mu_{hj}$ such that $|a_{hjlm}| \leq \mu_{hj}(l+m)^{-\kappa-1/2}$ for $l,m \geq 1.$
\end{condition}

We expand each non-zero functional coefficient $\beta_{hj}(u, v)$ using principal component functions $\{\psi_{jl}(u)\}_{l \geq 1}$ and $\{\phi_m(v)\}_{m \geq 1},$ which respectively provide the most rapidly convergent representation of $\{X_{tj}(u)\}$ and $\{Y_t(v)\}$ in the $L_2$ sense. Such condition prevents the coefficients $\{a_{hjlm}\}_{l,m \geq 1}$ from decreasing too slowly with parameter $\kappa$ controlling the level of smoothness in non-zero components of $\{\beta_{hj}(\cdot,\cdot)\}.$ 
See similar smoothness conditions in functional linear regression literature \cite[]{hall2007,kong2016}.

\begin{condition} \label{con_fof_eigen_min}
Denote the covariance matrix function by
\begin{equation} \nonumber   \widetilde \bSigma^X = \left( 
    \begin{matrix}
    \bSigma_0^X & \bSigma_1^{X} & \cdots & \bSigma_L^X\\
    \bSigma_1^{X} & \bSigma_0^X & \cdots & \bSigma_{L-1}^X\\
    \vdots & \vdots & \ddots & \vdots\\
    \bSigma_L^X & \bSigma_{L-1}^X & \cdots & \bSigma_0^X
    \end{matrix}
    \right) 
\end{equation}
and the diagonal matrix function by $\widetilde \bD_0^X = \bI_{L+1}\otimes\text{diag}(\Sigma^X_{0,11},\dots,\Sigma^X_{0,pp}).$ The infimum $\underline{\mu}$ of the functional Rayleigh quotient of $\widetilde \bSigma^X$ relative to $\widetilde \bD_0^X$ is bounded below by zero, i.e.
\begin{equation} \nonumber
\begin{split}
\underline{\mu} = \inf\limits_{\bPhi \in \Bar{\mathbb{H}}_0^{(L+1)p}} \frac{\langle\bPhi,\widetilde \bSigma^X(\bPhi)\rangle}{\langle\bPhi,\widetilde \bD_0^X(\bPhi)\rangle} > 0,
\end{split}
\end{equation}
where $\bPhi \in \Bar{\mathbb{H}}_0^{(L+1)p} = \{\bPhi \in \mathbb{H}^{(L+1)p}: \langle\bPhi,\widetilde \bD_0^X(\bPhi)\rangle \in (0,\infty)\}.$
\end{condition}

Condition~\ref{con_fof_eigen_min} can be interpreted as requiring the minimum eigenvalue of the correlation matrix function for $(\bX_{t-L}^{\T},\dots, \bX_{t}^{\T})^{\T}$ to be bounded below by zero. 
See also a similar condition in \cite{guo2019}.

Before presenting the consistency analysis of $\widehat\bbeta$ in Theorem~\ref{th_beta}, we show that the functional analogs of the restricted eigenvalue (RE) condition and the deviation condition in the lasso literature \cite[]{loh2012} are satisfied with high probability 
in Proposition~\ref{pr_fof_RE} below and Propositions~\ref{pr_fof_eigen}--\ref{pr_fof_max.error} in Appendix~\ref{ap.prop}, respectively.



\begin{proposition} \label{pr_fof_RE}
Suppose Conditions~\ref{con_stability}--\ref{con_eigen} and \ref{con_fof_eigen_min} hold. Then there exist some positive constants $C_\Gamma, c_1^*$ and $c_2^*$ such that, for  $n \gtrsim \log (pq_1)q_1^{4\alpha_1+2} (\cM_1^X)^2$, the matrix $\widehat{\bGamma} = 
(n-L)^{-1}\widehat{\bD}^{-1}\widehat{\bZ}^{\T}\widehat{\bZ}\widehat{\bD}^{-1} \in \mathbb{R}^{(L+1)pq_1 \times (L+1)pq_1}$ satisfies, with probability greater than $1-c_1^*(pq_1)^{-c_2^*},$
\begin{equation} \label{eq_fof_RE}
\begin{split}
    \btheta^{\T}\widehat{\bGamma}\btheta \geq \tau_2 \|\btheta\|^2 -\tau_1\|\btheta\|_1^2 \quad \forall \btheta \in \mathbb{R}^{(L+1)pq_1},
\end{split}
\end{equation}
where  $\tau_1=C_\Gamma \cM_1^X q_1^{\alpha_1 +1}\sqrt{\log (pq_1)/n }$ and  $\tau_2=\underline{\mu}.$
\end{proposition}
(\ref{eq_fof_RE}) can be viewed as the functional extension of RE condition under the FPCA framework. Intuitively, it provides some insight into the eigenstructure of the sample correlation matrix of a vector formed by estimated lagged FPC scores of $\{X_{tj}(\cdot)\}_{j=1}^p.$ In particular, for any $\btheta \in \mathbb{R}^{(L+1)pq_1}$ such that $\tau_1\|\btheta\|_1^2 /\tau_2 \|\btheta\|^2$ is relatively small, $\btheta^{\T}\widehat\bGamma\btheta/ \|\btheta\|^2$ is bounded away from 0.
Proposition~\ref{pr_fof_RE} formalize this intuition by showing (\ref{eq_fof_RE}) holds with high probability.
Furthermore, Propositions~\ref{pr_fof_eigen} and \ref{pr_fof_max.error} verify  the essential deviation bounds for model~(\ref{model_fully}), where further discussions can be found in Appendix~\ref{ap.prop}.

Now we are ready to present the main convergence result.


\begin{theorem} \label{th_beta}
    Suppose that Conditions~\ref{con_stability}--\ref{con_fof_eigen_min} hold with $\tau_2\geq 32\tau_1q_1q_2s.$
    If $n \gtrsim \log (pq_1q_2)(q_1^{4\alpha_1+4} \vee q_2^{4\alpha_2+4})(\cM_1^X+ \cM^Y)^2,$ 
    then there exist some positive constants $c_1^*$ and $c_2^*$ such that, for any regularization parameter, $\lambda_n \geq 2 C_0 s  q_1^{1/2}\big\{(\cM_1^X+ \cM^\epsilon) \vee \cM^Y\big\}\{(  q_1^{\alpha_1+3/2} \vee   q_2^{\alpha_2+3/2})\sqrt{\frac{\log(pq_1q_2)}{n}}+ q_1^{-\kappa+1/2}\}$ and $q_1^{\alpha_1/2} s \lambda_n \to 0$ as $n,p,q_1,q_2 \to \infty,$ 
the estimate $\widehat\bbeta$ satisfies
\begin{equation}  \label{res_fof}
\begin{split}
\|\widehat{\bbeta}-\bbeta\|_{1} \lesssim  \frac{q_1^{\alpha_1/2} s \lambda_n}{\underline{\mu}},
\end{split}
\end{equation}
with probability greater than $1-c_1^*(pq_1q_2)^{-c_2^*}.$
\end{theorem}

\begin{remark} \label{rmark_model1}
\begin{enumerate}[(a)]
\item The error bound of $\widehat\bbeta$ under functional $\ell_1$ norm is determined by sample size ($n$), number of functional variables $(p),$ functional sparsity level ($s$) as well as internal parameters, e.g., the convergence rate in (\ref{res_fof}) is better when truncated dimensions ($q_1,q_2$), functional stability measures ($\cM_1^X,\cM^\epsilon, \cM^Y$), decay rates of the lower bounds for eigenvalues ($\alpha_1,\alpha_2$) in Condition~\ref{con_eigen} are small and decay rate of the upper bounds for basis coefficients ($\kappa$) in Condition~\ref{con_fof_beta_a} and curvature ($\underline{\mu}$) in (\ref{eq_fof_RE}) are large.

\item The serial dependence contributes the additional term $(\cM_1^{X} + \cM^{\epsilon}) \vee \cM^Y$ in the error bound. Specifically, the presence of $\cM_1^{X} + \cM^{\epsilon}$ is due to 
Proposition~\ref{pr_FPCscores_XE} under the independence assumption between $\{\bX_t(\cdot)\}$ and $\{\epsilon_t(\cdot)\},$ 
which is used to verify the deviation bound in
Proposition~\ref{pr_fof_max.error}. 
Moreover, provided that our estimation is based on the representation in (\ref{model_fully_matrix}), formed by eigenfunctions $\{\phi_m(\cdot)\}$ of $\Sigma_0^Y,$ the term $\cM^Y$ comes from the consistency analysis of $\{\widehat \phi_m\}$ in Proposition~\ref{pr_fof_eigen}. 

\item Note that the VFAR model 
can be rowwisely viewed as a special case of model~(\ref{model_fully}). The serial dependence in the error bound of the VFAR estimate is captured by $\cM^X_1$ partially due to its presence in the deviation bounds on estimated cross-covariance between response $\{\bX_t(\cdot)\}$ and covariates $\{\bX_{t-h}(\cdot): 1 \leq h \leq L\}.$ By contrast, the serial dependence effect in (\ref{res_fof}) partially comes from estimated cross-covariance between covariates $\{\bX_{t-h}(\cdot): 0 \leq h \leq L\}$ and error $\{\epsilon_t(\cdot)\}$ instead of that between $\{\bX_{t-h}(\cdot): 0 \leq h \leq L\}$ and response $\{Y_t(\cdot)\}$ due to the fact that $\{Y_t(\cdot)\}$ is completely determined by $\{\bX_{t-h}(\cdot): 0 \leq h \leq L\}$ and $\{\epsilon_t(\cdot)\}$ via (\ref{model_fully}) given $\bbeta.$ Specially, if $\cM^{\epsilon} \vee \cM^Y  \lesssim \cM_1^X,$ $q_1 \asymp q_2$ and $\alpha_1=\alpha_2,$ the rate in (\ref{res_fof}) is consistent to that of the VFAR estimate in \cite{guo2019}.
\end{enumerate}
\end{remark}

\section{High-dimensional partially functional linear regression}
\label{sec_partial}
This section is organized in the same manner as Section~\ref{sec_fully}. 
We first present the three-step procedure to estimate sparse functional and scalar coefficients in model~(\ref{model_scalar}) and then study the estimation consistency in the high-dimensional regime.

\subsection{Estimation procedure}
\label{sec_estimation_sof}
Consider partially functional linear regression model in (\ref{model_scalar}), where $\eulB(\cdot) = (\beta_{1}(\cdot),\dots,\beta_{p}(\cdot))^{\T}$ are functional coefficients of functional covariates $\{\bX_t(\cdot)\}_{t=1}^n$ and $\gamma = (\gamma_1, \dots, \gamma_d)^{\T}$ are regression coefficients of scalar covariates $\{\bZ_t\}_{t=1}^n$. $\{\epsilon_t\}_{t=1}^n$ are mean-zero errors from sub-Gaussian linear process, independent of
$\{\bZ_t\}$ from sub-Gaussian linear process and 
$\{\bX_t(\cdot)\}$ from sub-Gaussian functional linear process. 
To estimate $\eulB(\cdot)$ and $\gamma$ under large $p$ and $d$ scenario, we assume some sparsity patterns in model~(\ref{model_scalar}), i.e. $\eulB(\cdot)$ is functional $s_1$-sparse, with support $S_1 = \{j\in\{1,\dots,p\}:\|\beta_{j}\|\neq0\}$ and cardinality $s_1 = |S_1|,$ and $\gamma$ is  $s_2$-sparse, with support $S_2 = \{j\in\{1,\dots,d\}:\gamma_{j}\neq0\}$ and cardinality $s_2 = |S_2|.$  Here $s_1$ and $s_2$ are much smaller than dimension parameters, $p$ and $d,$ respectively.

Under the Karhunen-Lo\`eve expansion of each $X_{tj}(\cdot)$ as described in Section~\ref{sec.m1.est}, model~(\ref{model_scalar}) can be rewritten as 
\begin{equation} \nonumber
     Y_t = \sum_{j=1}^p \sum_{l=1}^{q_j}\zeta_{tjl} \langle \psi_{jl},\beta_{j}\rangle + \sum_{j=1}^d Z_{tj}\gamma_{j} + r_t+ \epsilon_t,
\end{equation}
where $r_{t} = \sum_{j=1}^p \sum_{l=q_j+1}^{\infty}\zeta_{tjl} \langle \psi_{jl},\beta_{j}\rangle.$ Let $\mathcal{Y} = (Y_1, \dots, Y_n)^{\T} \in \mathbb{R}^{n},$ $\mathcal{Z} = (\mathcal{Z}_1,\dots,\mathcal{Z}_d) \in \mathbb{R}^{n \times d},$ $\mathcal{Z}_j = (Z_{1j}, \dots, Z_{nj})^{\T} \in \mathbb{R}^n,$ $\gamma = (\gamma_1, \dots, \gamma_d)^{\T} \in \mathbb{R}^{d},$ $\mathcal{X}_j \in \mathbb{R}^{n \times q_j}$ with its row vectors given by $\bzeta_{1j},\dots,\bzeta_{nj}$ and $\Psi_{j}=\int_{\cU}\bpsi_j(u)\beta_{j}(u) du\in\mathbb{R}^{q_j}.$ Then we can represent model~(\ref{model_scalar}) in the following matrix form,
\begin{equation} \label{model2_mat}
    \mathcal{Y}= \sum_{j=1}^p\mathcal{X}_{j}\Psi_{j} + \mathcal{Z}\gamma + R + E, 
\end{equation}
where $R = (r_1, \dots, r_n)^{\T} \in \mathbb{R}^{n}$ and $E = (\epsilon_1, \dots, \epsilon_n)^{\T} \in \mathbb{R}^{n}$ correspond to the truncation and random errors, respectively.

Our proposed three-step estimation procedure proceeds as follows. We start with performing FPCA on each $\{X_{tj}(\cdot)\}_{t=1}^n,$ and hence obtain estimated FPC scores $\{\widehat \zeta_{tjl}\}$ and eigenfunctions $\{\widehat \psi_{jl}(\cdot)\}.$ Motivated from (\ref{model2_mat}), we then develop a regularized least square approach by incorporating a standardized group lasso penalty for $\{\Psi_j\}_{j=1}^p$ and the lasso penalty for $\gamma,$ aimed to shrink all elements in $\Psi_j$ of unimportant functional covariates and coefficients of unimportant scalar covariates to be exactly zero. Specifically, we consider minimizing the following criterion over $\Psi_1, \dots, \Psi_p$ and $\gamma :$
\begin{equation}
\label{lasso_partial}
 \frac{1}{2}\|\mathcal{Y} - \sum_{j=1}^p \widehat{\mathcal{X}}_{j}\Psi_{j} - \mathcal{Z}\gamma\|^2  + \lambda_{n1}\sum_{j=1}^p\|\widehat{\mathcal{X}}_{j}\Psi_{j}\| + \widetilde\lambda_{n2}\|\gamma\|_1,
\end{equation}
where $\widehat{\mathcal{X}}_j$ is the estimate of $\mathcal{X}_j,$ and $\lambda_{n1}, \widetilde \lambda_{n2}$ are non-negative regularization parameters. Let the minimizers of (\ref{lasso_partial}) be $\widehat\Psi_1, \dots, \widehat\Psi_p$ and $\widehat\gamma.$ Finally, our estimated functional coefficients are given by $\widehat{\beta}_{j}(\cdot) = \widehat{\bpsi}_j(\cdot)^{\T}\widehat{\Psi}_{j}$ for $j = 1, \dots, p.$

\subsection{Theoretical properties}
We start with some notation that will be used in this section. For a block vector $B = (b_{1}^{\T}, \dots, b_{p}^{\T})^{\T} \in \eR^{pq}$ with the $j$-th block $b_{j} \in {\eR}^{q},$ we define its $q$-block versions of $\ell_1$ and elementwise $\ell_{\infty}$ norms by $\|B\|_{1}^{(q)} = \sum_{j} \|b_{j}\|$ and $\|B\|_{\max}^{(q)} = \max_{j} \|b_{j}\|,$ respectively.
To simplify our notation, we denote $\alpha_1$ in Condition~\ref{con_eigen} by $\alpha$ and assume the same truncated dimension across $j=1, \dots, p,$ denoted by $q$.  Let $\widehat{\mathcal{X}} = (\widehat{\mathcal{X}}_{1},\dots,\widehat{\mathcal{X}}_{p})\in \mathbb{R}^{n\times pq},$ $\Psi = (\Psi_{1}^{\T},\dots,\Psi_{p}^{\T})^{\T}\in\mathbb{R}^{pq},$ $\widehat{D} = \text{diag}(\widehat{D}_1,\dots,\widehat{D}_p)\in \mathbb{R}^{pq\times pq},$ where $\widehat{D}_j = \{n^{-1}\widehat{\mathcal{X}}_{j}^{\T}\widehat{\mathcal{X}}_{j}\}^{1/2}\in\mathbb{R}^{q\times q}$ for $j=1,\dots,p.$ Then our minimizing task in (\ref{lasso_partial}) is equivalent to 
\begin{equation}
\label{target_partial}
(\widehat B,\widehat \gamma) = \underset{ B \in \mathbb{R}^{pq},\gamma \in \mathbb{R}^d}{ \text{arg min}} \left\{
\frac{1}{2n}\|\mathcal{Y} -  \widehat \Omega B - \mathcal{Z}\gamma\|^2  + \lambda_{n1}\|B\|_1^{(q)}+ \lambda_{n2}\|\gamma\|_1\right \},
\end{equation}
where $\widehat{\Omega} = 
\widehat{\mathcal{X}}\widehat{D}^{-1}$ and $\lambda_{n2}=\widetilde\lambda_{n2}/n.$ 
Then $\widehat{\Psi} = \widehat{D}^{-1}\widehat{B}$ with its $j$-th row block given by $\widehat{\Psi}_{j}.$ 

\begin{condition} \label{con_partial_beta_b}
For $j \in S_1,$ $\beta_{j}(u) = \sum_{l=1}^\infty a_{jl}\psi_{jl}(u)$ and there exist some positive constants $\kappa > \alpha /2 + 1$ and $\mu_j$ such that $|a_{jl}| \leq \mu_jl^{-\kappa}$ for $l \geq 1.$
\end{condition}

Condition~\ref{con_partial_beta_b} controls the level of smoothness for non-zero coefficient functions in $\eulB(\cdot).$ See also Condition~\ref{con_fof_beta_a} for model~(\ref{model_fully}) and its subsequent discussion.

\begin{condition} \label{con_partial_eigen_min_xz}
For the mixed process $\{\bX_t(\cdot), \bZ_t\}_{t \in \mathbb{Z}},$ we denote
a diagonal matrix function by $\bD_0^X = \text{diag} (\Sigma^X_{0,11},\dots,\Sigma^X_{0,pp}).$ The infimum $\underline{\mu}^*$ is bounded below by zero, i.e.
\begin{equation} \nonumber
\begin{split}
\underline{\mu}^* = \inf\limits_{\bPhi \in \Bar{\mathbb{H}}_0^p ,\bnu \in \widetilde{\mathbb{R}}_0^d } \frac{\langle\bPhi,\bSigma_0^{X}(\bPhi)\rangle+ \langle\bPhi,\bSigma_0^{X,Z}\bnu\rangle+ \bnu^{\T}\bSigma_0^{Z,X}(\bPhi) + \bnu^{\T}\bSigma_0^{Z}\bnu}{{ \langle\bPhi,\bD_0^X(\bPhi)\rangle +\bnu^{\T}\bnu}} > 0,
\end{split}
\end{equation}
where $\Bar{\mathbb{H}}_0^p = \{\bPhi \in \mathbb{H}^p: \langle\bPhi,\bD_0^X(\bPhi)\rangle \in (0,\infty)\}.$
\end{condition}


This condition is similar to Condition~\ref{con_fof_eigen_min}. In the special case where each $X_{tj}(\cdot)$ is $b_j$-dimensional, $\underline{\mu}^*$ reduces to the minimum eigenvalue of the covariance matrix of $\big(\frac{\xi_{t11}}{\sqrt{\omega^X_{11}}}, \dots, \\ \frac{\xi_{t1b_1}}{\sqrt{\omega^X_{1b_1}}}, \dots, \frac{\xi_{tp1}}{\sqrt{\omega^X_{p1}}}, \dots, \frac{\xi_{tpb_p}}{\sqrt{\omega^X_{pb_p}}},Z_{t1}, \dots, Z_{td}\big)^{\T} \in {\mathbb R}^{\sum_{j=1}^p b_j +d}.$

We next present Proposition~\ref{pr_partial_RE} below and Propositions~\ref{pr_partial_max.error_X}--\ref{pr_partial_max.error_z} in Appendix~\ref{ap.prop} to respectively show that the RE and deviation conditions are satisfied with high probability. 
These results together with Proposition~\ref{pr_fof_eigen}(i) lead to theoretical guarantees for regularized estimates of model (\ref{model_scalar}).


\begin{proposition} \label{pr_partial_RE}
Suppose Conditions~\ref{con_stability}--\ref{con_eigen} and \ref{con_partial_eigen_min_xz} hold.
Let $\mathcal{S} = ( \widehat \Omega , \mathcal{Z}) \in \mathbb{R}^{n \times( pq+d)},$
 then there exist some positive constants 
$C_{Z\Gamma}, c_1^*$ and $c_2^*$ such that, for $n\gtrsim\log (pqd)q^{4\alpha+2}\cM_{X,Z}^2,$  with probability greater than $1-c_1^*({pq+d})^{-c_2^*},$
\begin{equation} \label{RE_scalar}
\begin{split}
\frac{1}{n}\btheta^{\T}\mathcal{S}^{\T} \mathcal{S}\btheta \geq \tau_2^* \|\btheta\|^2 -\tau_1^*\|\btheta\|_1^2, 
~~ \forall\btheta \in \mathbb{R}^{pq+d},
\end{split}
\end{equation}
where $\tau_{1}^* =  C_{Z\Gamma} \cM_{X,Z} q^{\alpha +1}\sqrt{\frac{\log (pq+d)}{n}}$ and $\tau_{2}^*=\underline{\mu}^*.$
\end{proposition}

Instead of verifying RE conditions on  $n^{-1} \widehat \Omega^{\T}\widehat \Omega$ and $n^{-1} \mathcal{Z}^{\T}\mathcal{Z}$ separately, since $\widehat\Omega$ is correlated with $\mathcal{Z},$ we define $\mathcal{S} = ( \widehat \Omega , \mathcal{Z})$ and verify (\ref{RE_scalar}), which requires $n^{-1}\btheta^{\T}\mathcal{S}^{\T} \mathcal{S}\btheta$ to be strictly positive as long as $\tau_1^*\|\btheta\|_1^2/\tau_2^*\|\btheta\|^2$ is relatively small. Let $\btheta = ( \Delta ^{\T},\delta^{\T})^{\T}$ with $\Delta = \widehat{B} - B$ and $\delta = \widehat \gamma - \gamma,$ applying Proposition~\ref{pr_partial_RE} with suitable choice of $\tau_2^*$ yields that, with high probability, $n^{-1}(\widehat \Omega \Delta + \mathcal{Z}\delta)^{\T}(\widehat \Omega \Delta + \mathcal{Z}\delta) \geq \frac{\tau_2^*}{4}(\|\Delta\| + \|\delta\|)^2,$  which plays a crucial role in the proof of Theorem~\ref{th_beta_partial} below.
Similar to Proposition~\ref{pr_fof_max.error}, Propositions~\ref{pr_partial_max.error_X} and \ref{pr_partial_max.error_z} in Appendix~\ref{ap.prop} verify that, with high probability, the essential deviation bounds hold for model~({\ref{model_scalar}}). 

Now we are ready to present the main theorem about the error bound for $\widehat B$ and $\widehat\gamma.$ 
\begin{theorem} \label{th_beta_partial}
    Suppose that Conditions~\ref{con_stability}--\ref{con_eigen},  \ref{con_partial_beta_b} and \ref{con_partial_eigen_min_xz} hold with $\tau_2^*\geq 64\tau_1^*q(s_1+s_2).$
    If $n\gtrsim\log (pqd)q^{4\alpha+2}\cM_{X,Z}^2,$ then, for any regularization parameters, $\lambda_n \asymp \lambda_{n1} \asymp \lambda_{n2} \geq 2 C_{0}^* s_1 (\cM_{X,Z} + \cM^\epsilon)[q^{\alpha+2}\{\log (pq+d)/n\}^{1/2}+ q^{-\kappa+1}]$ with $q^{\alpha/2} \lambda_{n}  (s_1+s_2) \to 0$ as $n,p,q,d \to \infty,$
the estimates $\widehat \eulB$ and $\widehat \gamma$ satisfy
\begin{equation} \label{res_partial}
\begin{split}
\|\widehat{\eulB}-\eulB\|_{1} + q^{\alpha/2}\|\widehat \gamma - \gamma\|_1 \lesssim  \frac{q^{\alpha/2} \lambda_{n} (s_1+s_2) }{\underline{\mu}^*},
\end{split}
\end{equation}
with probability greater than $1-c_1^*(pq+d)^{-c_2^*}.$
\end{theorem}


\begin{remark}
\begin{enumerate}[(a)]
\item The error bound in (\ref{res_partial}) is governed by both dimensionality parameters ($n, p, d, s_1, s_2$) and internal parameters ($\cM^{X},\cM^Z, \cM^{X,Z}, \cM^{\epsilon}, q, \alpha,\kappa, \underline{\mu}^*$). See also similar Remark~\ref{rmark_model1}~(a) for model~(\ref{model_fully}). 

\item Note that the sparse stochastic regression \cite[]{basu2015a,Wu2016} can be viewed as a special case of model (\ref{model_scalar}) without the functional part. Under such scenario, the absence of $\{\bX_t(\cdot)\}$ 
degenerates  (\ref{eq_partial_max.error_z}) in Proposition~\ref{pr_partial_max.error_z} to $n^{-1} \|\cZ^{\T}({\mathcal Y}-\cZ\gamma)\|_{\max} \leq \widetilde C_{0}(\cM_1^Z + \cM^\epsilon)(\log d/n)^{1/2}$ and simplifies the error bound to $\|\widehat \gamma - \gamma\|_1 \lesssim  {\lambda_{n2}  s_2}/{ \tau_2^*}$ with $\lambda_{n2} \geq 2\widetilde C_{0}(\cM_1^Z + \cM^\epsilon)(\log d/n)^{1/2}$ for some positive constant $\widetilde C_{0},$ which is of the same order as the rate in \cite{basu2015a}.
\item In another special scenario where scalar covariates are not included in (\ref{model_scalar}), 
the error bound reduces to $\|\widehat{\eulB}-\eulB\|_{1}  \lesssim  {q^{\alpha/2}  \lambda_{n1}  s_1 }/{\tau_2^*}$ with $\lambda_{n1} \geq 2C_{0}^* s_1 (\cM_1^X + \cM^\epsilon)\{q^{\alpha+2}\sqrt{\frac{\log(pq)}{n}}+ q^{-\kappa+1}\}.$ Interestingly, this rate is consistent to that of $\widehat \bbeta$ in Theorem~\ref{th_beta} under the special case where the non-functional response results in the absence of $\cM^Y$ and $q_2$ in the rate. 

\end{enumerate}
\end{remark}

\section{Simulation studies} \label{sec.sim}
We conduct a number of simulations to evaluate the finite-sample performance of our proposed $\ell_1/\ell_2$-penalized least squares estimators ($\ell_1/\ell_2$-LS) for models~(\ref{model_fully}) and (\ref{model_scalar}) in Sections~\ref{sec.sim.fully} and \ref{sec.sim.scalar}, respectively.


\subsection{High-dimensional functional linear lagged regression} \label{sec.sim.fully}
We consider model  (\ref{model_fully}) with $L=1$, where functional covariates $\{\bX_t(\cdot)\}_{t=1, \dots, n}$ are generated from a sparse VFAR model \cite[]{guo2019}. 
Specifically, we generate $X_{tj}(u) = \bzeta_{tj}^{\T} \bpsi(u)$ for $j = 1, \dots, p$ and $u \in \cU = [0,1],$ where $\bpsi(\cdot) = (\psi_1(\cdot), \dots, \psi_5(\cdot))^{\T}$ is a 5-dimensional Fourier basis function and $\bzeta_t = (\bzeta_{t1}^{\T},\dots, \bzeta_{tp}^{\T})^{\T} \in \mathbb{R}^{5p}$ are generated from a stationary block sparse vector autoregressive (VAR) model, $\bzeta_t = \bW \bzeta_{t-1} + \boldsymbol{\eta}_t.$ The transition matrix $\bW =(\bW_{jk})_{p\times p}  \in \mathbb{R}^{5p \times 5p}$ is block sparse such that $\sum_{k=1}^p I(\|\bW_{jk}\|_{\tF} \neq 0)=5$ for each $j,$ and $\boldsymbol \eta_t$ are sampled independently from $N(\boldsymbol 0, \bI_{5p}).$
The nonzero elements in $\bW$ are sampled from $N(0,1)$ and we rescale $\bW$ by $\iota\bW/\rho(\bW)$ with $\iota \sim$ Unif[0.5,1] to guarantee the stationarity of $\{\bzeta_t\}.$
For each $(h,j) \in S = \{0,1\} \times \{1, \dots, 5\},$ we generate non-zero functional coefficients $\beta_{hj} (u,v) = \sum_{l,m = 1}^5 b_{hjlm}\psi_l(u)\psi_m(v),$  where $b_{hjlm}$'s are sampled from
 Unif$(0,0.4)$ for $h=0$ and Unif$(0,0.15)$ for $h=1.$ 
The functional responses $\{Y_{t}(v): v \in \cV\}_{t=1, \dots,n}$ with $\cV=[0,1]$ are then generated from model~(\ref{model_fully}), where $\epsilon_t(v) = \sum_{m = 1}^5 e_{tm} \psi_m(v)$ with $e_{tm}$'s being independent $N(0,1)$ variables. 

In our simulations, we consider $n=75, 100, 150$ dependent observations for $p=40, 80$ and replicate each simulation 100 times. The truncated dimensions $q_{1j}$ for $j=1, \dots p$ and $q_2$ are selected by the ratio-based method \cite[]{Lam2012}. To select the regularization parameter $\lambda_n,$ there exists several possible methods such as AIC/BIC and cross-validation. The AIC/BIC requires to specify the effective degrees of freedom, which poses a challenging task for functional data under the high-dimensional setting and is left for future study. In this example, we generate two separate training and validation samples of the same size $n.$ For a sequence of $\lambda_n$ values, we implement the block fast iterative shrinkage-thresholding (FISTA) algorithm \cite[]{guo2019} to solve the optimization problem (\ref{lasso}) on the training data, obtain $\{\widehat\beta_{hj}^{(\lambda_n)}(\cdot,\cdot)\}_{h=0,1, j=1, \dots,p}$ as a function of $\lambda_n,$ calculate the squared error between observed and fitted responses on the validation set, i.e. $\sum_{t = 1}^n\|Y_t(\cdot) - \sum_{h=0}^L \sum_{j=1}^p \int_{\cU}X_{(t-h)j}(u) \widehat \beta_{hj}^{(\lambda_n)}(u,\cdot)du \|^2$ and choose the optimal $\widehat \lambda_n$ with the smallest error. 


We evaluate the performance of $\ell_1/\ell_2$-LS in terms of both model selection consistency and estimation accuracy. For model selection consistency, we plot the true positive rates against false positive rates,  defined as $\frac{\#\{(h,j):  ||\widehat \beta_{hj}^{(\lambda_n)}||_{\cS}\neq0 \text{ and } ||\beta_{hj}||_{\cS}\neq 0\}}{\#\{(h,j): ||\beta_{hj}||_{\cS}\neq 0\}}$ and $\frac{\#\{(h,j):  ||\widehat \beta_{hj}^{(\lambda_n)}||_{\cS}\neq0 \text{ and } ||\beta_{hj}||_{\cS}= 0\}}{\#\{(h,j): ||\beta_{hj}||_{\cS}= 0\}}$ respectively, over a grid of values of $\lambda_n$ to produce a ROC curve, and then calculate the {\it area under the ROC curve} (AUROC) with values closer to 1 indicating better performance in support recovery. The estimation accuracy is measured by the relative estimation error  
$\|\widehat\bbeta - \bbeta \|_\tF/\|\bbeta\|_\tF.$ For comparison, we also implement the ordinary least squares in the oracle case (OLS-O), which uses the true sparsity structure in the estimates and does not perform variable selection. Table~\ref{table.FoF} gives some numerical summaries. Several conclusions can be drawn. 
First, the model selection consistency and estimation accuracy are improved as $n$ increases or $p$ decreases.
Second, $\ell_1/\ell_2$-LS provides substantially improved estimation accuracy over $\text{OLS-O}$ especially in the ``large $p,$ small $n$" scenario. This is not surprising, since implementing OLS-O in the sense of (\ref{lasso}) with $\lambda_n=0$ still require to estimate $10 \times 5^2=250$ parameters, which is intrinsically a high-dimensional estimation problem.

\begin{table}[t]
	\caption{\label{table.FoF} The mean and standard error (in parentheses) of AUROCs and estimation errors for model~(\ref{model_fully}) over 100 simulation runs. }
	\begin{center}
		\resizebox{4in}{!}{
			\begin{tabular}{ccccc}		
				\hline
		\multirow{2}{*}{$n$}	&	\multirow{2}{*}{$p$}	&	\multicolumn{2}{c}	{$\ell_1/\ell_2$-$\LS$ }		&	$\OLS$	\\	
	&		&	AUROC	&	Estimation error	&	Estimation error	\\	\hline
\multirow{2}{*}{75}	&	40	&	0.849(0.006)	&	0.727(0.005)	&	1.116(0.011)	\\	
	&	80	&	0.834(0.007)	&	0.768(0.005)	&	1.121(0.012)	\\	\hline
\multirow{2}{*}{100}	&	40	&	0.898(0.005)	&	0.648(0.005)	&	0.777(0.006)	\\	
	&	80	&	0.879(0.007)	&	0.684(0.005)	&	0.787(0.006)	\\	\hline
\multirow{2}{*}{150}	&	40	&	0.953(0.004)	&	0.544(0.004)	&	0.550(0.004)	\\	
	&	80	&	0.942(0.004)	&	0.576(0.004)	&	0.547(0.004)	\\	\hline
			\end{tabular}
		}	
	\end{center}
\end{table}

\subsection{High-dimensional partially functional linear regression}
\label{sec.sim.scalar}

We now consider model (\ref{model_scalar}) with $p$-dimensional vector of functional covariates $\{\bX_t(\cdot)\}_{t=1, \dots, n}$ and $d$-dimensional  scalar covariates $\{\bZ_t\}_{t=1,\dots,n},$ which are jointly generated in a similar procedure as in Section~\ref{sec.sim.fully}. 
Let $X_{tj}(u) = \bzeta_{tj}^{\T} \bpsi(u)$  for $j = 1, \dots, p$ and $u \in [0,1],$ 
and $(\bzeta_t^{\T} , \bZ_t^{\T})^{\T} \in \mathbb{R}^{5p+d}$ are jointly generated from a stationary VAR(1) process with a block sparse transition matrix $\bW^* \in \mathbb{R}^{(5p+d)\times (5p+d)},$ whose $(j,k)$-th block is  $\bW^*_{jk}.$ 
In particular, for each $j=1, \dots, p,$ $\bW_{jk}^* \in \mathbb{R}^{5\times 5}$ ($k=1, \dots, p$) and $\bW_{jk}^* \in \mathbb{R}^{5}$ ($k=p+1, \dots, p+d$) such that $\sum_{k=1}^p I(\|\bW_{jk}^*\|_{\tF} \neq 0)=\sum_{k=p+1}^{p+d}I(\|\bW_{jk}^*\| \neq 0)=5.$
For each $j=p+1, \dots, p+d,$ $(\bW^*_{jk})^{\T} \in \mathbb{R}^{5}$ ($k = 1, \dots, p$) and $\bW^*_{jk} \in \mathbb{R}$ ($k = p+1, \dots, p+d$) such that $\sum_{k=1}^p I\big(\|(\bW^*_{jk})^{\T}\| \neq 0\big)=\sum_{k=p+1}^{p+d}I(|\bW_{jk}^*| \neq 0)=5.$
For each $j \in S_1 = \{1,\dots, 5\},$ the non-zero functional coefficients are generated by $\beta_{j}(u) = \sum_{l= 1}^5 b_{jl}\psi_l(u),$ where $b_{jl}$'s are uniformly sampled from $[0,0.15].$ For each $k \in S_2 = \{1,\dots,10\},$ the non-zero scalar coefficients $\gamma_k$'s are uniformly sampled from $[0.5,1].$  
Finally, we generate responses $\{Y_t\}_{t = 1,\dots,n}$ from model~(\ref{model_scalar}), where $\epsilon_t$'s are sampled from $N(0,1).$

We simulate the data under six different settings, where $n \in \{75,100,150\}$ and $p=d \in \{40,80\},$ and replicate each simulation 100 times. 
For a sequence of pairs of $(\lambda_{n1}, \lambda_{n2}),$ following the procedure in Section~\ref{sec_estimation_sof}, we truncate each functional covariate with $q_{j}$ chosen by the ratio-based method, apply the block FISTA algorithm to minimize the criterion (\ref{lasso_partial}) on the training data and obtain $\{\widehat \beta_j^{(\lambda_{n1}, \lambda_{n2})}(\cdot)\}_{j=1, \dots,p}$ and
$\{\widehat \gamma_k^{(\lambda_{n1}, \lambda_{n2})}\}_{k=1, \dots,d}.$
The optimal regularization parameters $(\widehat \lambda_{n1},\widehat \lambda_{n2})$ are selected by minimizing the prediction error on the validation data with size $n,$ i.e. $\sum_{t = 1}^n \big\{Y_t - \sum_{j=1}^p \int_{\cU}X_{tj}(u)\widehat \beta_{j}^{(\lambda_{n1},\lambda_{n2})}(u)du - \sum_{k=1}^d Z_{tk}\widehat\gamma_{k}^{(\lambda_{n1},\lambda_{n2})}\big\}^2.$

We examine the performance of $\ell_1/\ell_2$-LS based on AUROCs and estimation errors, and compare it with the performance of OLS-O, where the sparsity structures in the estimates are determined by the true model in advance. The numerical results are summarized in Table~\ref{table.SoF}, where the relative estimation errors for functional and scalar coefficients are $\|\widehat \eulB - \eulB \|/\|\eulB \|$ and $\|\widehat \gamma - \gamma\|/\|\gamma\|,$ respectively. A few trends are apparent. First, as expected, we obtain improved overall support recovery and estimation accuracies as $n$ increases or $p$ and $d$ decrease. Second, although $\ell_1/\ell_2$-LS is outperformed by OLS-O with lower estimation errors for scalar coefficients, it provides more accurate estimates of functional coefficients relative to OLS-O, since, in the oracle case, the number of unknown parameters is still relatively large especially when $n$ is small.

\begin{table}[t]
	\caption{\label{table.SoF} The mean and standard error (in parentheses) of AUROCs and estimation errors for model~(\ref{model_scalar})  over 100 simulation runs.}
	\begin{center}
		\resizebox{4.8in}{!}{
			\begin{tabular}{ccccccc}
			\hline
				\multirow{2}{*}{$n$}		&	\multirow{2}{*}{$p=d$}		&	\multicolumn{3}{c}	{$\ell_1/\ell_2$-$\LS$ }				&			\multicolumn{2}{c}	{$\OLS$}						\\	
			&			&	AUROC	&	$\|\widehat \eulB - \eulB \|/\|\eulB \|$	&	$\|\widehat \gamma - \gamma\|/\|\gamma\|$	&			$\|\widehat \eulB - \eulB \|/\|\eulB \|$	&	$\|\widehat \gamma - \gamma\|/\|\gamma\|$					\\	\hline
\multirow{2}{*}{75}			&		40	&	0.901(0.004)	&	1.034(0.013)	&	0.283(0.005)	&			1.741(0.034)	&		0.196(0.005)				\\	
			&		80	&	0.868(0.004)	&	1.051(0.012)	&	0.363(0.008)	&			1.750(0.039)	&		0.198(0.005)				\\	\hline
\multirow{2}{*}{100}			&		40	&	0.919(0.003)	&	0.999(0.007)	&	0.235(0.005)	&			1.376(0.024)	&		0.151(0.004)				\\	
			&		80	&	0.902(0.004)	&	1.025(0.008)	&	0.283(0.005)	&			1.417(0.025)	&		0.151(0.004)				\\	\hline
\multirow{2}{*}{150}			&		40	&	0.945(0.003)	&	0.938(0.008)	&	0.185(0.004)	&			1.006(0.018)	&		0.113(0.003)				\\	
			&		80	&	0.937(0.004)	&	0.972(0.009)	&	0.216(0.004)	&			1.061(0.018)	&		0.113(0.003)				\\	\hline

			\end{tabular}
		}	
	\end{center}
\end{table}

\section{Discussion}
\label{sec.discuss}
We identify several directions for future study. First, it is possible to extend our established finite sample theory for stationary functional linear processes with sub-Gaussian errors to that with more general noise distributions, e.g. generalized sub-exponential process, or even non-stationary functional processes. 
Second, it is of interest to develop useful non-asymptotic results under other commonly adopted dependence framework, e.g. moment-based dependence measure \cite[]{hormann2010}
and different types of mixing conditions \cite[]{Bbosq1}. However, moving from standard asymptotic analysis to non-asymptotic analysis would pose complicated theoretical challenges.
Third, from a frequency domain perspective, it is interesting to study the non-asymptotic behaviour of smoothed periodogram estimators \cite[]{panaretos2013} for spectral density matrix function, served as the frequency domain analog of the sample covariance matrix function. Under a high-dimensional regime, it is also interesting to develop the functional thresholding strategy to estimate sparse spectral density matrix functions. These topics are beyond the scope of the current paper and will be pursued elsewhere.



\section*{Appendix}
\setcounter{equation}{0}
\renewcommand{\theequation}{A.\arabic{equation}}

\appendix
\section{Additional theoretical results}
\label{ap.prop}

We first present the following Propositions~\ref{pr_fof_eigen} and \ref{pr_fof_max.error}, in which we show that the essential deviation bounds for model~(\ref{model_fully}) are satisfied with high probability.

\begin{proposition} 
\label{pr_fof_eigen}
Suppose that Conditions~\ref{con_stability}--\ref{con_eigen} hold. Then there exist some positive constants $C_{\psi},$  $C_{\omega},$ $C_{\phi},$ $c_1^*$ and $c_2^*$ such that
(i) for $n \gtrsim \log (pq_1)q_1^{4\alpha_1+2} (\cM_1^X)^2,$
\begin{equation} \label{eq_fof_eigen_X}
\begin{split}
\mathop{\max}\limits_{1\leq j \leq p,1\leq l \leq q_1}\Bigg|\frac{\{\widehat{\omega}_{jl}^X\}^{-1/2} -\{\omega_{jl}^X\}^{-1/2}}{\{\omega_{jl}^X\}^{-1/2}}\Bigg| \leq C_{\omega}\mathcal{M}_1^X\sqrt{\frac{\log(pq_1)}{n}},\\\mathop{\max}\limits_{1\leq j \leq p,1\leq l \leq q_1}\|\widehat{\psi}_{jl} - \psi_{jl}\| \leq C_{\psi}\mathcal{M}_1^Xq_1^{\alpha_1+1}\sqrt{\frac{\log(pq_1)}{n}},
\end{split}
\end{equation}
with probability greater than $1-c_1^*\{pq_1 \}^{-c_2^*};$
(ii) for $n \gtrsim\log (q_2)q_2^{4\alpha_2+2} (\cM^Y)^2,$ 
\begin{equation} \label{eq_fof_eigen_Y}
\mathop{\max}\limits_{1\leq m \leq q_2}\|\widehat{\phi}_{m} - \phi_{m}\| \leq C_{\phi}\mathcal{M}^Yq_2^{\alpha_2+1}\sqrt{\frac{\log(q_2)}{n}},
\end{equation}
with probability greater than $1-c_1^*\{q_2 \}^{-c_2^*}.$
\end{proposition}

\begin{proposition} \label{pr_fof_max.error}
Suppose that Conditions~\ref{con_stability}--\ref{con_fof_beta_a} hold. Then there exist some positive constants $C_0, c_1^*$ and $c_2^*$ such that, for $n \gtrsim \log (pq_1q_2)(q_1^{4\alpha_1+4} \vee q_2^{4\alpha_2+4})(\cM_1^X+ \cM^Y)^2,$ 
\begin{equation} \label{eq_fof_max.error}
\begin{split}
 &(n-L)^{-1}\big\|\widehat{\bD}^{-1}\widehat{\bZ}^{\T}(\widehat{\bU}-\widehat{\bZ}\widehat{\bD}^{-1}\bB)\big\|_{\max}^{(q_1,q_2)}\\
     \leq&   C_0 s  q_1^{1/2}\big\{(\cM_1^X+ \cM^\epsilon) \vee \cM^Y\big\}\big\{(  q_1^{\alpha_1+3/2} \vee   q_2^{\alpha_2+3/2})\sqrt{\frac{\log(pq_1q_2)}{n}}+ q_1^{-\kappa+1/2}\big\},
\end{split}
\end{equation}
with probability greater than $1-c_1^*(pq_1q_2)^{-c_2^*}.$
\end{proposition}

(\ref{eq_fof_eigen_X}) and (\ref{eq_fof_eigen_Y}) in Proposition~\ref{pr_fof_eigen} control deviation bounds for relevant estimated eigenpairs of $X_{tj}(\cdot)$ and $Y_t(\cdot)$ under the FPCA framework. 
(\ref{eq_fof_max.error}) in Proposition~\ref{pr_fof_max.error} ensures that 
the sample cross-covariance between estimated lagged-and-normalized FPC scores and estimated errors consisting of truncated and random errors due to (\ref{model_fully_matrix}), are nicely concentrated around zero.

We next provide Propositions~\ref{pr_partial_max.error_X} and \ref{pr_partial_max.error_z}, where the essential deviation bounds for model~(\ref{model_scalar}) hold with high probability.

\begin{proposition} \label{pr_partial_max.error_X}
Suppose Conditions~\ref{con_stability}--\ref{con_eigen} and \ref{con_partial_beta_b} hold.   Then there exist some positive constants $C_{0}^*,$ $c_1^*$ and $c_2^*$ such that, for  $n\gtrsim \log (pq)q^{4\alpha+2}(\cM_{1}^X)^2,$  
\begin{equation} \label{eq_partial_max.error_X}
    \frac{1}{n} \| \widehat \Omega^{\T}(\mathcal{Y}  - \widehat \Omega B - \mathcal{Z}\gamma) \|_{\max}^{(q)} \leq C_{0}^* s_1 (\cM_1^X + \cM^\epsilon)\{q^{\alpha+2}\sqrt{\frac{\log(pq)}{n}}+ q^{-\kappa+1}\}.
\end{equation}
 with probability greater than $1-c_1^*(pq)^{-c_2^*}.$
\end{proposition}

\begin{proposition} \label{pr_partial_max.error_z}
Suppose Conditions~\ref{con_stability}--\ref{con_eigen} and \ref{con_partial_beta_b} hold.  Then there exist some positive constants $C_{0}^*, c_1^*$ and $c_2^*$ such that, for $n \gtrsim \log (pqd)q^{3\alpha+2}\cM_{X,Z}^2,$ 
\begin{equation} \label{eq_partial_max.error_z}
    \frac{1}{n} \| \mathcal{Z}^{\T}(\mathcal{Y} - \widehat \Omega B - \mathcal{Z}\gamma) \|_{\max} \leq C_{0}^* s_1 (\cM_{X,Z}+ \cM^\epsilon)\{q^{\alpha+1}\sqrt{\frac{\log(pq+d)}{n}}+ q^{-\kappa+1/2}\},
\end{equation} 
with probability greater than $1-c_1^*(pq+d)^{-c_2^*}.$
\end{proposition}

Intuitively, (\ref{eq_partial_max.error_X}) in Proposition~\ref{pr_partial_max.error_X} (or (\ref{eq_partial_max.error_z}) in Proposition~\ref{pr_partial_max.error_z})
indicates the sample cross-covariance between estimated normalized FPC scores (or scalar covariates) and estimated errors is nicely concentrated around zero.

\section{Proofs of theoretical results in Section~\ref{sec_non-asymptotic results}}
\label{ap.thm.sec2}

We provide proofs of theorems and propositions stated in Section~\ref{sec_non-asymptotic results} in Appendices~\ref{ap.th.secA1}--\ref{ap.pro.secA2}, followed by the supporting technical lemmas and their proofs in Appendix~\ref{ap.lemma.secA3}.
Throughout, we use $C_0, C_1, \dots,$ $c, c_1, \dots,$ $\tilde c_1, \tilde c_2, \dots,$ $\rho, \rho_1,  \rho_2, \dots$ to denote positive constants.
 For a matrix $\bB \in {\eR}^{p \times q},$ we denote its operator norm by  $||\bB||={\sup}_{||\bx||_2\leq 1}||\bB\bx||_2$. 
For $\phi_1, \phi_2 \in \mathbb{H}$ and $K \in \mathbb{S},$ we respectively denote $\int_\cU K(u,v) \phi_1(u)du,$ $\int_\cV K(u,v) \phi_2(v)dv$ and $\int_\cU\int_\cV K(u,v)\phi_1(u) \phi_2(v)dudv$ by $\langle  \phi_1, K \rangle,$ $\langle K, \phi_2 \rangle$ and $\langle \phi_1, \langle K, \phi_2 \rrangle$. For a fixed $\bPhi \in \mathbb{H}^p,$ we denote $\cM(\bbf^{X}, \bPhi)= 2\pi \cdot\text{ess} \sup_{\theta \in [-\pi,\pi]}|\langle\bPhi,\bbf^{X}_{\theta}(\bPhi)\rangle|.$

\subsection{Proofs of theorems}
\label{ap.th.secA1}
\paragraph{Proof of Theorem~\ref{th_XY}}
Part~(i): Define $\bY = (\langle\bPhi_1,\bX_{1}\rangle,\dots,\langle\bPhi_1,\bX_{n}\rangle)^{\T},$ then we obtain $|\langle\bPhi_1,(\widehat{\bSigma}_{0}^X - \bSigma_{0}^X)(\bPhi_1)\rangle|=\frac{1}{n} |\bY^{\T}\bY - \mathbb{E}(\bY^{\T}\bY)|.$ Our proof is organised as follows: We first introduce the $M$-truncated sub-Gaussian process $\bX_{M,L,t}(u) = \sum_{l=0}^{L}\bA_l(\bvarepsilon_{M,t-l}),$ where  $\varepsilon_{M,tj}(\cdot)  = \sum_{l=1}^M \sqrt{\omega_{jl}^{\varepsilon}}a_{tjl} \phi_{jl}(\cdot)$ for $j=1,\dots,p.$ 
We then apply the inequality in Lemma~\ref{lm_iidsub} on $\bX_{\infty, L,t}=\bX_{L,t}(u) = \sum_{l=0}^{L}\bA_l(\bvarepsilon_{t-l})$ by proving $\|\bPi_{M,L}\| \leq \cM(\bbf^X_{M,L}, \bPhi_1)$ and $\lim_{M \rightarrow \infty}\cM(\bbf^X_{M,L}, \bPhi_1) = \cM(\bbf^X_{L}, \bPhi_1).$  Finally, we will show that such inequality still holds as $L \to \infty.$ 

When $L$ and $M$ are both fixed, we first define $\bY_{M,L} = (\langle\bPhi_1,\bX_{M,L,1}\rangle,\dots, \langle\bPhi_1,\bX_{M,L,n}\rangle)^{\T}.$ Then $ \bY_{M,L}^{\T} \bY_{M,L}$ can be represented in the same form as $\langle\be_{M}, \bK(\be_{M})\rangle$ in Lemma~\ref{lm_iidsub}, where $\be_M = (\bvar_{M,n}^{\T}, \dots, \bvar_{M,1-L}^{\T})^{\T} \in \mathbb{H}^{(n+L)p}.$ We rewrite $\bY_{M,L}$ as
\begin{equation} \nonumber
    \bY_{M,L} = \int\int(\bI_n \otimes \bPhi_1(u)^{\T})\bW_L(u,v)\bTheta_M(v)\ba_{M,L}dudv = \bGamma_{M,L}\ba_{M,L},
\end{equation}
where 
\begin{equation} \nonumber
    \bW_L = \left( 
    \begin{matrix}
    \boldsymbol{0} & \boldsymbol{0} & \cdots & \boldsymbol{0} & \bA_0  & \cdots & \bA_{L-1} & \bA_L \\
    \boldsymbol{0} & \boldsymbol{0} & \cdots & \bA_0 & \bA_1& \cdots & \bA_L & \boldsymbol{0} \\
    \vdots & \vdots & \ddots  & \vdots & \vdots & \ddots & \vdots & \vdots  \\
    \bA_0 & \bA_1& \cdots & \cdots & \cdots & \bA_L & \cdots & \boldsymbol{0}
    \end{matrix}
    \right), 
\end{equation}
$\bTheta_M(u) = \bI_{n+L}\otimes \text{diag}(\bvarphi_{M,1}^{\T}, \dots, \bvarphi_{M,p}^{\T})$ with $\bvarphi_{M,i} = \big(\sqrt{\omega_{i1}^{e}}\phi_{i1},\dots,\sqrt{\omega_{iM}^{e}}\phi_{iM}\big)^{\T}$ and $\ba_{M,L} = (a_{n11},\dots, a_{n1M},\dots, a_{np1},\dots, a_{npM}, \dots, a_{(1-L)p1}, \dots, a_{(1-L)pM})^{\T} \in \mathbb{R}^{(n+L)pM}.$ Then we can write $ \bY_{M,L}^{\T} \bY_{M,L}=\ba_{M,L}^{\T}\bPi_{M,L}\ba_{M,L}$ with $\bPi_{M,L} = {\bGamma_{M,L}}^{\T}\bGamma_{M,L}.$ Lemma~\ref{lm_Q} implies that $\|\var(\bY_{M,L})\| = \|\bGamma_{M,L}\bGamma_{M,L}^{\T}\| \leq \cM(\bbf^X_{M,L}, \bPhi_1),$ 
where $\cM(\bbf^X_{M,L}, \bPhi_1)= 2\pi \cdot\text{ess} \sup_{\theta \in [-\pi,\pi]} \langle\bPhi_1,\bbf^X_{M,L,\theta}(\bPhi_1)\rangle$ and $\bbf^X_{M,L,\theta}(\cdot)$ is the spectral density matrix operator of process $\{\bX_{M,L,t}(\cdot)\}_{t \in \mathbb{Z}}.$ 

Define $\bY_L =\bY_{\infty,L}= (\langle\bPhi_1,\bX_{L,1}\rangle,\dots,\langle\bPhi_1,\bX_{L,n}\rangle)^{\T}.$
By Lemma~\ref{lm_middle matrix_ML},  (\ref{subG_2}) in Lemma~\ref{lm_iidsub} and $\text{rank}(\bGamma_{\infty,L}^{\T}\bGamma_{\infty,L}) = n$, we obtain 
\begin{equation*}
    \begin{split}
        &P\{|\langle\bPhi_1,(\widehat{\bSigma}_{L,0}^X - \bSigma_{L,0}^X)(\bPhi_1) | > \cM(\bbf^X_{L}, \bPhi_1)\eta\}
    \\= &P\{|\bY_L^{\T}\bY_L - \mathbb{E}\bY_L^{\T}\bY_L | > n\cM(\bbf^X_{L}, \bPhi_1)\eta\} 
    \leq 2 \exp \left\{ -cn \min \left(\eta^2,\eta\right)\right\},
    \end{split}
\end{equation*}
where $\cM(f_L^{X}, \bPhi_1)= 2\pi \cdot\text{ess} \sup_{\theta \in [-\pi,\pi]} \langle\bPhi_1,\bbf^{X}_{L,\theta}(\bPhi_1)\rangle$ and $\bbf^X_{L,\theta}(\cdot)$ is the spectral density matrix operator of $\{\bX_{L,t}(\cdot)\}_{t \in \mathbb{Z}}.$

Next, we need to show that this result still holds as $L \to \infty.$
Lemmas~\ref{lm_Y} and \ref{lm_M_L} imply that $\lim_{L \to \infty} \mathbb{E}\left\{\left|\langle\bPhi_1,(\widehat{\bSigma}_{L,0}^X - \widehat{\bSigma}_{0}^X)(\bPhi_1)\rangle \right|\right\} = 0,$ $\lim_{L \to \infty} \langle\bPhi_1,\bSigma_{L,0}^X(\bPhi_1)\rangle  = \langle\bPhi_1,\bSigma_{0}^X(\bPhi_1)\rangle$ and $\lim_{L \to \infty} \cM(f_L^X, \bPhi_1) =\cM(\bbf^X, \bPhi_1) .$ Combining the above results and following the similar argument in the proof of Lemma~\ref{lm_iidsub}, we obtain
\begin{equation} \label{eq_th1_Mxphi}
    \begin{split}
         P\left\{\left|\langle\bPhi_1,(\widehat{\bSigma}_{0}^X - \bSigma_{0}^X)(\bPhi_1)\rangle \right| > \cM(\bbf^X, \bPhi_1)\eta\right\} \\\leq 2 \exp \left\{ -cn \min \left(\eta^2,\eta\right)\right\}.
    \end{split}
\end{equation}
Provided that $\cM(\bbf^X, \bPhi_1) \leq \cM_{k_1}^X\langle\bPhi_1,\bSigma_{0}^X(\bPhi_1)\rangle,$ we obtain 
\begin{equation} \nonumber
    \begin{split}
         P\left\{\left|\frac{\langle\bPhi_1,(\widehat{\bSigma}_{0}^X - \bSigma_{0}^X)(\bPhi_1)\rangle}{\langle\bPhi_1,\bSigma_{0}^X(\bPhi_1)\rangle}\ \right| > \cM_{k_1}^X\eta\right\} \leq 2 \exp \left\{ -cn \min \left(\eta^2,\eta\right)\right\},
    \end{split}
\end{equation}
which completes the proof of (\ref{eq_th_X}).
\noindent
Part~(ii): For fixed vectors $\bPhi_1 \in \mathbb{H}^p$ and $\bPhi_2 \in \mathbb{H}^d,$ we denote $\cM(\bbf^{X,Y}, \bPhi_1, \bPhi_2)= 2\pi \cdot\text{ess} \sup_{\theta \in [-\pi,\pi]} |\langle\bPhi_1,\bbf^{X,Y}_{\theta}(\bPhi_2)\rangle|.$ Define $\bM_t(\cdot) = [(\bX_t(\cdot))^{\T},(\bY_t(\cdot))^{\T}]^{\T}.$ Letting $\bPhi = (\bPhi_1^{\T},\bPhi_2^{\T})^{\T},$ we have
\begin{equation}\nonumber
\begin{split}
    \langle\bPhi_1,(\widehat{\bSigma}_0^{X,Y}-\bSigma_0^{X,Y})(\bPhi_2)\rangle =&\frac{1}{2}[\langle\bPhi,(\widehat{\bSigma}_0^M-\bSigma_0^M)(\bPhi)\rangle
- \langle\bPhi_1,(\widehat{\bSigma}_0^X-\bSigma_0^X)(\bPhi_1)\rangle \\&-\langle\bPhi_2,(\widehat{\bSigma}_0^Y-\bSigma_0^Y)(\bPhi_2)\rangle].
\end{split}
\end{equation}
Applying (\ref{eq_th1_Mxphi}) on $\{\bX_t(\cdot)\}$ and $\{\bY_t(\cdot)\},$ we obtain that
\begin{equation*}
\begin{split}
P\left\{\left|\langle\bPhi_1,(\widehat{\bSigma}_0^X-\bSigma_0^X)(\bPhi_1)\rangle\right| > \mathcal{M}(\bbf^{X},\bPhi_1)\eta \right\}\leq 2\exp\{-cn\min(\eta^2,\eta)\},
\\
P\left\{\left|\langle\bPhi_2,(\widehat{\bSigma}_0^Y-\bSigma_0^Y)(\bPhi_2)\rangle\right| > \mathcal{M}(\bbf^{Y},\bPhi_2)\eta \right\}\leq 2\exp\{-cn\min(\eta^2,\eta)\}.
\end{split}
\end{equation*}
For $\{\bM_t(\cdot)\},$ we have
$
\mathcal{M}(f^{M},\bPhi) 
\leq\mathcal{M}(\bbf^{X},\bPhi_1) + \mathcal{M}(\bbf^{Y},\bPhi_2) + 2\mathcal{M}(\bbf^{X,Y}, \bPhi_1, \bPhi_2).
$
This, together with (\ref{eq_th1_Mxphi}) implies that
\begin{small}
\begin{equation}\nonumber
\begin{split}
P\left\{\left|\langle\bPhi,(\widehat{\bSigma}_0^M-\bSigma_0^M)(\bPhi)\rangle\right| > \{\mathcal{M}(\bbf^{X},\bPhi_1) + \mathcal{M}(\bbf^{Y},\bPhi_2) + 2\mathcal{M}(\bbf^{X,Y}, \bPhi_1, \bPhi_2)\}\eta \right\} 
\\ \leq 2\exp\{-cn\min(\eta^2,\eta)\}.
\end{split}
\end{equation}
\end{small}
Combining the above results, we obtain 
\begin{small}
\begin{equation}\label{eq_th1_Mxyphi}
\begin{split}
    P\left\{\left|\langle\bPhi_1,(\widehat{\bSigma}_0^{X,Y}-\bSigma_0^{X,Y})(\bPhi_2)\rangle\right| > \{\mathcal{M}(\bbf^{X},\bPhi_1) + \mathcal{M}(\bbf^{Y},\bPhi_2) + \mathcal{M}(\bbf^{X,Y}, \bPhi_1, \bPhi_2)\}\eta \right\}  
\\ \leq 6\exp\{-cn\min(\eta^2,\eta)\}.
\end{split}
\end{equation}
\end{small}

For $h>0,$ let $\bU_{1,t} = \bX_t + \bX_{t+h},$ $\bU_{2,t} = \bX_t - \bX_{t+h},$ $\bV_{1,t} = \bY_t + \bY_{t+h}$ and $\bV_{2,t} = \bY_t - \bY_{t+h}.$ Accordingly, we have that
\begin{equation}\nonumber
\begin{split}
    \langle\bPhi_1,\bSigma_l^{U_1,V_1}(\bPhi_2)\rangle &= 2\langle\bPhi_1,\bSigma_l^{X,Y}(\bPhi_2)\rangle + \langle\bPhi_1,\bSigma_{l-h}^{X,Y}(\bPhi_2)\rangle + \langle\bPhi_1,\bSigma_{l+h}^{X,Y}(\bPhi_2)\rangle,
    \\
     \langle\bPhi_1,\bSigma_l^{U_2,V_2}(\bPhi_2)\rangle &= 2\langle\bPhi_1,\bSigma_l^{X,Y}(\bPhi_2)\rangle - \langle\bPhi_1,\bSigma_{l-h}^{X,Y}(\bPhi_2)\rangle - \langle\bPhi_1,\bSigma_{l+h}^{X,Y}(\bPhi_2)\rangle,
\end{split}
\end{equation}
and 
\begin{equation}\nonumber
\begin{split}
     \bbf_\theta^{U_1,V_1} &= (2 + \text{exp}(-ih\theta) + \text{exp}(ih\theta))\bbf_\theta^{X,Y},
     \\
     \bbf_\theta^{U_2,V_2} &= (2 - \text{exp}(-ih\theta) - \text{exp}(ih\theta))\bbf_\theta^{X,Y}.
\end{split}
\end{equation}
Combining these with the definition of $\cM(\bbf^{X,Y},\bPhi_1,\bPhi_2)$ yields
\begin{equation}\nonumber
\begin{split}
     &4\langle\bPhi_1,(\widehat{\bSigma}_{h}^{X,Y}-\bSigma_{h}^{X,Y})(\bPhi_2)\rangle\\ =& \langle\bPhi_1,(\widehat{\bSigma}_0^{U_1,V_1}-\bSigma_0^{U_1,V_1})(\bPhi_2)\rangle -\langle\bPhi_1,(\widehat{\bSigma}_0^{U_2,V_2}-\bSigma_0^{U_2,V_2})(\bPhi_2)\rangle,
\end{split}
\end{equation}
and 
$$\cM(\bbf^{U_1,V_1},\bPhi_1,\bPhi_2) \leq 4 \cM(\bbf^{X,Y},\bPhi_1,\bPhi_2).$$
By similar arguments, we obtain $\cM(\bbf^{U_i},\bPhi_1) \leq 4 \cM(\bbf^{X},\bPhi_1)$ and $\cM(\bbf^{V_i},\bPhi_2) \leq 4 \cM(\bbf^{Y},\bPhi_2),$ for $i = 1,2.$ Then it follows from (\ref{eq_th1_Mxyphi}) that
\begin{footnotesize}
\begin{equation} \nonumber
\begin{split}
    &P\left\{\left|\langle\bPhi_1,(\widehat{\bSigma}_h^{X,Y}-\bSigma_h^{X,Y})(\bPhi_2)\rangle\right| > 2\{\mathcal{M}(\bbf^{X},\bPhi_1) + \mathcal{M}(\bbf^{Y},\bPhi_2) + \mathcal{M}(\bbf^{X,Y}, \bPhi_1, \bPhi_2)\}\eta \right\}  
\\ 
\leq &\sum_{i = 1}^2 P\left\{\left|\langle\bPhi_1,(\widehat{\bSigma}_0^{U_i,V_i}-\bSigma_0^{U_i,V_i})(\bPhi_2)\rangle\right| > \{\mathcal{M}(\bbf^{U_i},\bPhi_1) + \mathcal{M}(\bbf^{V_i},\bPhi_2) + \mathcal{M}(\bbf^{U_i,V_i}, \bPhi_1, \bPhi_2)\}\eta \right\}  
\\ \leq & 12\exp\{-cn\min(\eta^2,\eta)\}.
\end{split}
\end{equation}
\end{footnotesize}
Provided that $\mathcal{M}(\bbf^{X,Y},\bPhi_1,\bPhi_2) \leq \cM_{k_1,k_2}^{X,Y}(\langle\bPhi_1,\bSigma_0^X(\bPhi_1)\rangle +\langle\bPhi_2,\bSigma_0^Y(\bPhi_2)\rangle)$ and $\mathcal{M}(\bbf^{X},\bPhi_1) \leq \cM_{k_1}^X\langle\bPhi_1,\bSigma_0^X(\bPhi_1)\rangle,$ we obtain 
\begin{equation*}
\begin{split}
P\left\{\left|\frac{\langle\bPhi_1,(\widehat{\bSigma}_0^{X,Y}-\bSigma_0^{X,Y})(\bPhi_2)\rangle}{\langle\bPhi_1,\bSigma_0^X(\bPhi_1)\rangle +\langle\bPhi_2,\bSigma_0^Y(\bPhi_2)\rangle}\right| >  \left(\mathcal{M}_{k_1}^{X} +\mathcal{M}_{k_2}^{Y} +\mathcal{M}_{k_1,k_2}^{X,Y}\right)\eta \right\} \\ \leq 6\exp\left\{-cn\min(\eta^2,\eta)\right\},
\end{split}
\end{equation*}

\begin{equation*}
\begin{split}
P\left\{\left|\frac{\langle\bPhi_1,(\widehat{\bSigma}_h^{X,Y}-\bSigma_h^{X,Y})(\bPhi_2)\rangle}{\langle\bPhi_1,\bSigma_0^X(\bPhi_1)\rangle +\langle\bPhi_2,\bSigma_0^Y(\bPhi_2)\rangle}\right| >  2\left(\mathcal{M}_{k_1}^{X} +\mathcal{M}_{k_2}^{Y} +\mathcal{M}_{k_1,k_2}^{X,Y}\right)\eta \right\} \\ \leq 12\exp\{-cn\min(\eta^2,\eta)\}.
\end{split}
\end{equation*}

Letting $c_2 = c/4,$ we complete the proof of (\ref{eq_th_XY}).
$\square$
\\
\paragraph{Proof of Theorem~\ref{th_Dev}}
Under FPCA framework, for each $k = 1,\dots, d,$ we have $Y_{tk}(\cdot) = \sum_{m=1}^\infty\xi_{tkm}\phi_{km}(\cdot)$ with eigenpairs $(\omega_{km}^Y, \phi_{km}),$ and for each $j = 1,\dots, p,$ we have $X_{tj}(\cdot)=\sum_{l=1}^{\infty}\zeta_{tjl}\psi_{jl}(\cdot)$ with eigenpairs $(\omega_{jl}^X, \psi_{jl}).$ 

Denote $\cM_{X,Y} = \mathcal{M}_1^{X} +\mathcal{M}_1^{Y}+\mathcal{M}_{1,1}^{X,Y}.$ Let $\bPhi_1 = (0,\dots,0,\{\omega_{jl}^X\}^{-\frac{1}{2}}\psi_{jl},0,\dots,0)^{\T}$ and $\bPhi_2 = (0,\dots,0,\{\omega_{km}^Y\}^{-\frac{1}{2}}\phi_{km},0,\dots,0)^{\T}.$ Following the similar argument in the proof of Theorem~2 in \cite{guo2019} with $2\sqrt{\omega_0^X\omega_0^Y} \leq \omega_0^X + \omega_0^Y$ and Theorem~\ref{th_XY}, we can prove
\begin{equation}\nonumber
P\left\{\| \widehat{\Sigma}_{h,jk}^{X,Y}-\Sigma_{h,jk}^{X,Y}\|_\cS > (\omega_0^X + \omega_0^Y)\cM_{X,Y}\eta\right\} \leq c_1\exp\{-c_3n\min(\eta^2,\eta)\}.
\end{equation}
By the definition of $\| \widehat{\bSigma}_h^{X,Y}-\bSigma_h^{X,Y}\|_{\max} = \max_{1 \leq j \leq p, 1 \leq k \leq d} \| \widehat{\Sigma}_{h,jk}^{X,Y}-\Sigma_{h,jk}^{X,Y}\|_\cS $, we have that
\begin{equation}\nonumber
P\left\{\| \widehat{\bSigma}_h^{X,Y} -\bSigma_h^{X,Y}\|_{\max} > (\omega_0^X + \omega_0^Y)\cM_{X,Y}\eta\right\} \leq c_1pd\exp\{-c_3n\min(\eta^2,\eta)\}.
\end{equation}
Let $\eta = \rho\sqrt{\log (pd)/n} \leq 1$ and $\rho^2c_3 > 1,$ which can be achieved for sufficiently large $n.$ We obtain that 
\begin{equation}\nonumber
P\left\{\| \widehat{\bSigma}_h^{X,Y} -\bSigma_h^{X,Y}\|_{\max} > (\omega_0^X + \omega_0^Y)\cM_{X,Y}\rho\sqrt{\frac{\log (pd)}{n}}\right\} \leq c_1(pd)^{1-c\rho^2},
\end{equation}
which implies (\ref{eq_Dev_max}).
$\square$
\\

Before presenting the proof of Theorem~\ref{th_FPCscores_XY}, we provide some useful inequalities for estimated eigenpairs under the FPCA framework. For $\{\bX_t(\cdot)\}_{t\in \mathbb{Z}},$ let $\delta_{jl}^X = \min_{1\leq l' \leq l}\{\omega_{jl'}^X - \omega_{j(l'+1)}^X\}$ and $\widehat\Delta_{jl}^X = \widehat\Sigma_{0,jl}^X - \Sigma_{0,jl}^X$ for $j = 1,\dots,p$ and $l = 1,2,\dots.$ It follows from (4.43) and Lemma~4.3 of \cite{Bbosq1} that
\begin{equation} \label{con_6}
    \sup\limits_{l \geq 1} |\widehat\omega_{jl}^X - \omega_{jl}^X| \leq \| \widehat\Delta_{jj}^X\|_\cS \quad \text{and} \quad \sup\limits_{l \geq 1}\delta_{jl}^X\|\widehat\psi_{jl} - \psi_{jl} \| \leq 2\sqrt{2} \|\widehat\Delta_{jj}^X \|_\cS.
\end{equation}
Similarly, for process $\{\bY_t(\cdot)\}_{t \in \mathbb{Z}},$ let $\delta_{km}^Y = \min_{1\leq m' \leq m}\{\omega_{km'}^Y - \omega_{k(m'+1)}^Y\}$ and $\widehat{\Delta}_{km}^Y = \widehat{\Sigma}_{0,km}^Y - \Sigma_{0,km}^Y$ for $k = 1, \dots, d$ and $m = 1,2,\dots,$ we have
\begin{equation} \label{con_7}
    \sup\limits_{m \geq 1} |\widehat{\omega}_{km}^Y - \omega_{km}^Y| \leq \|\widehat{\Delta}_{kk}^Y\|_\cS \quad \text{and} \quad \sup\limits_{m \geq 1}\delta_{km}^Y\|\widehat{\phi}_{km} - \phi_{km}\| \leq 2\sqrt{2}\|\widehat{\Delta}_{kk}^Y\|_\cS.
\end{equation}

\paragraph{Proof of Theorem~\ref{th_FPCscores_XY}}
Recall that $\widehat{\sigma}_{h,jklm}^{X,Y}  = \frac{1}{n-h}\sum_{t=1}^{n-h}\widehat{\zeta}_{tjl}\widehat{\xi}_{(t+h)km}$ and $\sigma_{h,jklm}^{X,Y} =\\ \cov(\zeta_{tjl},\xi_{(t+h)km}) = \langle \psi_{jl},\langle\Sigma_{h,jk}^{X,Y},\phi_{km}\rrangle.$ Let $\widehat{r}_{jl} = \widehat{\psi}_{jl} - \psi_{jl},$ $\widehat{w}_{km} = \widehat{\phi}_{km} - \phi_{km}$ and $\widehat{\Delta}_{h,jk}^{X,Y} = \widehat{\Sigma}_{h,jk}^{X,Y}-\Sigma_{h,jk}^{X,Y},$ then
        \begin{equation} \nonumber
		\begin{split}
		\widehat{\sigma}_{h,jklm}^{X,Y}  - \sigma_{h,jklm}^{X,Y}  &= \langle \widehat{r}_{jl}, \langle \widehat{\Sigma}_{h,jk}^{X,Y}, \widehat{w}_{km}\rrangle + \left(\langle \widehat{r}_{jl}, \langle \widehat{\Delta}_{h,jk}^{X,Y},\phi_{km} \rrangle + \langle \psi_{jl}, \langle \widehat{\Delta}_{h,jk}^{X,Y},\widehat{w}_{km} \rrangle \right) 
		\\& \quad + \left(\langle\widehat{r}_{jl}, \langle \Sigma_{h,jk}^{X,Y}, \phi_{km} \rrangle + \langle \psi_{jl}, \langle \Sigma_{h,jk}^{X,Y},\widehat{w}_{km}\rrangle \right) + \langle \psi_{jl},\langle \widehat{\Delta}_{h,jk}^{X,Y}, \phi_{km}\rrangle
		\\& = I_1 + I_2 + I_3 + I_4.
		\end{split}
		\end{equation}
Let $\Omega_{jk,\eta}^{X,Y} = \left\{ \|\widehat{\Delta}_{h,jk}^{X,Y}\|_\cS \leq (\omega_0^X + \omega_0^Y)\cM_{X,Y}\eta\right\},$ $\Omega_{jj,\eta}^X = \left\{ \|\widehat{\Delta}_{jj}^X\|_\cS  \leq 2\mathcal{M}_1^{X}\omega_0^X\eta \right\},$ $\Omega_{kk,\eta}^Y =\\ \left\{ \|\widehat{\Delta}_{kk}^Y\|_\cS  \leq 2\mathcal{M}_1^{Y}\omega_0^Y\eta \right\}$ and $\Omega_1= \left\{ \|\widehat{\Delta}_{h,jk}^{X,Y}\|_\cS \leq (\omega_0^X + \omega_0^Y)\right\}.$ By Theorem~\ref{th_Dev} and Lemma~\ref{lm_guo2019_elementwise}, we have 
		\begin{equation} \nonumber
		\begin{split} 
		P((\Omega_{jk,\eta}^{X,Y})^C) &\leq c_1\exp\{-c_3n\min(\eta^2,\eta)\},
\\
		P((\Omega_{jj,\eta}^X)^C) &\leq 4\exp\{-\tilde{c}_1 n\min(\eta^2,\eta)\},
\\
		P((\Omega_{kk,\eta}^Y)^C) &\leq 4\exp\{-\tilde{c}_1 n\min(\eta^2,\eta)\},
\\
		P((\Omega_1)^C) &\leq c_1\exp\{-c_3n(\cM_{X,Y})^{-2}\}.
		\end{split}
		\end{equation}
On the event of $\Omega_1 \cap \Omega_{\eta,jk}^{X,Y} \cap \Omega_{jj,\eta}^X \cap \Omega_{kk,\eta}^Y ,$ by Condition~\ref{con_eigen}, (\ref{con_6}), (\ref{con_7}), Lemma~\ref{lm_X_traceclass} and the fact that $({\omega_0^X \omega_0^Y})^{1/2} \leq 1/2(\omega_0^X + \omega_0^Y),$ we obtain that
		\begin{equation}\nonumber
		\begin{split}
		\left| \frac{I_1}{\sqrt{\omega_{jl}^X \omega_{km}^Y}}\right| &\leq c_0^{-1}(\alpha_1\alpha_2)^{1/2} l^{\alpha_1/2}m^{\alpha_2/2}\|\widehat{r}_{jl}\|(\|\widehat{\Delta}_{h,jk}^{X,Y}\|_\cS +\|\Sigma_{h,jk}^{X,Y}\|_\cS)\|\widehat{w}_{km}\|
		\\& \lesssim  l^{3\alpha_1/2+1}m^{3\alpha_2/2+1} \|\widehat{\Delta}_{jj}^X\|_\cS \|\widehat{\Delta}_{kk}^Y\|_\cS(\|\widehat{\Delta}_{h,jk}^{X,Y}\|_\cS + ({\omega_0^X \omega_0^Y})^{1/2})
		\\& \lesssim  (l^{3\alpha_1+2} \vee m^{3\alpha_2+2})\mathcal{M}_1^{X}\mathcal{M}_1^{Y}\eta^2,
		\\& \lesssim  (l^{3\alpha_1+2}\vee m^{3\alpha_2+2})(\cM_1^X + \cM_1^Y)^2\eta^2,
		\end{split}
		\end{equation}
		\begin{equation}\nonumber
		\begin{split}
		\left| \frac{I_2}{\sqrt{\omega_{jl}^X \omega_{km}^Y}}\right| &\leq c_0^{-1}(\alpha_1\alpha_2)^{1/2} l^{\alpha_1/2}m^{\alpha_2/2}\|\widehat{\Delta}_{h,jk}^{X,Y}\|_\cS (\|\widehat{r}_{jl}\| + \|\widehat{w}_{km}\|)
		\\& \lesssim l^{\alpha_1/2}m^{\alpha_2/2}\|\widehat{\Delta}_{h,jk}^{X,Y}\|_\cS(l^{\alpha_1+1}\|\widehat{\Delta}_{jj}^X\|_\cS + m^{\alpha_2+1}\|\widehat{\Delta}_{kk}^Y\|_\cS)
		\\& \lesssim  (l^{2\alpha_1+1} \vee m^{2\alpha_2+1})\cM_{X,Y}(\cM_{X} \vee \cM_{Y})\eta^2,
		\\& \lesssim  (l^{2\alpha_1+1} \vee m^{2\alpha_2+1})\cM_{X,Y}^2\eta^2,
		\end{split}
		\end{equation}
By Theorem~\ref{th_XY},	
		\begin{equation}\nonumber
		P\{\left| \frac{I_4}{\sqrt{\omega_{jl}^X \omega_{km}^Y}}\right| \geq 2\cM_{XY}\eta\}
		\leq c_1\exp\{-c_2n\min(\eta^2,\eta)\}.
		\end{equation}
		Next, we consider the term $I_3=\langle\widehat{r}_{jl}, \langle \Sigma_{h,jk}^{X,Y}, \phi_{km} \rrangle + \langle \psi_{jl}, \langle \Sigma_{h,jk}^{X,Y},\widehat{w}_{km}\rrangle .$ By Condition~\ref{con_eigen}, Lemmas~\ref{lm_sigma_h_xy} and \ref{lm_guo2019_phi_g} for $\{\bX_t\}_{t \in \mathbb{Z}}$ and $\{\bY_t\}_{t \in \mathbb{Z}},$ we obtain that
		\begin{small}
		\begin{equation} \nonumber
		\begin{split} 
		&\left| \frac{I_3}{\sqrt{\omega_{jl}^X \omega_{km}^Y}}\right| \\\lesssim & \mathcal{M}_1^{X}l^{\alpha_1+1}\eta + \mathcal({M}_1^X)^2l^{(5\alpha_1+4)/2}\eta^2 + \mathcal{M}_1^{Y}m^{\alpha_2+1}\eta + (\mathcal{M}_1^Y)^2m^{(5\alpha_2+4)/2}\eta^2
		\\ \lesssim &(l^{\alpha_1+1}\vee m^{\alpha_2+1})(\cM_1^X + \cM_1^Y)\eta
+ (l^{(5\alpha_1+4)/2} \vee m^{(5\alpha_2+4)/2})(\cM_1^X + \cM_1^Y)^2\eta^2
		\end{split}
		\end{equation}
				\end{small}
		holds with probability greater than $1-16 \exp\{-\tilde c_4 n\min(\eta^2,\eta)\} - 8 \exp\{-\tilde c_4 n \\ (\{\cM_1^X\}^2l ^{2(\alpha_1+1)}\vee \{\cM_1^Y\}^2m^{2(\alpha_2+1)})^{-1}\}.$

		Combining the above results, we obtain that there exists positive constants $\rho_1,$ $\rho_2,$ $\tilde c_7$ and $\tilde c_8$ such that 
		\begin{small}
		\begin{equation} \nonumber
		\begin{split} 
		&P\left\{\left| \frac{ \widehat{\sigma}_{h,jklm}^{X,Y}  - \sigma_{h,jklm}^{X,Y} }{\sqrt{\omega_{jl}^X \omega_{km}^Y}}\right| \geq \rho_1 \cM_{X,Y} (l^{\alpha_1+1}\vee m^{\alpha_2+1})\eta + \rho_2\mathcal{M}_{X,Y}^2(l^{3\alpha_1+2} \vee m^{3\alpha_2+2})\eta^2\right\}
		\\& \quad \leq \tilde c_8 \exp\{-\tilde c_7 n \min(\eta^2,\eta)\} + \tilde c_8 \exp\{-\tilde c_7 \mathcal{M}_{X,Y}^{-2}n (l^{2(\alpha_1+1)} \vee m^{2(\alpha_2+1)})^{-1}\},
		\end{split}
		\end{equation}
		\end{small}
		where $\mathcal{M}_{X,Y} = \mathcal{M}_1^{X} +\mathcal{M}_1^{Y} +\mathcal{M}_{1,1}^{X,Y}.$ Applying the Boole's inequality, we obtain that 
		\begin{footnotesize}
		\begin{equation} \nonumber
		\begin{split} 
		&P\left\{\underset{ 1 \leq m \leq M_2}{\underset{1 \leq l \leq M_1}{\underset{ 1 \leq k \leq d}{\underset{1 \leq j \leq p}{\max}}}}\left| \frac{ \widehat{\sigma}_{h,jklm}^{X,Y}  - \sigma_{h,jklm}^{X,Y} }{\sqrt{\omega_{jl}^X \omega_{km}^Y}}\right| \geq  \rho_1 \cM_{X,Y}(l^{\alpha_1+1}\vee m^{\alpha_2+1})\eta + \rho_2\mathcal{M}_{X,Y}^2(l^{3\alpha_1+2} \vee m^{3\alpha_2+2})\eta^2\right\}
		\\& \quad \leq pdM_1M_2\{\tilde c_8 \exp\{-\tilde c_7  n \min(\eta^2,\eta)\} + \tilde c_8  \exp\{-\tilde c_7 \mathcal{M}_{X,Y}^{-2}n (l^{2(\alpha_1+1)} \vee m^{2(\alpha_2+1)^{-1}})\}.
		\end{split}  
		\end{equation}
		\end{footnotesize}
		Letting $\eta = \rho_3\sqrt{\frac{\log(pdM_1M_2)}{n}} < 1$ and $\rho_1 + \rho_2\rho_3\mathcal{M}_{X,Y}(M_1^{2\alpha_1+1} \vee M_2^{2\alpha_2+1}) \eta\leq \rho_4,$ there exist some constants $c_5, c_6 >0$ such that
		\begin{small}
		\begin{equation} \nonumber
		\begin{split} 
		&P\left\{\underset{1 \leq l \leq M_1,1 \leq m \leq M_2}{\underset{1 \leq j \leq p, 1 \leq k \leq d}{\max}}\left| \frac{ \widehat{\sigma}_{h,jklm}^{X,Y}  - \sigma_{h,jklm}^{X,Y} }{\sqrt{\omega_{jl}^X \omega_{km}^Y}}\right| \geq \rho_3\rho_4 \mathcal{M}_{X,Y}(M_1^{\alpha_1+1} \vee M_2^{\alpha_2+1})\sqrt{\frac{\log(pdM_1M_2)}{n}}\right\}
		\\& \quad \leq c_5(pdM_1M_2)^{c_6}.
		\end{split}  
		\end{equation}
		\end{small}
$\square$	

\subsection{Proofs of propositions}
\label{ap.pro.secA2}
\paragraph{Proof of Proposition~\ref{pr_FPCscores_XZ}}
Under a mixed-process scenario consisting of $\{\bX_t(\cdot)\}$ and $d$-dimensional time series $\{\bZ_t\},$ we obtain the concentration bound on $\widehat{\bSigma}_h^{X,Z},$
\begin{small}
\begin{equation} \label{eq_th_XZ}
P\left\{\left|\frac{\langle\bPhi_1,(\widehat{\bSigma}_h^{X,Z}-\bSigma_h^{X,Z})\bnu\rangle}{\langle\bPhi_1,\bSigma_0^X(\bPhi_1)\rangle +\bnu^{\T}\bSigma_0^Z\bnu}\right| >  \left(\mathcal{M}_{k_1}^{X} +\mathcal{M}_{k_2}^{Z} +\mathcal{M}_{k_1,k_2}^{X,Z}\right)\eta \right\}  \leq c_1\exp\{-c_2n\min(\eta^2,\eta)\}.
\end{equation}
\end{small}
Provided with Lemma~\ref{lm_basu_cov}, the above result can be proved in similar way to (\ref{eq_th_XY}) in Theorem~\ref{th_XY}, hence we omit it here. 

Denote $\sigma_{0,kk}^Z= \sqrt{\var(Z_k)},$ $(\sigma_{0}^Z)^2 = \max_{1 \leq k \leq d} {\var(Z_k)} < \infty$ and $\cM_{X,Z} = \mathcal{M}_1^{X} +\mathcal{M}_1^{Z}+\mathcal{M}_{1,1}^{X,Z}.$ Letting $\bPhi_1 = (0,\dots,0,\{\omega_{jl}^X\}^{-\frac{1}{2}}\psi_{jl},0,\dots,0)^{\T}$ and $\bnu = (0,\dots,0,\{\sigma_{0,kk}^Z\}^{-1},0,\dots,0)^{\T},$ we obtain that $\Delta_{h,jkl} = \langle\bPhi_1,(\widehat{\bSigma}_h^{X,Z}-\bSigma_h^{X,Z})\bnu\rangle = (\omega_{jl}^X )^{-1/2}(\sigma_{0,kk}^Z)^{-1}\langle\psi_{jl},\widehat \Sigma_{h,jk}^{X,Z}-\Sigma_{h,jk}^{X,Z}\rangle$ and 
$\langle\bPhi_1,\bSigma_0^X(\bPhi_1)\rangle =\bnu^{\T}\bSigma_0^Z\bnu = 1.$ Then $\|\Sigma_{h,jk}^{X,Z}-\Sigma_{h,jk}^{X,Z}\|^2 = \sum_{l = 1}^\infty \omega_{jl}^X(\sigma_{0,kk}^Z)^2\Delta_{h,jkl}^2.$ By Jensen's inequality, we have that
\begin{equation} \nonumber
\begin{split}
\mathbb{E}\Big\{ \big\|\widehat{\Sigma}_{h,jk}^{X,Z} - \Sigma_{h,jk}^{X,Z} \big\|_{\mathcal{S}}^{2q}\Big\}
&\le (\sigma_{0,kk}^Z)^{2q} \Big( \sum_{l= 1}^{\infty} \omega_{jl}^X \Big)^ {q-1} \sum_{l = 1}^{\infty}
\omega_{jl}^X  \mathbb{E}\big|\Delta_{h,jkl}  \big|^{2q}
\\&\le  \{\sigma_0^Z\}^{2q}\{\omega_0^X\}^{q} \sup_{l} \tE\big|\Delta_{h,jkl}  \big|^{2q}.
\end{split}
\end{equation}
By (\ref{eq_th_XZ}), we obtain that
\begin{equation} \nonumber
P\left\{\left|\Delta_{h,jkl}\right| >  2\cM_{X,Z}\eta \right\}  \leq c_1\exp\{-c_2n\min(\eta^2,\eta)\}.
\end{equation}
Combining the above results and following the similar argument in the proof of Theorem~2 in \cite{guo2019} yields 
\begin{equation}\nonumber
P\left\{\| \widehat{\Sigma}_{h,jk}^{X,Z}-\Sigma_{h,jk}^{X,Z}\| > 2\cM_{X,Z}\sigma_0^Z\sqrt{\omega_0^X}\eta\right\} \leq c_1\exp\{-c_3n\min(\eta^2,\eta)\}.
\end{equation}
Then with the fact that $2\sqrt{(\sigma_0^Z)^2\omega_0^X} \leq (\sigma_0^Z)^2 + \omega_0^X,$ we obtain 
\begin{equation}\label{eq_Dev_XZ}
P\left\{\| \widehat{\Sigma}_{h,jk}^{X,Z}-\Sigma_{h,jk}^{X,Z}\| > ((\sigma_0^Z)^2 + \omega_0^X)\cM_{X,Z}\eta\right\} \leq c_1\exp\{-c_3n\min(\eta^2,\eta)\}.
\end{equation}
This also implies (\ref{eq_Dev_max_XZ}).

Recall that $\widehat\varrho_{h,jkl}^{X,Z} = \frac{1}{n-h}\sum_{t=1}^{n-h}\widehat{\zeta}_{tjl}Z_{(t+h)k} $ and $\varrho_{h,jkl}^{X,Z} =\cov({\zeta}_{tjl},Z_{(t+h)k}).$ Let $\widehat{r}_{jl} = \widehat{\psi}_{jl} - \psi_{jl}$ and $\widehat{\Delta}_{h,jk}^{X,Z} = \widehat{\Sigma}_{h,jk}^{X,Z}-\Sigma_{h,jk}^{X,Z}.$ We have 
		\begin{equation} \nonumber
		\begin{split}
		\widehat{\varrho}_{h,jkl}^{X,Z}  - \varrho_{h,jkl}^{X,Z}  &= \langle \widehat{r}_{jl},  \widehat{\Delta}_{h,jk}^{X,Z}\rangle + \langle\widehat{r}_{jl},  \Sigma_{h,jk}^{X,Z}\rangle  + \langle \psi_{jl}, \widehat{\Delta}_{h,jk}^{X,Z} \rangle
		\\& = I_1 + I_2 + I_3. 
		\end{split}
		\end{equation}
Let $\Omega_{jk,\eta}^{X,Z} = \left\{ \|\widehat{\Delta}_{h,jk}^{X,Z}\| \leq (\omega_0^X + (\sigma_0^Z)^2)\cM_{X,Z}\eta\right\},$ $\Omega_{jj,\eta}^X = \left\{ \|\widehat{\Delta}_{jj}^X\|_\cS  \leq 2\mathcal{M}_1^{X}\omega_0^X\eta \right\}$ and $\Omega_1= \left\{ \|\widehat{\Delta}_{h,jk}^{X,Z}\|\leq (\omega_0^X + (\sigma_0^Z)^2)\right\}.$ By (\ref{eq_Dev_XZ}) and Lemma \ref{lm_guo2019_elementwise}, we have 
		\begin{equation} \nonumber
		\begin{split} 
		P((\Omega_{jk,\eta}^{X,Y})^C) &\leq c_1\exp\{-c_3n\min(\eta^2,\eta)\},
\\
		P((\Omega_{jj,\eta}^X)^C) &\leq 4\exp\{-\tilde{c}_1 n\min(\eta^2,\eta)\},
\\
		P((\Omega_1)^C) &\leq c_1\exp\{-c_3n(\cM_{X,Z})^{-2}\}.
		\end{split}
		\end{equation}
On the event of $\Omega_1 \cap \Omega_{\eta,jk}^{X,Z} \cap \Omega_{jj,\eta}^X  ,$ by Condition~\ref{con_eigen}, (\ref{con_6}), Lemma~\ref{lm_X_traceclass} and $(\sigma_0^Z)^2 < \infty$, we obtain that
		\begin{equation}\nonumber
		\begin{split}
		\left| \frac{I_1}{\sqrt{\omega_{jl}^X}}\right| &\lesssim  l^{\alpha_1/2}\|\widehat{\Delta}_{h,jk}^{X,Z}\| \|\widehat{r}_{jl}\|  \lesssim l^{3\alpha_1/2+1}\|\widehat{\Delta}_{h,jk}^{X,Z}\|\|\widehat{\Delta}_{jj}^X\|_\cS 
		\\& \lesssim  l^{3\alpha_1/2+1} \cM_{X,Z}\cM_1^{X} \eta^2.
		\end{split}
		\end{equation}
By Condition~\ref{con_eigen}, Lemma~\ref{lm_guo2019_phi_g} and $\|\Sigma_{h,jk}^{XZ}\| \leq \omega_0^{1/2} \sigma_{0,kk}^Z,$ we obtain that
		\begin{equation} \nonumber
		\begin{split} 
		\left| \frac{I_2}{\sqrt{\omega_{jl}^X}}\right| &\lesssim \mathcal{M}_1^{X}l^{\alpha_1+1}\eta + \mathcal({M}_1^X)^2l^{(5\alpha_1+4)/2}\eta^2 
		\end{split}
		\end{equation}
		holds with probability greater than $1- 8\exp\{-\tilde c_4 n\min(\eta^2,\eta)\} - 4 \exp\{-\tilde c_4n \\(\{\cM_1^X\}^{-2}l ^{-2(\alpha_1+1)})\}.$
By (\ref{eq_th_XZ}) and the fact that $\sqrt{(\sigma_0^Z)^2\omega_0^X} \leq 1/2 \{(\sigma_0^Z)^2 + \omega_0^X\},$ we obtain that
		\begin{equation}\nonumber
		P\{\left| \frac{I_3}{\sqrt{\omega_{jl}^X }}\right| \geq 2\cM_{X,Z}\sigma_0^Z\eta\}
		\leq c_1\exp\{-c_2n\min(\eta^2,\eta)\}.
		\end{equation}

		Combining the above results, we obtain that there exists positive constants $\rho_5,$ $\rho_6,$ $\tilde c_9$ and $\tilde c_{10}$ such that 
		\begin{small}
		\begin{equation} \nonumber
		\begin{split} 
		&P\left\{\left| \frac{ \widehat{\varrho}_{h,jkl}^{X,Z}  - \varrho_{h,jkl}^{X,Z} }{\sqrt{\omega_{jl}^X }}\right| \geq \rho_5 \cM_{X,Z}l^{\alpha_1+1}\eta + \rho_6\mathcal{M}_{X,Z}^2l^{(5\alpha_1+4)/2} \eta^2\right\}
		\\& \quad \leq \tilde c_{10} \exp\{-\tilde c_9 n \min(\eta^2,\eta)\} + \tilde c_{10} \exp\{-\tilde c_9 \mathcal{M}_{X,Z}^{-2}n l^{-2(\alpha_1+1)} \}.
		\end{split}
		\end{equation}
		\end{small}
		Letting $\eta = \rho_7\sqrt{\frac{\log(pdM_1)}{n}} < 1$ and $\rho_5 + \rho_6\rho_7\mathcal{M}_{X,Z}M_1^{1.5\alpha_1+1}  \eta\leq \rho_8,$ there exist some constants $c_7, c_8 >0$ such that
		\begin{small}
		\begin{equation} \nonumber
		\begin{split} 
		&P\left\{\underset{1 \leq l \leq M_1}{\underset{1 \leq j \leq p, 1 \leq k \leq d}{\max}}\left| \frac{ \widehat{\varrho}_{h,jkl}^{X,Z}  - \varrho_{h,jkl}^{X,Z} }{\sqrt{\omega_{jl}^X }}\right| \geq \rho_7\rho_8 \mathcal{M}_{X,Z}M_1^{\alpha_1+1} \sqrt{\frac{\log(pdM_1)}{n}}\right\} \quad \leq c_7(pdM_1)^{c_8},
		\end{split}  
		\end{equation}
		\end{small}
which implies (\ref{eq_FPCscores_XY_partial}).
$\square$	
\\

\paragraph{Proof of Proposition~\ref{pr_FPCscores_XE}}
To simplify our notation, we will denote $\widehat{\sigma}_{h,jlm}^{X,\epsilon}$ and $\sigma_{h,jlm}^{X,\epsilon}$ by $\widehat{\sigma}_{h,jlm}$ and $\sigma_{h,jlm}$ in subsequent proofs. Recall that $\widehat{\sigma}_{h,jlm} = \langle \widehat \psi_{jl},\langle \widehat{\Sigma}_{h,j}^{X,\epsilon}, \widehat\phi_{m}\rrangle$ and $\sigma_{h,jlm} = \langle\psi_{jl},\langle\Sigma_{h,j}^{X,\epsilon},\phi_{m}\rangle.$ Since we assume $\{\bX_t(\cdot)\}$ and $\{\epsilon_t(\cdot)\}$ are independent processes, $\sigma_{h,jlm} = 0.$

Let $\widehat{r}_{jl} = \widehat{\psi}_{jl} - \psi_{jl},$ $\widehat{w}_{m} = \widehat{\phi}_{m} - \phi_{m}$ and $\widehat{\Delta}_{h,j}^{X,\epsilon} = \widehat{\Sigma}_{h,j}^{X,\epsilon}-\Sigma_{h,j}^{X,\epsilon}.$
\begin{equation} \nonumber
    \begin{split}
        \widehat{\sigma}_{h,jlm}  &= \langle \widehat{r}_{jl}, \langle \widehat{\Sigma}_{h,j}^{X,\epsilon}, \widehat{w}_{m}\rrangle + \left(\langle \widehat{r}_{jl}, \langle \widehat{\Delta}_{h,j}^{X,\epsilon},\phi_{m} \rrangle + \langle \psi_{jl}, \langle \widehat\Delta_{h,j}^{X,\epsilon},\widehat{w}_{m} \rrangle \right) 
        \\& \quad + \langle \psi_{jl},\langle \widehat{\Delta}_{h,j}^{X,\epsilon}, \phi_{m}\rrangle
        \\& = I_1 + I_2 + I_3.
    \end{split}
\end{equation}
Denote $\Omega_{j,\eta}^{X,\epsilon} = \left\{ \|\widehat{\Delta}_{h,j}^{X,\epsilon}\|_\cS \leq (\omega_0^X + \omega_0^\epsilon)\cM_{X,\epsilon}\eta\right\},$ $\Omega_{jj,\eta}^X = \left\{ \|\widehat{\Delta}_{jj}^X\|_\cS  \leq 2\mathcal{M}_1^{X}\omega_0^X\eta \right\},$  $\Omega_{\eta}^Y = \\ \left\{ \|\widehat{\Delta}^Y\|_\cS  \leq 2\mathcal{M}^Y\omega_0^Y\eta \right\}$ and $\Omega_1= \left\{ \|\widehat{\Delta}_{h,j}^{X,\epsilon}\|_\cS \leq (\omega_0^X + \omega_0^\epsilon)\right\}.$ By Theorem~\ref{th_Dev} and Lemma~\ref{lm_guo2019_elementwise}, we have 
\begin{equation} \nonumber
    \begin{split} 
        P((\Omega_{j,\eta}^{X,\epsilon})^C) &\leq c_1\exp\{-c_3n\min(\eta^2,\eta)\},
\\ 
        P((\Omega_{jj,\eta}^X)^C) &\leq 4\exp\{-\tilde{c}_1 n\min(\eta^2,\eta)\},
\\
        P((\Omega_{\eta}^Y)^C) &\leq 4\exp\{-\tilde{c}_1 n\min(\eta^2,\eta)\},
\\
        P((\Omega_1)^C) &\leq c_1\exp\{-c_3n(\cM_{X,\epsilon})^{-2}\}.
    \end{split}
\end{equation}
On the event of $\Omega_1 \cap \Omega_{j,\eta}^{X,\epsilon} \cap \Omega_{jj,\eta}^X \cap \Omega_{\eta}^Y ,$ by Condition~\ref{con_eigen}, (\ref{con_6}), (\ref{con_7}) and Lemma~\ref{lm_X_traceclass}, we obtain that
\begin{equation}\nonumber
    \begin{split}
        \left| \frac{I_1}{\sqrt{\omega_{jl}^X \omega_{m}^Y}}\right| &\leq c_0^{-1}(\alpha_1\alpha_2)^{1/2} l^{\alpha_1/2}m^{\alpha_2/2}\|\widehat{r}_{jl}\|(\|\widehat{\Delta}_{h,j}^{X,\epsilon}\|_\cS  + \|{\Sigma}_{h,j}^{X,\epsilon}\|_\cS)\|\widehat{w}_{m}\|
        \\& \lesssim  (l^{3\alpha_1+2} \vee m^{3\alpha_2+2})\cM_1^X\cM^Y\eta^2,
    \end{split}
\end{equation}
\begin{equation}\nonumber
    \begin{split}
        \left| \frac{I_2}{\sqrt{\omega_{jl}^X \omega_{m}^Y}}\right| & \lesssim  l^{\alpha_1/2}m^{\alpha_2/2}\|\widehat{\Delta}_{h,j}^{X,\epsilon}\|_\cS(l^{\alpha_1+1}\|\widehat{\Delta}_{jj}^X\|_\cS + m^{\alpha_2+1}\|\widehat{\Delta}^Y\|_\cS)
        \\& \lesssim  (l^{2\alpha_1+1} \vee m^{2\alpha_2+1})\cM_{X,\epsilon}\{\mathcal{M}_1^X +\mathcal{M}^Y\}\eta^2,
    \end{split}
\end{equation}
\begin{equation}\nonumber
        \left| \frac{I_3}{\sqrt{\omega_{jl}^X \omega_{m}^Y}}\right| \leq c_0^{-1}(\alpha_1\alpha_2)^{1/2} l^{\alpha_1/2}m^{\alpha_2/2}\|\widehat{\Delta}_{h,j}^{X,\epsilon}\|_\cS \lesssim  (l^{\alpha_1} \vee m^{\alpha_2})\cM_{X,\epsilon}\eta.
\end{equation}

Combining the above results, we obtain that there exists positive constants $\rho_9,$ $\rho_{10},$ $\tilde c_{11}$ and $\tilde c_{12}$ such that 
\begin{equation} \nonumber
    \begin{split} 
        &P\left\{\left| \frac{ \widehat{\sigma}_{h,jlm}}{\sqrt{\omega_{jl}^X \omega_{m}^Y}}\right| \geq \rho_9 \mathcal{M}_{X, \epsilon}(l ^{\alpha_1}\vee m^{\alpha_2})\eta + \rho_{10}\mathcal{M}_{X, \epsilon}(\mathcal{M}_1^X +\mathcal{M}^Y)(l ^{3\alpha_1+2}\vee m^{3\alpha_2+2})\eta^2\right\}
        \\& \quad \leq \tilde c_{12} \exp\{-\tilde c_{11} n \min(\eta^2,\eta)\} + \tilde c_{12} \exp\{-\tilde c_{11} \mathcal{M}_{X,\epsilon}^{-2}n \}.
    \end{split}
\end{equation}
Letting $\eta = \rho_{11}\sqrt{\frac{\log(pM_1M_2)}{n}} < 1$ and $\rho_9 + \rho_{10}\rho_{11}\{\mathcal{M}_1^X + \mathcal{M}^Y\}(M_1^{2\alpha_1+2} \vee M_2^{2\alpha_2+2})\eta\leq \rho_{12},$ there exist some constants $c_{9},c_{10} > 0$ such that
\begin{small}
\begin{equation} \nonumber
\begin{split}  
&P\left\{\underset{1 \leq l \leq M_1, 1 \leq m \leq M_2}{\underset{1 \leq j \leq p}{\max}}\left| \frac{ \widehat{\sigma}_{h,jlm} - \sigma_{h,jlm}}{\sqrt{\omega_{jl}^X \omega_{m}^Y}}\right| \geq \rho_{11}\rho_{12} (\cM_1^X + \cM^\epsilon)(M_1^{\alpha_1} \vee M_2^{\alpha_2})\sqrt{\frac{\log(pM_1M_2)}{n}}\right\} 
\\& \leq  c_{9}(pM_1M_2)^{c_{10}},
\end{split}  
\end{equation}
\end{small}
which completes the proof. 
$\square$
\\
\paragraph{Proof of Proposition~\ref{pr_FPCscores_XE_partial}}
Recall that $\widehat\varrho_{h,jl}^{X,\epsilon} = \frac{1}{n-h}\sum_{t=1}^{n-h}\widehat{\zeta}_{tjl}\epsilon_{t+h} $ and $\varrho_{h,jl}^{X,\epsilon} =\cov({\zeta}_{tjl},\epsilon_{t+h}).$ Let $\widehat{r}_{jl} = \widehat{\psi}_{jl} - \psi_{jl}$ and $\widehat{\Delta}_{h,j}^{X,\epsilon} = \widehat{\Sigma}_{h,j}^{X,\epsilon}-\Sigma_{h,j}^{X,\epsilon}.$ We have 
		\begin{equation} \nonumber
		\begin{split}
		\widehat{\varrho}_{h,jl}^{X,\epsilon}  - \varrho_{h,jl}^{X,\epsilon}  &= \langle \widehat{r}_{jl},  \widehat{\Delta}_{h,j}^{X,\epsilon}\rangle  + \langle \psi_{jl}, \widehat{\Delta}_{h,j}^{X,\epsilon} \rangle
		\\& = I_1 + I_2. 
		\end{split}
		\end{equation}
Let $\Omega_{j,\eta}^{X,\epsilon} = \left\{ \|\widehat{\Delta}_{h,j}^{X,\epsilon}\| \leq (\omega_0^X + (\sigma_0^\epsilon)^2)\cM_{X,\epsilon}\eta\right\}$and $\Omega_{jj,\eta}^X = \left\{ \|\widehat{\Delta}_{jj}^X\|_\cS  \leq 2\mathcal{M}_1^{X}\omega_0^X\eta \right\}.$ By (\ref{eq_Dev_XZ}) and Lemma~\ref{lm_guo2019_elementwise}, we have 
		\begin{equation} \nonumber
		\begin{split} 
		P((\Omega_{j,\eta}^{X,\epsilon})^C) &\leq c_1\exp\{-c_3n\min(\eta^2,\eta)\},
\\
		P((\Omega_{jj,\eta}^X)^C) &\leq 4\exp\{-\tilde{c}_1 n\min(\eta^2,\eta)\}.
		\end{split}
		\end{equation}
On the event of $ \Omega_{\eta,j}^{X,\epsilon} \cap \Omega_{jj,\eta}^X  ,$ by Condition~\ref{con_eigen}, (\ref{con_6}) and Lemma~\ref{lm_X_traceclass}, we obtain that
		\begin{equation}\nonumber
		\begin{split}
		\left| \frac{I_1}{\sqrt{\omega_{jl}^X}}\right| &\lesssim  l^{\alpha_1/2}\|\widehat{\Delta}_{h,j}^{X,\epsilon}\| \|\widehat{r}_{jl}\|  \lesssim l^{3\alpha_1/2+1}\|\widehat{\Delta}_{h,j}^{X,\epsilon}\|\|\widehat{\Delta}_{jj}^X\|_\cS 
		\\& \lesssim  l^{3\alpha_1/2+1} \cM_{X,\epsilon}\cM_1^{X} \eta^2.
		\end{split}
		\end{equation}
By (\ref{eq_th_XZ}) and the fact that $\sqrt{(\sigma_0^\epsilon)^2\omega_0^X} \leq 1/2 \{(\sigma_0^\epsilon)^2 + \omega_0^X\},$ we obtain that
		\begin{equation}\nonumber
		P\{\left| \frac{I_2}{\sqrt{\omega_{jl}^X }}\right| \geq 2\cM_{X,\epsilon}\sigma_0^\epsilon\eta\}
		\leq c_1\exp\{-c_2n\min(\eta^2,\eta)\}.
		\end{equation}

		Combining the above results, we obtain that there exists positive constants $\rho_{13},$ $\rho_{14},$ $\tilde c_{13}$ and $\tilde c_{14}$ such that 
		\begin{small}
		\begin{equation} \nonumber
		\begin{split} 
		&P\left\{\left| \frac{ \widehat{\varrho}_{h,jl}^{X,\epsilon}  - \varrho_{h,jl}^{X,\epsilon} }{\sqrt{\omega_{jl}^X }}\right| \geq \rho_{13} \cM_{X,\epsilon}\eta + \rho_{14}l^{3\alpha_1/2+1}\mathcal{M}_{X,\epsilon}\cM_1^X \eta^2\right\} \leq \tilde c_{14} \exp\{-\tilde c_{13} n \min(\eta^2,\eta)\}.
		\end{split}
		\end{equation}
		\end{small}
		Letting $\eta = \rho_{15}\sqrt{\frac{\log(pM_1)}{n}} < 1$ and $\rho_{13} + \rho_{14}\rho_{15}\mathcal{M}_1^X M_1^{3\alpha_1/2+1}  \eta\leq \rho_{16},$ there exist some constants $c_{11}, c_{12} >0$ such that
		\begin{small}
		\begin{equation} \nonumber
		\begin{split} 
		&P\left\{\underset{1 \leq l \leq M_1}{\underset{1 \leq j \leq p}{\max}}\left| \frac{ \widehat{\varrho}_{h,jl}^{X,\epsilon}  - \varrho_{h,jl}^{X,\epsilon} }{\sqrt{\omega_{jl}^X }}\right| \geq \rho_{15}\rho_{16} \mathcal{M}_{X,\epsilon} \sqrt{\frac{\log(pM_1)}{n}}\right\} \quad \leq c_{11}(pM_1)^{c_{12}},
		\end{split}  
		\end{equation}
		\end{small}
which implies (\ref{eq_FPCscores_XE_partial}).
$\square$

\subsection{Technical lemmas and their proofs} \label{ap.lemma.secA3}
\begin{lemma}
\label{lemma_compare}
The non-functional version of our proposed cross-spectral stability measure satisfies
\begin{equation} \nonumber
    \mathop{\text{ess sup}}\limits_{\theta \in [-\pi,\pi],\bnu_1\in \mathbb{R}_0^p,\bnu_2\in \mathbb{R}_0^d}\frac{\left|\bnu_1^{\T} \bbf_\theta^{X,Y}\bnu_2\right|}{\sqrt{\bnu_1^{\T}\bnu_1}\sqrt{\bnu_2^{\T}\bnu_2}} \leq  \widetilde{\mathcal{M}}^{X,Y},
\end{equation}
where $\widetilde{\mathcal{M}}^{X,Y}$ is defined in (\ref{eq_basuM}).
\end{lemma}
{\bf Proof.} For any fixed $\theta \in [-\pi,\pi],$ we perform singular value decomposition on $\bbf_\theta^{X,Y} = \bU\bD\bV^{\T},$ where $\bD$ is a diagonal matrix with singular values $\{\sigma_i\}$ of $\bbf_\theta^{X,Y}$ on the diagonal. Then
\begin{eqnarray*}
    \mathop{\max}\limits_{\bnu_1\in \widetilde{\mathbb{R}}_0^p,\bnu_2\in \widetilde{\mathbb{R}}_0^d}\frac{\left|\bnu_1^{\T} \bbf_\theta^{X,Y}\bnu_2\right|}{\sqrt{\bnu_1^{\T}\bnu_1}\sqrt{\bnu_2^{\T}\bnu_2}}  &=& \mathop{\max}\limits_{\bx\in \widetilde{\mathbb{R}}_0^p,\by\in \widetilde{\mathbb{R}}_0^d} \frac{\left|\bx^{\T} \bD\by\right|}{\sqrt{\bx^{\T}\bx}\sqrt{\by^{\T}\by}} \quad (\bx= \bU^{\T}\bnu_1, \by= \bV^{\T}\bnu_2)
    \\&=& \mathop{\max}\limits_{\bx\in \widetilde{\mathbb{R}}_0^p,\by\in \widetilde{\mathbb{R}}_0^d}\frac{\sum x_iy_i\sigma_i}{\sqrt{\bx^{\T}\bx}\sqrt{\by^{\T}\by}}
    \\&\leq& \mathop{\max}\limits_{\bx\in \widetilde{\mathbb{R}}_0^p,\by\in \widetilde{\mathbb{R}}_0^d}\frac{\sqrt{\sum x_i^2\sum(y_i\sigma_i)^2}}{\sqrt{\sum x_i^2}\sqrt{\sum y_i^2}} = \mathop{\max}\limits_{\by\in \widetilde{\mathbb{R}}_0^d}\sqrt{\frac{{\sum(y_i\sigma_i)^2}}{{\sum y_i^2}} }
    \\&\leq& \max (\sigma_i)= \mathop{\max}\limits_{\bnu\in \widetilde{\mathbb{R}}_0^d}\sqrt{\frac{\bnu^{\T} \{\bbf^{X,Y}_\theta\}^* \bbf^{X,Y}_\theta\bnu}{\bnu^{\T}\bnu}}.
\end{eqnarray*}
This holds almost everywhere for $\theta \in [-\pi,\pi],$ which completes our proof.
$\square$

\begin{lemma}  \label{lm_X_traceclass}
Suppose that Conditions~\ref{con_sub_coefficient} and \ref{con_e} hold, then 
$\omega_0^{X} =O(1).$
\end{lemma}
{\bf Proof.} Recall that $\bX_t(u) = \sum\limits_{l=0}^{\infty}\int \bA_l(u,v)\bvarepsilon_{t-l}(v)dv$ and $\bvarepsilon_{t}(\cdot)$'s are i.i.d. mean-zero functional processes. 
Let $\bA_{l,j}$ denote the $j$-th row of $\bA_l$. Then
\begin{equation} \nonumber
\begin{split}
         & \quad \max_{1\leq j\leq p}\int\Sigma_{0,jj}^{X}(u,u)du 
\\&= \max_{ j}\int\mathbb{E}\left\{X_{tj}(u)X_{tj}(u)\right\}du
    \\& = \max_{ j}\sum_{l=0}^\infty \int\int \bA_{l,j}(u,v)\bSigma_0^{\varepsilon}(v,v)\{\bA_{l,j}(u,v)\}^{\T} dvdu
    \\& = \max_{ j}\sum_{l=0}^\infty \sum_{k,k' = 1}^p \int\int \Sigma_{0,kk'}^\varepsilon(v,v) A_{l,jk}(u,v)A_{l,jk'}(u,v)dvdu 
    \\ & \leq  \max_{ j}\sum_{l=0}^\infty \sum_{k,k' = 1}^p\int \Sigma_{0,kk'}^\varepsilon(v,v)dv \int\int A_{l,jk}(u,v)A_{l,jk'}(u,v)dvdu
    \\ & \leq \max_j \left(\max_{k,k'} \int \Sigma_{0,kk'}^\varepsilon(v,v)dv\right) \left( \sum_{l=0}^\infty \sum_{k,k' = 1}^p\int \int A_{l,jk}(u,v)A_{l,jk'}(u,v)dvdu\right)
    \\&  \leq  \omega_0^\varepsilon \max_{j}\sum_{l=0}^\infty \left(\sum_{k= 1}^p\sum_{k'=1}^p \sqrt{\int \int (A_{l,jk}(u,v))^2dvdu \int \int (A_{l,jk'}(u,v))^2dvdu }\right) \quad 
    \\&=  \omega_0^\varepsilon \max_{j}\sum_{l=0}^\infty (\sum_{k=1}^p \|A_{l,jk}\|_{\cS})^2 = \omega_0^\varepsilon \sum_{l=0}\|\bA_l\|_\infty^2 \leq \omega_0^\varepsilon \{\sum_{l=0}\|\bA_l\|_\infty\}^2 = O(1),
\end{split}
\end{equation}
which completes our proof.
$\square$
\\

Before presenting Lemma~{\ref{lm_hwinequality}}, we define sub-Gaussian distribution and sub-Gaussian norm as follows. A centered random variable $x$ with variance proxy $\sigma^2$ is sub-Gaussian if for any $t >0,$ $P(|x|>t) \leq 2 \exp (-t^2/(2 \sigma^2)).$ The sub-Gaussian norm of $x$ is defined by
$
\|x\|_{\psi_2} = \inf \{K>0 : \mathbb{E}\exp(x^2/K^2)  \leq 2\}.$

\begin{lemma} \label{lm_hwinequality}
Let $\bx = (x_1, \dots, x_n) \in \mathbb{R}^n$ be a random vector with independent mean zero sub-Gaussian coordinates. Without loss of generality, we assume that $\mathbb{E}x_i^2 = 1$ for $i = 1, \dots, n.$ Let $\bA$ be an $ n \times n$ matrix. Then there exists some universal constant $c>0$ such that for any given $\eta>0,$ 
\begin{equation} \label{lm_hw}
     P\left(|\bx^{\T} \bA \bx - \mathbb{E}\bx^{\T} \bA \bx | \geq \|\bA\|\eta\right) \leq 2 \exp \left\{ -c \min \left(\frac{ \eta^2}{\text{rank}(\bA)},\eta\right)\right\}.
\end{equation}
\end{lemma}
\textbf{Proof}. It follows from Theorem~1.1 of \cite{rudelson2013} and $\|x_i\|_{\psi_2}=1$
for $i = 1, \dots, n,$ that there exists a constant $c>0$ such that
\begin{equation} \nonumber
     P\left(|\bx^{\T} \bA \bx - \mathbb{E}\bx^{\T} \bA \bx | \geq t\right) \leq 2 \exp \left\{ -c \min \left(\frac{ t^2}{\|\bA\|_\tF^2},\frac{ t}{\|\bA\|}\right)\right\}.
\end{equation}
By $\|\bA\|_\tF \leq \sqrt{\text{rank}(\bA)}\|\bA\|$ and letting $ t  = \eta\|\bA\|,$ we obtain 
(\ref{lm_hw}). $\square$
\\
\begin{lemma} \label{lm_subG definition}
Suppose that sub-Gaussian process $\{\varepsilon_{tj}(\cdot)\}_{t \in \mathbb{Z}}$ follows Definition \ref{def_subGaussian}. Under Karhunen-Lo\`eve expansion $\varepsilon_{tj}(\cdot) = \sum_{l=1}^\infty \xi_{tjl}\phi_{jl}(\cdot)  = \sum_{l=1}^\infty \sqrt{\omega_{jl}^{\varepsilon}}a_{tjl} \phi_{jl}(\cdot)$ with $\mathbb{E}(a_{tjl}) = 0$ and $\mathbb{E}(a_{tjl}^2) = 1$ for $t \in \eZ$ and $j=1\dots,p,$ $a_{tjl}$ follows sub-Gaussian distribution with $\|a_{tjl}\|_{\psi_2}=1,$ that is for all $\eta >0,$ $t \in \eZ,$ $j=1,\dots,p$ and $ l \geq 1,$ 
\begin{equation} \nonumber
    P[|a_{tjl}| > \eta] \leq 2 \exp (-\eta^2/2).
\end{equation}
\end{lemma}
\textbf{Proof}. 
By Definition~\ref{def_subGaussian}, for all $ x \in \mathbb{H},$ $\mathbb{E}\{e^{\langle x, X \rangle}\} \leq e^{\alpha^2\langle x,\Sigma_0(x)\rangle/2}.$ Combining with the choice of $x=c\phi_{jl}(\cdot)$ for $c > 0$ and orthonormality of $\{\phi_{jl}(\cdot)\}$ yields
\begin{equation} \nonumber
    \mathbb{E}(e^{c\sqrt{\omega_{jl}^{\varepsilon}} a_{tjl}}) \leq e^{\alpha^2 c^2 {\omega_{jl}^{\bvar}} /2}.
\end{equation}
Without loss of generality, we assume $\alpha = 1.$ By Markov's inequality and the above result, we have that for all $c > 0,$
\begin{equation} \nonumber
    P(a_{tjl}> \eta) \leq P(e^{c\sqrt{\omega_{jl}^{\varepsilon}}a_{tjl}} > e^{c\sqrt{\omega_{jl}^{\bvar}}\eta}) \leq \frac{\mathbb{E}(e^{c\sqrt{\omega_{jl}^{\varepsilon}}a_{tjl}})}{e^{c\sqrt{\omega_{jl}^{\varepsilon}}\eta}} \leq e^{c^2{\omega_{jl}^{\varepsilon}}/2 - c\sqrt{\omega_{jl}^{\varepsilon}}\eta}.
\end{equation}
Choosing $c =\eta/\sqrt{\omega_{jl}^{\varepsilon}},$ we have 
$
    P(a_{tjl}> \eta) \leq e^{- \frac{\eta^2}{2}}.
$
In the same manner with the choice of $x=-c\phi_{jl}(\cdot)$ for $c > 0,$ we can prove $P(a_{tjl}< - \eta) \leq e^{- \frac{\eta^2}{2}}.$ Combining the above results,
$
P[|a_{tjl}| > \eta] = P(a_{tjl}> \eta) + P(a_{tjl}< - \eta) \leq  2e^{- \frac{\eta^2}{2}}
$
which completes the proof. $\square$
\\

Before presenting Lemma~\ref{lm_iidsub} below, we give some definitions: 
\\
(i) Suppose that $\be = (e_1, \dots, e_N)^{\T} \in \mathbb{H}^{N}$ is formed by $N$ independent mean zero sub-Gaussian processes with $e_{i}(\cdot)  = \sum_{l=1}^\infty \sqrt{\omega_{il}^{e}}a_{il} \phi_{il}(\cdot)$ under the Karhunen-Lo\`eve expansion. Define 
$\bvarphi_{M,i} = \big(\sqrt{\omega_{i1}^{e}}\phi_{i1},\dots,\sqrt{\omega_{iM}^{e}}\phi_{iM}\big)^{\T}.$ 
\\
(ii) Suppose $\bK=(K_{ij})_{N \times N}$ with each $K_{ij} \in {\mathbb S}$. 
For any nonempty subset $G \subset \mathbb{Z}_{+}$ = \{1,2,\dots\} with $ |G|<\infty,$ write $G=\{g_1, \dots, g_{|G|}\}$ with $g_1 < \dots < g_{|G|}$ and $\bphi_{G,i}=(\phi_{ig_1}, \dots, \phi_{ig_{|G|}})^{\T}$ for each $i = 1,\dots, N$.
Let $\bvarPhi_G=\text{diag}(\bphi_{G,1}^{\T}, \dots, \bphi_{G,N}^{\T}),$ 
then we define $$\text{rank}(\bK) = \sup_{G\subset \mathbb{Z}_{+}, |G| < \infty} \text{rank}\left( \int\int\bvarPhi_G^{\T}(u)\bK(u,v)\bvarPhi_G(v)dudv\right).$$

\begin{condition} \label{con_Pi}
Let
$\bPi_M = \int\int\bTheta_M^{\T}(u)\bK(u,v)\bTheta_M(v)dudv$  with $\bTheta_M=\text{diag}(\bvarphi_{M,1}^{\T}, \\\dots, \bvarphi_{M,N}^{\T})$ and $\bK=(K_{ij})_{N \times N}$ with each $K_{ij} \in {\mathbb S}.$ It satisfies that
$\|\bPi_M\| \leq b_M$ and $\lim_{M \to \infty}b_M = b.$
\end{condition}

\begin{lemma}\label{lm_iidsub}
Suppose that $\max_{1 \leq i \leq N}\int_\cU \Sigma_{ii}^{e}(u,u)du < \infty$ and $\bK$ satisfies Condition~\ref{con_Pi}. Then, there exists some universal constant $c>0$ such that for any given $\eta>0,$ 
\begin{equation} \label{subG_2}
    P\left(|\langle\be, \bK(\be)\rangle - \mathbb{E}\langle\be, \bK(\be)\rangle | \geq b\eta\right) \leq 2 \exp \left\{ -c \min \left(\frac{\eta^2}{\text{rank}(\bK)},\eta\right)\right\}.
\end{equation}
\end{lemma}
\textbf{Proof}.
We organize our proof as follows: First, we truncate $e_{i}(\cdot)$ to $M$-dimensional process $e_{M,i}(\cdot)=\sum_{l=1}^M \sqrt{\omega_{il}^{e}}a_{il} \phi_{il}(\cdot),$ then apply Hanson-Wright inequality in Lemma~\ref{lm_hwinequality} and finally show that the inequality still hold under the infinite-dimensional setting. 

Rewrite $\be_M=(e_{M,1}, \dots, e_{M,N})^{\T}$ with $e_{M,i}= \ba_{M,i}^{\T}\bvarphi_{M,i}$ and $\ba_{M,i} = (a_{i1},\dots,a_{iM})^{\T}.$ 
Let $\ba_M = (\ba_{M,1}^{\T}, \dots, \ba_{M,N}^{\T})^{\T} \in \mathbb{R}^{NM},$ then we have
$ \langle\be_{M}, \bK(\be_{M})\rangle = \ba_{M}^{\T} \bPi_M \ba_{M}.$ 
By Lemma~\ref{lm_subG definition}, elements in $\ba_{M} \in \mathbb{R}^{NM}$ are i.i.d. sub-Gaussian with $\mathbb{E}(a_{il}) = 0$ and $\mathbb{E}(a_{il}^2) = 1.$ Combining this with Lemma~\ref{lm_hwinequality} yields
\begin{equation}\label{HW_M}
\begin{split}
    &\quad P\left(|\langle\be_{M}, \bK(\be_{M})\rangle - \mathbb{E}\langle\be_{M}, \bK(\be_{M})\rangle | \geq b_M\eta \right)   
    \\& 
    \leq P\left(|\ba_{M}^{\T} \bPi_M \ba_{M} - \mathbb{E}\ba_{M}^{\T} \bPi_M \ba_{M} | \geq \|\bPi_M\|\eta\right) 
    \\& \leq 
    2 \exp \left\{ -c \min \left(\frac{\eta^2}{\text{rank}(\bPi_M)},\eta\right)\right\}.
\end{split}
\end{equation} 

It follows from Lemma~\ref{lm_e} that $ \langle\be_{M}, \bK(\be_{M})\rangle$ converges in probability to $\langle\be, \bK(\be)\rangle$ and $\lim_{M \to \infty}\mathbb{E}\langle\be_{M}, \bK(\be_{M})\rangle=  \mathbb{E}\langle\be, \bK(\be)\rangle.$ These results together with 
Condition~\ref{con_Pi} imply that 
\begin{equation} \nonumber
    \langle\be_{M}, \bK(\be_{M})\rangle - \mathbb{E}\langle\be_{M}, \bK(\be_{M})\rangle - b_M\eta
\end{equation}
converges in distribution to 
\begin{equation} \nonumber
    \langle\be, \bK(\be)\rangle - \mathbb{E}\langle\be, \bK(\be)\rangle - b\eta.
\end{equation}
Finally, by the fact that $\text{rank}(\bPi_M) \leq  \text{rank}(\bK)$ and taking $M \rightarrow \infty$ on both sides of (\ref{HW_M}), we obtain (\ref{subG_2}), which completes the proof. 
$\square$
\\

\begin{lemma} \label{lm_e}
Under the same assumption and notation in Lemma~\ref{lm_iidsub} and its proof, we have
\begin{equation} \label{lm_e2} 
    \lim_{M \to \infty} \mathbb{E}\left\{\|\be_{M} - \be\|^2\right\} =  0
\end{equation}
and 
\begin{equation} \label{lm_eke}
    \lim_{M \to \infty}\mathbb{E}\langle\be_{M}, \bK(\be_{M})\rangle=  \mathbb{E}\langle\be, \bK(\be)\rangle.
\end{equation}
\end{lemma}
\textbf{Proof}. Since $\|\be_M - \be\|^2 = \sum_{i = 1}^N\|e_{M,i} - e_{i}\|^2 = \sum_{i = 1}^N\|\sum_{l=M+1}^\infty \sqrt{\omega_{il}^{e}}a_{il} \phi_{il}\|^2,$ it suffices to show  $\lim_{M \to \infty} \mathbb{E}\left\{\|\sum_{l=M+1}^\infty \sqrt{\omega_{il}^{e}}a_{il} \phi_{il}\|^2\right\} =  0.$ By $\mathbb{E}(a_{il}a_{il'}) = 1\{l=l'\}$  and the orthonormality of $\{\phi_{il}\},$ we have 
\begin{equation} \nonumber
    \mathbb{E}\left\{\int\left(\sum_{l=M+1}^\infty \sqrt{\omega_{il}^{\varepsilon}}a_{il} \phi_{il}(u)\right)^2du\right\} = \sum_{l=M+1}^\infty\omega_{il}^{\varepsilon}.
\end{equation}
This together with Condition~\ref{con_e} implies that above goes to zero as $M \to \infty,$ which completes the proof of (\ref{lm_e2}).

By triangle inequality, we have 
\begin{equation}  \label{ieq_quad}
    |\mathbb{E}\langle\be_{M}, \bK(\be_{M})\rangle - \mathbb{E}\langle\be, \bK(\be)\rangle| 
     \leq |\mathbb{E}\langle\be_{M}, \bK(\be_{M} - \be)\rangle| + |\mathbb{E}\langle(\be_{M}-\be), \bK(\be)\rangle|.
\end{equation}
By Jensen's inequality and Lemma~\ref{lm_curveinequality}, we have
\begin{equation} \nonumber
\begin{split}
    |\mathbb{E}\langle\be_{M}, \bK(\be_{M} - \be)\rangle|^2&\leq \|\bK\|_\tF^2\mathbb{E}(\|\be_{M}\|^2)\mathbb{E}(\|\be_{M} - \be\|^2),
    \\
    |\mathbb{E}\langle(\be_{M}-\be), \bK(\be)\rangle|^2&\leq \|\bK\|_\tF^2\mathbb{E}(\|\be\|^2)\mathbb{E}(\|\be_{M} - \be\|^2).
    \end{split}
\end{equation}
From (\ref{lm_e2}), we have $\lim_{M \to \infty} \mathbb{E}\left\{\|\be_{M} - \be\|^2\right\} =  0$ and $\lim_{M \to \infty}\mathbb{E}\{\|\be_{M}\|^2\} = \mathbb{E}\{\|\be\|^2\}.$ Combining these with $\mathbb{E}(\|\be\|^2) \leq {N \max_{1 \leq i \leq N}\int_\cU \Sigma_{ii}^{e}(u,u)du} < \infty$ and $\|\bK\|_\tF < \infty $ implies the right side of (\ref{ieq_quad}) goes to zero when $M \to \infty,$ which completes the proof of (\ref{lm_eke}). $\square$
\\
\begin{lemma} \label{lm_middle matrix_ML}
Suppose Conditions~\ref{con_stability}, \ref{con_sub_coefficient} and \ref{con_e} hold for stationary sub-Gaussian process  $\{\bX_t(\cdot)\}_{t\in\mathbb{Z}}.$  Let $\bX_{M,L,t}(u) = \sum\limits_{l=0}^{L}\bA_l(\bvarepsilon_{M,t-l}).$  
Then, for any $\bPhi_1\in \mathbb{H}_0^p$ with $\|\bPhi_1\|_0 \leq k$ and $k=1,\dots,p,$
\begin{equation} \nonumber
    \lim_{M \to \infty} \cM(\bbf_{M,L}^X, \bPhi_1) =\cM(\bbf_L^X, \bPhi_1).
\end{equation}
\end{lemma}
\textbf{Proof}. By the definitions of $\cM(\bbf_{M,L}^X, \bPhi)$ and $\bbf^X_{M,L,\theta}(\bPhi)$ in the proof of Theorem~\ref{th_XY} in Appendix~\ref{ap.th.secA1}, we have
\begin{equation} \nonumber
    \begin{split}
        & \lim_{M \to \infty} |\cM(\bbf_{M,L}^X, \bPhi_1) - \cM(\bbf_{L}^X, \bPhi_1)| 
        \\&=  2\pi  \lim_{M \to \infty}  \left|\mathop{\text{ess sup}}\limits_{\theta \in [-\pi,\pi]} |\langle\bPhi_1,\bbf^X_{M,L,\theta}(\bPhi_1)\rangle| -   \mathop{\text{ess sup}}\limits_{\theta \in [-\pi,\pi]} |\langle\bPhi_1,\bbf^X_{L,\theta}(\bPhi_1)\rangle| \right|
        \\&  \leq 2\pi  \lim_{M \to \infty}  \mathop{\text{ess sup}}\limits_{\theta \in [-\pi,\pi]}\left| |\langle\bPhi_1,\bbf^X_{M,L,\theta}(\bPhi_1)\rangle| -  |\langle\bPhi_1,\bbf^X_{L,\theta}(\bPhi_1)\rangle| \right|
        \\& \leq \|\bPhi_1\|^2 \lim_{M \to \infty} \left\|\sum_{h \in \mathbb{Z}}(\bSigma_{M,L,h}^X - \bSigma_{L,h}^X  )\right\|_\tF  \quad \text{(by Lemma~\ref{lm_curveinequality} and $|\exp(-ih\theta)|=1$)}
        \\& \leq \|\bPhi_1\|^2 \lim_{M \to \infty} \sum_{h \in \mathbb{Z}}\left\|\bSigma_{M,L,h}^X - \bSigma_{L,h}^X\right\|_\tF.
    \end{split}
\end{equation}
Provided that $\|\bPhi_1\|^2 <\infty,$ it suffices to prove that $ \sum_{h= -\infty}^{\infty}\left\|\bSigma_{M,L,h}^X - \bSigma_{L,h}^X\right\|_\tF < \infty$ and $\lim_{M \to \infty}\left\|\bSigma_{M,L,h}^X - \bSigma_{L,h}^X\right\|_\tF = 0.$

By triangle inequality and Lemma~\ref{lm_sigma_sub}, we obtain that
\begin{equation} \nonumber
    \begin{split}
        \sum_{h= -\infty}^{\infty}\left\|\bSigma_{M,L,h}^X - \bSigma_{L,h}^X\right\|_\tF \leq 
        \sum_{h= -\infty}^{\infty} \|\bSigma_{M,L,h}^X\|_\tF + \sum_{h= -\infty}^{\infty}\|\bSigma_{L,h}^X\|_\tF < \infty.
    \end{split}
\end{equation}

We next prove $\lim_{M \to \infty}\left\|\bSigma_{M,L,h}^X - \bSigma_{L,h}^X\right\|_\tF = 0.$ 
Write
\begin{equation}\nonumber
    \begin{split}
       \bSigma_{M,L,h}^{X}(u,v)  &= \mathbb{E}\left\{\bX_{M,L,t-h}(u)\bX_{M,L,t}^{\T}(v)\right\} \\&= \sum_{l=0}^{L-h} \int \bA_{l+h}(u,u')\bSigma_0^{\varepsilon_M}(u',v')\{\bA_{l}(v,v')\}^{\T} du'dv',
       \\ \bSigma_{L,h}^{X}(u,v)  &= \mathbb{E}\left\{\bX_{L,t-h}(u)\bX_{L,t}^{\T}(v)\right\}\\&= \sum_{l=0}^{L-h} \int \bA_{l+h}(u,u')\bSigma_0^{\varepsilon}(u',v')\{\bA_{l}(v,v')\}^{\T} du'dv'.
    \end{split}
\end{equation}
Then,
\begin{equation}\nonumber
    \begin{split}
        &\lim_{M \to \infty}\left\|\bSigma_{M,L,h}^X - \bSigma_{L,h}^X\right\|_\tF 
        \\= & \lim_{M \to \infty}\left\|\sum_{l=0}^{L-h} \int \bA_{l+h}(u,u')\{\bSigma_0^{\varepsilon_M}(u',v')-\bSigma_0^{\varepsilon}(u',v')\}\{\bA_{l}(v,v')\}^{\T} du'dv'\right\|_\tF
        \\ \leq & \sum_{l=0}^{L-h}\| \bA_l\|_\tF\|\bA_{l+h}\|_\tF\lim_{M \to \infty}\|\bSigma_0^{\varepsilon_M}-\bSigma_0^{\varepsilon}\|_\tF \quad (\text{by Lemma~\ref{lm_curveinequality}})
        \\ \leq & \sum_{l=0}^{L-h}\| \bA_l\|_\tF\|\bA_{l+h}\|_\tF\lim_{M \to \infty}\left\{\sum_{j,k}\|\Sigma_{h,jk}^{\varepsilon_M}-\Sigma_{h,jk}^{\varepsilon}\|_\cS^2\right\}^{1/2}
        \\ \leq & \sum_{l=0}^{L-h}\| \bA_l\|_\tF\|\bA_{l+h}\|_\tF\lim_{M \to \infty}\sum_{j,k}\|\Sigma_{0,jk}^{\varepsilon_M}-\Sigma_{0,jk}^{\varepsilon}\|_\cS
        \\ = & 0 \quad (\text{by Lemmas~\ref{lm_sigma_sub} and \ref{lm_sigma_e}})
    \end{split}
\end{equation}
which completes the proof. $\square$
\\
\begin{lemma} \label{lm_Q}
Suppose that conditions in Lemma~\ref{lm_middle matrix_ML} hold.
For any $\bPhi_1\in \mathbb{H}_0^p$ with $\|\bPhi_1\|_0 \leq k$ and $k=1,\dots,p,$ define $\bY = (\langle\bPhi_1,\bX_1\rangle,\dots,\langle\bPhi_1,\bX_n\rangle )^{\T}.$ Then
\begin{equation} \nonumber
    \|\var(\bY)\| \leq \cM(\bbf^X, \bPhi_1) \leq \cM_k^X\langle\bPhi_1,\bSigma_0^X(\bPhi_1)\rangle.
\end{equation}
\end{lemma}
\textbf{Proof}. The proof follows from the proof of Theorem~1 in \cite{guo2019} and hence the proof is omitted here. $\square$
\\
\begin{lemma} \label{lm_Y}
Suppose that conditions in Lemma~\ref{lm_middle matrix_ML} hold.  Let $\bX_{L,t}(u) = \sum\limits_{l=0}^{L}\bA_l(\bvarepsilon_{t-l}).$ For any $\bPhi_1\in \mathbb{H}_0^p$ with $\|\bPhi_1\|_0 \leq k$ $(k=1,\dots,p),$ define $\bY_L = (\langle\bPhi_1,\bX_{L,1}\rangle,\dots,\langle\bPhi_1,\bX_{L,n}\rangle)^{\T}$ and $\bY = (\langle\bPhi_1,\bX_1\rangle,\dots,\langle\bPhi_1,\bX_n\rangle )^{\T},$ then
\begin{equation} \label{lm_Y_l2}
    \lim_{L \to \infty} \mathbb{E}\left\{\|\bY_L-\bY\|^2\right\} = 0
\end{equation}
and 
\begin{equation}\label{lm_YY_expectation}
    \lim_{L \to \infty} \mathbb{E}\left[\bY_L^{\T}\bY_L\right] = \mathbb{E}\left[\bY^{\T}\bY\right].
\end{equation}
\end{lemma}
\textbf{Proof of (\ref{lm_Y_l2})}. By definitions of $\bY_L$ and $\bY,$ we have that
\begin{equation} \nonumber
    \begin{split}
        \mathbb{E}\left\{\|\bY_L-\bY\|^2\right\} =\sum_{t=1}^n \mathbb{E}\left\{|\langle\bPhi_1,\bX_{L,t}-\bX_t\rangle|^2\right\}
    \end{split}
\end{equation}
By Lemma~\ref{lm_curveinequality}, we have $\mathbb{E}\left\{|\langle\bPhi_1,\bX_{L,t}-\bX_t\rangle|^2\right\} \leq \|\bPhi_1\|^2 \mathbb{E}\{\|\bX_{L,t}-\bX_t\|^2\}.$ With the fact $\|\bPhi_1\|^2 < \infty,$ it suffices to prove that $\lim_{L \to \infty} \mathbb{E}\{\|\bX_{L,t}-\bX_t\|^2\} = 0$ for $ t=1,\dots,n.$ By Lemma~\ref{lm_sigma_e}, we have $\mathbb{E}(\|\bvarepsilon_{t-l}\|) \leq \sqrt{p\omega_0^{\bvarepsilon}}.$ This together with Lemma~\ref{lm_curveinequality} implies that
\begin{equation}\nonumber
    \begin{split}
        \mathbb{E}(\|\bX_{L,t}-\bX_t\|^2) &=\mathbb{E}\left\{\left\|\sum_{l=L+1}^\infty\int \bA_l(u,v)\bvarepsilon_{t-l}(v)dv\right\|^2\right\} 
        \\ &\leq \mathbb{E}
        \left(\sum_{l_1=L+1}^\infty\sum_{l_2 = L+1}^\infty\|\bA_{l_1}\|_\tF\|\bA_{l_2}\|_\tF\|\bvarepsilon_{t-l_1}\|\|\bvarepsilon_{t-l_2}\|\right) 
        \\ &\leq p\omega_0^{\bvarepsilon} \left(\sum_{l = L+1}^\infty\|\bA_{l}\|_\tF\right)^2.
    \end{split}
\end{equation}
By Lemma~\ref{lm_sigma_sub}, we have $\sum_{l = 0}^\infty\|\bA_{l}\|_\tF < \infty.$ This together with the above yields 
\begin{equation} \label{lm_x_expectation}
    \lim_{L \to \infty}\mathbb{E}\{\|\bX_{L,t}-\bX_t\|^2\} = 0,
\end{equation}
which completes the proof of (\ref{lm_Y_l2}).

\noindent
\textbf{Proof of (\ref{lm_YY_expectation})}. Next we show that $ \lim_{L \to \infty} \mathbb{E}\left[\bY_L^{\T}\bY_L\right] - \mathbb{E}\left[\bY^{\T}\bY\right] = 0.$ 
Write
\begin{equation}\nonumber
\begin{split}
& \quad \left|\mathbb{E}\left[\bY_L^{\T}\bY_L\right] - \mathbb{E}\left[\bY^{\T}\bY\right]\right| \\&= n\left| \langle\bPhi_1,(\bSigma_{L,0}^X - \bSigma_0^X)(\bPhi_1) \rangle\right| 
\\& =n\left| \int\bPhi_1^{\T}(u)\mathbb{E}\left(\bX_{L,t}(u)\bX_{L,t}^{\T}(v)- \bX_t(u)\bX_t^{\T}(v)\right )\bPhi_1(v)dudv \right|
\\& \leq n\left|\int\bPhi_1^{\T}\mathbb{E}\left(\bX_{L,t}(\bX_{L,t}-\bX_t)^{\T} \right)\bPhi_1dudv \right| 
\\ & \quad+ n\left| \int\bPhi_1^{\T}\mathbb{E}\left((\bX_{L,t}-\bX_t)\bX_t^{\T} \right)\bPhi_1dudv \right|.
\end{split}
\end{equation}
By Jensen's inequality and Lemma~\ref{lm_curveinequality}, we have
\begin{equation}\nonumber
    \begin{split}
        \left|\int\bPhi_1^{\T}\mathbb{E}\left(\bX_{L,t}(\bX_{L,t}-\bX_t)^{\T} \right)\bPhi_1dudv \right| ^2 &\leq \|\bPhi_1\|^4\mathbb{E}\{\|\bX_{L,t}\|^2\}\mathbb{E}\{\|\bX_{L,t}-\bX_t\|^2\},
       \\ \left| \int\bPhi_1^{\T}\mathbb{E}\left((\bX_{L,t}-\bX_t)\bX_t^{\T} \right)\bPhi_1dudv \right|^2 &\leq \|\bPhi_1\|^4\mathbb{E}\{\|\bX_{t}\|^2\}\mathbb{E}\{\|\bX_{L,t}-\bX_t\|^2\}.
    \end{split}
\end{equation}
Combining the above results with (\ref{lm_x_expectation}), we complete the proof of (\ref{lm_YY_expectation}).
$\square$
\\

\begin{lemma} \label{lm_M_L}
Suppose that conditions in Lemma~\ref{lm_middle matrix_ML} hold.  Let $\bX_{L,t}(u) = \sum\limits_{l=0}^{L}\bA_l(\bvarepsilon_{t-l}).$  
Then, for any  $\bPhi_1\in \mathbb{H}_0^p$ with $\|\bPhi_1\|_0 \leq k$ and $k=1,\dots,p,$
\begin{equation} \nonumber
    \lim_{L \to \infty} \cM(\bbf_L^X, \bPhi_1) =\cM(\bbf^X, \bPhi_1).
\end{equation}
\end{lemma}
\textbf{Proof}. By definitions of $\cM(\bbf^X, \bPhi)$ and $\bbf^X_\theta(\bPhi),$ we have
\begin{equation} \nonumber
    \begin{split}
        & \lim_{L \to \infty} |\cM(\bbf_L^X, \bPhi_1) - \cM(\bbf^X, \bPhi_1)| 
        \\&=  2\pi  \lim_{L \to \infty}  \left|\mathop{\text{ess sup}}\limits_{\theta \in [-\pi,\pi]} |\langle\bPhi_1,\bbf^X_{L,\theta}(\bPhi_1)\rangle| -   \mathop{\text{ess sup}}\limits_{\theta \in [-\pi,\pi]} |\langle\bPhi_1,\bbf^X_\theta(\bPhi_1)\rangle| \right|
        \\&  \leq 2\pi  \lim_{L \to \infty}  \mathop{\text{ess sup}}\limits_{\theta \in [-\pi,\pi]}\left| |\langle\bPhi_1,\bbf^X_{L,\theta}(\bPhi_1)\rangle| -  |\langle\bPhi_1,\bbf^X_\theta(\bPhi_1)\rangle| \right|
        \\& \leq \|\bPhi_1\|^2 \lim_{L \to \infty} \left\|\sum_{h \in \mathbb{Z}}(\bSigma_{L,h}^X - \bSigma_{h}^X  )\right\|_\tF  \quad (\text{by Lemma~\ref{lm_curveinequality} and $|\exp(-ih\theta)|=1$})
        \\& \leq \|\bPhi_1\|^2 \lim_{L \to \infty} \sum_{h \in \mathbb{Z}}\left\|\bSigma_{L,h}^X - \bSigma_{h}^X\right\|_\tF.
    \end{split}
\end{equation}
With $\|\bPhi_1\|^2 <\infty,$ it suffices to prove $\sum_{h= -\infty}^{\infty}\left\|\bSigma_{L,h}^X - \bSigma_{h}^X\right\|_\tF < \infty$ and $\lim_{L \to \infty}\left\|\bSigma_{L,h}^X - \bSigma_{h}^X\right\|_\tF = 0.$

By triangle inequality and Lemma~\ref{lm_sigma_sub}, we obtain that
\begin{equation} \nonumber
    \begin{split}
        \sum_{h= -\infty}^{\infty}\left\|\bSigma_{L,h}^X - \bSigma_{h}^X\right\|_\tF \leq 
        \sum_{h= -\infty}^{\infty} \|\bSigma_{L,h}^X\|_\tF + \sum_{h= -\infty}^{\infty}\|\bSigma_{h}^X\|_\tF < \infty.
    \end{split}
\end{equation}
We next prove $\lim_{L \to \infty}\left\|\bSigma_{L,h}^X - \bSigma_{h}^X\right\|_\tF = 0.$ Write
\begin{equation}\nonumber
    \begin{split}
       & \bSigma_h^{X}(u,v)  = \mathbb{E}\left(\bX_{t-h}(u)\bX_{t}^{\T}(v)\right) = \sum_{l=0}^\infty \int \bA_{l+h}(u,u')\bSigma_0^{\varepsilon}(u',v')\{\bA_{l}(v,v')\}^{\T} du'dv',
       \\& \bSigma_{L,h}^{X}(u,v)  = \mathbb{E}\left(\bX_{L,{t-h}}(u)\bX_{L,t}^{\T}(v)\right) = \sum_{l=0}^{L-h} \int \bA_{l+h}(u,u')\bSigma_0^{\varepsilon}(u',v')\{\bA_{l}(v,v')\}^{\T} du'dv'.
    \end{split}
\end{equation}
Then,
\begin{equation}\nonumber
    \begin{split}
        \lim_{L \to \infty}\left\|\bSigma_{L,h}^X - \bSigma_{h}^X\right\|_\tF &= \lim_{L \to \infty}\left\|\sum_{l=L-h+1}^{\infty} \int \bA_{l+h}(u,u')\bSigma_0^{\varepsilon}(u',v')\{\bA_{l}(v,v')\}^{\T} du'dv'\right\|_\tF
        \\& \leq p\omega_0^{\varepsilon}\lim_{L \to \infty}\sum_{l=L-h+1}^{\infty}\| \bA_l\|_\tF\|\bA_{l+h}\|_\tF \quad (\text{by Lemmas~\ref{lm_curveinequality} and \ref{lm_sigma_e}})
        \\& \leq p\omega_0^{\varepsilon}\lim_{L \to \infty}\sum_{l=L-h+1}^{\infty}\| \bA_l\|_\tF\sum_{l=L-h+1}^{\infty}\|\bA_{l+h}\|_\tF 
        \\& = 0 \quad (\text{by Lemma~\ref{lm_sigma_sub}}),
    \end{split}
\end{equation}
which completes the proof.
$\square$
\\
\begin{lemma} \label{lm_curveinequality}
(i)Let $\bA=(A_{ij})_{p \times q}$ with each $A_{ij}\in \mathbb{S}$ and $\bB=(B_1, \dots, B_q)^{\T} \in \mathbb{H}^q.$
\begin{equation}  \label{B6}
    \left\|\int \int\bA(u,v)\bB(v)dudv\right\| \leq \|\bA\|_\tF\|\bB\|.
\end{equation}
Similarly, we have 
\begin{equation} \label{B7}
\begin{split}
    \|\bA(u)\bB(u)\| &\leq \|\bA\|_\tF\|\bB\|,
    \\ \|\bA(u)\bB(v)\| &\leq \|\bA\|_\tF\|\bB\|,
    \\ \left\|\int\bA(u,v)\bB(v)dv\right\| &\leq \|\bA\|_\tF\|\bB\|,
\end{split}
\end{equation}
(ii)Let $\bA=(A_{ij})_{p \times q}$ with each $A_{ij}\in \mathbb{S}$ and $\bB=(B_{jk})_{q \times r}$ with each $B_{jk} \in \mathbb{S}.$ Then we have
\begin{equation} \label{B8}
\begin{split}
    \\ \left\|\int\bA(u,z)\bB(z,v)dz\right\|_\tF &\leq \|\bA\|_\tF\|\bB\|_\tF.
\end{split}
\end{equation}
\end{lemma}
\noindent \textbf{Proof of (\ref{B6})}. Let $C = \int \int\bA(u,v)\bB(v)dudv,$ we have $|C_{i}| = |\sum_k\int \int A_{ik}(u,v)B_{k}(v)dudv| \leq \sum_k\|A_{ik}\|_\cS \|B_{k}\|.$
\begin{equation} \nonumber
    \begin{split}
        \|C\|^2 = & \sum_{i} |C_{i}|^2 \leq \sum_{i}(\sum_k\|A_{ik}\|_\cS \|B_{k}\|)^2
        \\& \leq \sum_{i}(\sum_k\|A_{ik}\|_\cS^2)(\sum_k\|B_{k}\|^2) \quad (\text{by Cauchy-Schwarz inequality})
        \\ &\leq \sum_{i,k}\|A_{ik}\|_\cS^2\sum_{k}\|B_{k}\|^2 = \|\bA\|_\tF^2\|\bB\|^2.
    \end{split}
\end{equation} 

\noindent
\textbf{Proof of (\ref{B7})}. Let $C(u) = \bA(u)\bB(u),$ then $C_{i}(u) = \sum_k A_{ik}(u)B_{k}(u).$
\begin{equation} \nonumber
    \begin{split}
        \|C\|^2 = & \sum_{i} \int C_{i}(u)^2du = \sum_{i} \int \left\{\sum_k A_{ik}(u)B_{k}(u)\right\}^2du
        \\& \leq  \sum_{i} \left\{\sum_k  \int  A_{ik}(u)B_{k}(u) du\right\}^2
        \\& \leq \sum_{i}(\sum_k\|A_{ik}\| \|B_{k}\|)^2 \leq \|\bA\|_\tF^2\|\bB\|^2.
    \end{split}
\end{equation}
By similar arguments, we can prove the other two inequalities in (\ref{B7}).

\noindent
\textbf{Proof of (\ref{B8})}. Let $\bC(u,v) = \int\bA(u,z)\bB(z,v)dz,$ then  $C_{ij}(u,v) = \sum_k \int A_{ik}(u,z)B_{kj}(z,v)dz.$
\begin{equation} \nonumber
    \begin{split}
        \|\bC\|_\tF^2 = & \sum_{i,j} \int \int C_{ij}(u,v)^2dudv = \sum_{i,j} \int\int \left\{\sum_k \int A_{ik}(u,z)B_{kj}(z,v)dz\right\}^2dudv
        \\& \leq  \sum_{i,j} \left\{\sum_k  \int \int \int  A_{ik}(u,z)B_{kj}(z,v) dzdudv\right\}^2
        \\& \leq \sum_{i,j}(\sum_k\|A_{ik}\|_\cS \|B_{kj}\|_\cS)^2 \leq \|\bA\|_\tF^2\|\bB\|_\tF^2.
    \end{split}
\end{equation}
$\square$
\\

\begin{lemma} \label{lm_sigma_sub}
Suppose that conditions in Lemma~\ref{lm_middle matrix_ML} hold. Then we have
\begin{equation} \nonumber
    \sum_{l = 0}^\infty\|\bA_{l}\|_\tF < \infty
\end{equation}
and
\begin{equation} \nonumber
    \sum_{h \in \mathbb{Z}}\|\bSigma_h^{X}\|_\tF \leq 2p\omega_0^{\varepsilon}\left\{\sum_{l = 0}^\infty\|\bA_{l}\|_\tF\right\}^2 < \infty.
\end{equation}
\end{lemma}
\textbf{Proof}. It follows from Condition~\ref{con_sub_coefficient} that
\begin{equation} \nonumber
    \begin{split}
        \sum_{l = 0}^\infty\|\bA_{l}\|_\tF  &= \sum_{l = 0}^\infty\left\{\sum_{j,k}\|A_{l,jk}\|^2_\cS\right\}^{1/2}
        \\& \leq \sum_{l=0}^\infty \sum_{j}\|\bA_{l}\|_{\infty} < \infty.
    \end{split}
\end{equation}

Provided that $\bX_t(u) = \sum\limits_{l=0}^{\infty}\int \bA_l(u,v)\bvarepsilon_{t-l}(v)dv$ and $\bvarepsilon_{t}(\cdot)$'s are i.i.d. mean zero sub-Gaussian processes, we have
\begin{equation} \nonumber
\begin{split}
        \bSigma_h^{X}(u,v) & = \mathbb{E}\left\{\bX_{t-h}(u)\bX_{t}^{\T}(v)\right\}
        \\& = \sum_{l=0}^\infty \int \bA_{l+h}(u,u')\mathbb{E}\left\{\bvarepsilon_{t-l}(u')\bvarepsilon_{t-l}^{\T}(v')\right\}\{\bA_{l}(v,v')\}^{\T} du'dv'
    \\& = \sum_{l=0}^\infty \int \bA_{l+h}(u,u')\bSigma_0^{\varepsilon}(u',v')\{\bA_{l}(v,v')\}^{\T} du'dv'.
\end{split}
\end{equation}
This together with the fact that $\bSigma_{-h}^{X}(u,v) = \left\{\bSigma_h^{X}(v,u) \right\}^{\T} $ implies that
\begin{equation} \nonumber
    \begin{split}
        \sum_{h \in \mathbb{Z}}\|\bSigma_h^{X}(u,v)\|_\tF &\leq 2\sum_{h =0}^\infty\|\bSigma_h^{X}(u,v)\|_\tF 
        \\& =  2\sum_{h =0}^\infty\|\sum_{l=0}^\infty \int \bA_{l+h}(u,u')\bSigma_0^{\varepsilon}(u',v')\{\bA_{l}(v,v')\}^{\T} du'dv'\|_\tF
        \\& \leq 2\sum_{h =0}^\infty\sum_{l=0}^\infty \|\int \bA_{l+h}(u,u')\bSigma_0^{\varepsilon}(u',v')\{\bA_{l}(v,v')\}^{\T} du'dv'\|_\tF
        \\& \leq 2\sum_{h =0}^\infty\sum_{l=0}^\infty \|\bA_l\|_\tF\|\bA_{l+h}\|_\tF\|\bSigma_0^{\varepsilon} \|_\tF \quad(\text{by Lemma~\ref{lm_curveinequality}})
        \\&  \leq  2p\omega_0^{\varepsilon}\left\{\sum_{l = 0}^\infty\|\bA_{l}\|_\tF\right\}^2 < \infty \quad(\text{by Lemme \ref{lm_sigma_e}}),
    \end{split}
\end{equation}
which completes the proof.
$\square$
\\
\begin{lemma} \label{lm_sigma_e}
For a $p$-dimensional vector process $\{\bX_t(\cdot)\}_{t \in \mathbb{Z}},$ whose lag-$h$ autocovariance matrix function is $\bSigma_h =\left(\Sigma_{h,jk}\right)_{1 \leq j,k \leq p}$ with each $\Sigma_{h,jk}\in \mathbb{S}$ and $\omega_0 = \max_{1 \leq j \leq p}\int \Sigma_{0,jj}(u,u)du < \infty ,$ we have 
\begin{equation} \nonumber
    \|\Sigma_{h,jk}\|_\cS \leq \omega_0, \quad  \|\bSigma_h\|_\tF\leq p\omega_0, \quad \mathbb{E}(\|\bX_{t}\|) \leq \sqrt{p\omega_0} \quad\text{and}\quad \mathbb{E}(\|\bX_{t}\|^2) \leq {p\omega_0}.
\end{equation}
Let $X_{M,tj}(\cdot)= \sum_{l=1}^M \xi_{tjl} \phi_{jl}(\cdot) $ be the $M$-truncated process, we have
\begin{equation} \label{lm_sigma_e_truncation}
    \lim_{M \to \infty} \|\Sigma_{h,jk}^{X_{M}} - \Sigma_{h,jk}^X\|_\cS = 0.
\end{equation}
\end{lemma}
\textbf{Proof}. By $\Sigma_{h,jk} = \sum_{l,m=1}^\infty \mathbb{E}(\xi_{tjl}\xi_{(t+h)km})\phi_{jl}(u)\phi_{km}(v),$ orthonormality of $\{\phi_{jl}\}$ and Cauchy--Schwarz inequality, we obtain 
\begin{equation} \nonumber
    \begin{split}
        \|\Sigma_{h,jk}\|_\cS^2 &= \int \left\{\sum_{l,m=1}^\infty \mathbb{E}(\xi_{tjl}\xi_{(t+h)km})\phi_{jl}(u)\phi_{km}(v)\right\}^2dudv 
        \\&= \sum_{l,m=1}^\infty\mathbb{E}(\xi_{tjl}\xi_{(t+h)km})^2 \leq \sum_{l,m=1}^\infty\mathbb{E}(\xi_{tjl}^2)\mathbb{E}(\xi_{(t+h)km}^2) \leq \omega_0^2.
    \end{split}
\end{equation}
This implies that $\|\bSigma_h\|_\tF^2 = \sum_{j,k}\|\Sigma_{h,jk}\|_\cS^2 \leq p^2\omega_0^2.$ By the similar arguments, we have
\begin{equation} \nonumber
\begin{split}
    \|\Sigma_{h,jk}^{X_{M}} - \Sigma_{h,jk}^X\|_\cS^2 &= \int \left\{\sum_{l,m=M+1}^\infty \mathbb{E}(\xi_{tjl}\xi_{(t+h)km})\phi_{jl}(u)\phi_{km}(v)\right\}^2dudv 
        \\&= \sum_{l,m=M+1}^\infty\mathbb{E}(\xi_{tjl}\xi_{(t+h)km})^2 \leq \sum_{l,m=M+1}^\infty\mathbb{E}(\xi_{tjl}^2)\mathbb{E}(\xi_{(t+h)km}^2). 
    \end{split}
\end{equation}
Since $\sum_{l=0}^\infty\mathbb{E}(\xi_{tjl}^2) \leq \omega_0 < \infty,$ the above goes to zero when $M \to \infty,$ completing the proof of (\ref{lm_sigma_e_truncation}). 

Provided that $X_{tj}(\cdot)= \sum_{l=1}^\infty \xi_{tjl} \phi_{jl}(\cdot),$ orthonormality of $\{\phi_{jl}\}$ and Jensen's inequality, we have 

\begin{equation}\nonumber
    \begin{split}
        \mathbb{E}(\|\bX_{t}\|)& = \mathbb{E}\left\{\sqrt{\sum_{j=1}^p \int X_{tj}^2(u)du}\right\} \leq \sqrt{\sum_{j=1}^p \mathbb{E}\left\{\int X_{tj}^2(u)du\right\}}
        \\& \leq \sqrt{\sum_{j=1}^p \sum_{l=0}^{\infty}\mathbb{E} (\xi_{tjl}^2 )} \leq \sqrt{p\omega_0}.
    \end{split}
\end{equation}
Similarly, we obtain that $\mathbb{E}(\|\bX_{t}\|^2) = \mathbb{E}\left\{{\sum_{j=1}^p \int X_{tj}^2(u)du}\right\}
         = {\sum_{j} \sum_l \mathbb{E}(\xi_{tjl}^2) } \leq {p\omega_0}.$
$\square$
\\
\begin{lemma} \label{lm_sigma_h_xy}
For process $\{\bX_t(\cdot)\}_{t \in \mathbb{Z}}$ and $\{\bY_t(\cdot)\}_{t \in \mathbb{Z}},$ we have that
\begin{equation} \nonumber
    \|\Sigma_{h,jk}^{X,Y}\|_\cS \leq \sqrt{\omega_0^X \omega_0^Y},  
\end{equation}
and 
\begin{equation} \nonumber
    \|\langle\Sigma_{h,jk}^{X,Y}, \phi_{km} \rangle\| \leq \sqrt{\omega_0^X \omega_{km}^Y} \quad \text{and} \quad \|\langle\Sigma_{h,jk}^{X,Y}, \psi_{jl} \rangle\|\leq \sqrt{\omega_{jl}^X \omega_0^Y}.
\end{equation}
\end{lemma}
\textbf{Proof}. This lemma can be proved in similar way to Lemma~8 of \cite{guo2019} and hence the proof is omitted here.
$\square$

\section{Proofs of theoretical results in Section~\ref{sec_fully}}
\label{ap.thm.sec3}


We present the proof of Theorem~\ref{th_beta} in Appendix~\ref{ap.th.secB1} and proofs of Propositions~\ref{pr_fof_RE}--\ref{pr_fof_max.error} in Appendix~\ref{ap.pro.secB2}, followed by the supporting technical lemmas and their proofs in 
Appendix~\ref{ap.lemma.secB3}. 
 For a matrix $\bA \in {\eR}^{p \times q},$ we denote its  elementwise maximum norm by  $||\bA||_{\max}=\max_{i,j}|A_{ij}|$. 
To simplify our notation, for a square-block matrix $\bB = (\bB_{jk})_{1 \leq j \leq p_1, 1 \leq k \leq p_2} \in \eR^{p_1q \times p_2 q}$ with the $(j,k)$-th block $\bB_{jk} \in {\eR}^{q \times q},$ we use $\|\bB\|_{\max}^{(q)}$ and $\|\bB\|_{1}^{(q)}$ to denote its block versions of elementwise $\ell_{\infty}$ and matrix $\ell_1$ norms.

\subsection{Proof of Theorem~\ref{th_beta}}
\label{ap.th.secB1}
Denote the minimizer of (\ref{target}) by $\widehat{\bB}\in \mathbb{R}^{(L+1)pq_1\times q_2}.$ Then
\begin{equation}\nonumber
\frac{1}{2(n-L)}\|\widehat{\bU}-\widehat{\bZ}\widehat{\bD}^{-1}\widehat{\bB}\|_\tF^2 + \lambda_n\|\widehat{\bB}\|^{(q_1,q_2)}_1 \leq \frac{1}{2(n-L)}\|\widehat{\bU}-\widehat{\bZ}\widehat{\bD}^{-1}\bB\|_\tF^2 + \lambda_n\|\bB\|^{(q_1,q_2)}_1
\end{equation}
Letting $\bDelta = \widehat{\bB} - \bB$ and $S^c$ be the complement of $S$ in the set $\{0,\dots, L\} \times\{1,\dots,p\},$ we write
\begin{equation}\nonumber
\begin{split}
& \frac{1}{2}\llangle\bDelta,\widehat{\bGamma}\bDelta\rrangle\\
 \leq  &\frac{1}{n-L} \llangle\bDelta,\widehat{\bD}^{-1}\widehat{\bZ}^{\T}(\widehat{\bU}-\widehat{\bZ}\widehat{\bD}^{-1}\bB)\rrangle + \lambda_n(\|\bB\|_1^{(q_1,q_2)} - \|\bB+\bDelta\|_1^{(q_1,q_2)}) 
 \\
 \leq & \frac{1}{n-L} \llangle\bDelta,\widehat{\bD}^{-1}\widehat{\bZ}^{\T}(\widehat{\bU}-\widehat{\bZ}\widehat{\bD}^{-1}\bB)\rrangle + \lambda_n(\|\bDelta_S\|_1^{(q_1,q_2)} - \|\bDelta_{S^c}\|_1^{(q_1,q_2)}), 
\end{split}
\end{equation}
where $\widehat{\bGamma} = 
(n-L)^{-1}
\widehat{\bD}^{-1}\widehat{\bZ}^{\T}\widehat{\bZ}\widehat{\bD}^{-1}.$ By Proposition~\ref{pr_fof_max.error} and $\lambda_n \geq 2 C_0  s  q_1^{1/2}((\cM_1^X+ \cM^\epsilon) \vee \cM_1^Y)\{(  q_1^{\alpha_1+3/2} \vee   q_2^{\alpha_2+3/2})\sqrt{\frac{\log(pq_1q_2)}{n}}+ q_1^{-\kappa+1/2}\},$ we have
\begin{equation}\nonumber
\begin{split}
& \quad \frac{1}{n-L}\vert\llangle\bDelta,\widehat{\bD}^{-1}\widehat{\bZ}^{\T}(\widehat{\bU}-\widehat{\bZ}\widehat{\bD}^{-1}\bB)\rrangle| \\&\leq\frac{1}{n-L} \|\widehat{\bD}^{-1}\widehat{\bZ}^{\T}(\widehat{\bU}-\widehat{\bZ}\widehat{\bD}^{-1}\bB)\|_{\max}^{(q_1,q_2)}\|\bDelta\|_1^{(q_1,q_2)}
\\&\leq\frac{\lambda_n}{2}(\|\bDelta_S\|_1^{(q_1,q_2)} +\|\bDelta_{S^c}\|_1^{(q_1,q_2)}).
\end{split}
\end{equation}
This implies that
\begin{equation}\nonumber
\begin{split}
0\leq\frac{1}{2}\llangle\bDelta,\widehat{\bGamma}\bDelta\rrangle \leq  \frac{3\lambda_n}{2}\|\bDelta_S\|_1^{(q_1,q_2)} - \frac{\lambda_n}{2}\|\bDelta_{S^c}\|_1^{(q_1,q_2)}\leq \frac{3}{2}\lambda_n\|\bDelta\|_1^{(q_1,q_2)}.
\end{split}
\end{equation}
Therefore $\|\bDelta\|_1^{(q_1,q_2)}\leq 4\|\bDelta_S\|_1^{(q_1,q_2)} \leq 4\sqrt{s}\|\bDelta\|_\tF.$ By Proposition~$\ref{pr_fof_RE}$ and $\tau_2\geq 32\tau_1q_1q_2s,$ we obtain 
\begin{equation}\nonumber
\begin{split}
\llangle\bDelta,\widehat{\bGamma}\bDelta\rrangle\geq\tau_2\|\bDelta\|_\tF^2 - \tau_1q_1q_2\{\|\bDelta\|_1^{(q_1,q_2)}\}^2\geq(\tau_2 - 16\tau_1q_1q_2s)\|\bDelta\|_\tF^2\geq\frac{\tau_2}{2}\|\bDelta\|_\tF^2.
\end{split}
\end{equation}
Therefore,
\begin{equation}\nonumber
\begin{split}
\frac{\tau_2}{4}\|\bDelta\|_\tF^2 \leq \frac{3}{2} \lambda_n \|\bDelta\|_1^{(q_1,q_2)} \leq 6\lambda_n s^{1/2}\|\bDelta\|_\tF,
\end{split}
\end{equation}
which implies that 
\begin{equation}\label{eq_delta_1^q}
\begin{split}
\|\bDelta\|_\tF\leq \frac{24s^{1/2}\lambda_n}{\tau_2} \; \text{and} \; \|\bDelta\|_1^{(q_1,q_2)}\leq \frac{96s\lambda_n}{\tau_2}.
\end{split}
\end{equation}

Here, we aim to prove the upper bound of $\|\widehat{\bbeta}-\bbeta\|_{1}.$ For each $(h,j) \in S$ we have,
\begin{equation} \nonumber
\begin{split}
\widehat\beta_{hj}-\beta_{hj} & = \widehat{\bpsi}_j(u)^{\T}\widehat{\bPsi}_{hj}\widehat{\bphi}(v) - \bpsi_j(u)^{\T}\bPsi_{hj}\bphi(v) + R_{hj}(u,v)
\\& =   (\widehat{\bpsi}_j(u)-\bpsi_j(u))^{\T}\widehat{\bPsi}_{hj}\widehat{\bphi}(v) + \bpsi_j(u)^{\T}\widehat{\bPsi}_{hj}(\widehat{\bphi}(v) -\bphi(v))  
\\&  \quad + \bpsi_j(u)^{\T}(\widehat{\bPsi}_{hj}-\bPsi_{hj})\bphi(v) + R_{hj}(u,v),
\end{split}
\end{equation}
where $R_{hj}(u,v) = -\sum_{l=q_1+1}^\infty\sum_{m=q_2+1}^\infty a_{hjlm}\psi_{jl}(u)\phi_m(v).$ Therefore,
\begin{equation} \label{beta1}
\begin{split}
\|\widehat{\bbeta}-\bbeta\|_{1} \leq &\sum\limits_{h,j}\|(\widehat{\bpsi}_j(u)-\bpsi_j(u))^{\T}\widehat{\bPsi}_{hj}\widehat{\bphi}(v)\|_\cS + \sum\limits_{h,j}\|\bpsi_j(u)^{\T}\widehat{\bPsi}_{hj}(\widehat{\bphi}(v) -\bphi(v))\|_\cS 
\\& + \sum\limits_{h,j}\|\bpsi_j(u)^{\T}(\widehat{\bPsi}_{hj}-\bPsi_{hj})\bphi(v)\|_\cS +
\sum\limits_{h,j}\|R_{hj}(u,v)\|_\cS.
\end{split}
\end{equation}
Due to the orthonormality of $\{\psi_{jl}(\cdot)\}$ and $\{\phi_{m}(\cdot)\}$ and the estimated eigenfunctions $\{\widehat{\psi}_{jl}(\cdot)\}$ and $\{\widehat{\phi}_{m}(\cdot)\},$
\begin{equation} \nonumber
\begin{split}
&\|(\widehat{\bpsi}_j(u)-\bpsi_j(u))^{\T}\widehat{\bPsi}_{hj}\widehat{\bphi}(v)\|_\cS  \leq q_1^{1/2}\|\widehat{\bPsi}_{hj}\|_\tF\max\limits_l\|\widehat{\psi}_{jl}-\psi_{jl}\|,
\\& \|\bpsi_j(u)^{\T}\widehat{\bPsi}_{hj}(\widehat{\bphi}(v) -\bphi(v))\|_\cS  \leq q_2^{1/2}\|\widehat{\bPsi}_{hj}\|_\tF\max\limits_m\|\widehat{\phi}_{m}-\psi_{m}\|,
\\& \|\bpsi_j(u)^{\T}(\widehat{\bPsi}_{hj}-\bPsi_{hj})\bphi(v)\|_\cS  = \|\widehat{\bPsi}_{hj}-\bPsi_{hj}\|_\tF.
\end{split}
\end{equation} 
To bound the first three terms of (\ref{beta1}), we start with  the upper bound of $\sum\limits_{h,j}\|\widehat{\bPsi}_{hj}-\bPsi_{hj}\|_\tF= \|\widehat{\bPsi} - \bPsi\|_1^{(q_1,q_2)}$ and $\sum\limits_{h,j}\|\widehat{\bPsi}_{hj}\|_\tF= \|\widehat{\bPsi}\|_1^{(q_1,q_2)}.$ From Condition~\ref{con_fof_beta_a}, for $(h,j) \in S,$ $\|\bPsi_{hj}\|_\tF = \{\sum_{l=1}^{q_1}\sum_{m=1}^{q_2}\mu_{hj}^2(l+m)^{-2\kappa-1}\}^{1/2} \leq \{\mu_{hj}^2\int_1^{q_2}\int_1^{q_1}(x+y)^{-2\kappa-1}dxdy\}^{1/2} = O(\mu_{hj}).$ For $(h,j) \in S^c,$ $\bPsi_{hj} = 0.$ Hence,
\begin{equation}\label{eq_psi}
\begin{split}
\|\bPsi\|_1^{(q_1,q_2)} = \sum\limits_{h,j}\|\bPsi_{hj}\|_\tF =  O(s).
\end{split}
\end{equation} 
By the definition of $\omega_0^X,$ Condition~\ref{con_eigen} and Proposition \ref{pr_fof_eigen}, we have $\|\bD\|_{\max} \leq \sqrt{\omega_0^X},$ $\|\bD^{-1}\|_{\max} \leq \alpha_1^{1/2}c_0^{-1/2}q_1^{\alpha_1/2}$ and $$\|\widehat{\bD}^{-1}-\bD^{-1}\|_{\max} \leq \alpha_1^{1/2}c_0^{-1/2}q_1^{\alpha_1/2} C_\omega\cM_1^X\sqrt{\frac{\log(pq_1)}{n}}.$$
Recall that $\widehat{\bPsi} - \bPsi = \widehat{\bD}^{-1}\widehat{\bB} -\bD^{-1}\bB =  \bD^{-1}(\widehat{\bB} - \bB) + (\widehat{\bD}^{-1}-\bD^{-1})\widehat{\bB}.$ Then 
\begin{eqnarray*}	
\|\widehat{\bPsi} - \bPsi\|_1^{(q_1,q_2)} &\leq&  \|\bD^{-1}\|_{\max}\|\widehat{\bB} - \bB\|_1^{(q_1,q_2)} + \|\widehat{\bD}^{-1}-\bD^{-1}\|_{\max}\|\widehat{\bB}\|_1^{(q_1,q_2)}
\\& \leq &  \|\bD^{-1}\|_{\max}\|\widehat{\bB} - \bB\|_1^{(q_1,q_2)} + \|\widehat{\bD}^{-1}-\bD^{-1}\|_{\max}\|\widehat{\bB} - \bB\|_1^{(q_1,q_2)}
\\& &+ \|\widehat{\bD}^{-1}-\bD^{-1}\|_{\max}\|\bB\|_1^{(q_1,q_2)}
\\& \leq &   \|\bD^{-1}\|_{\max}\|\widehat{\bB} - \bB\|_1^{(q_1,q_2)} + \|\widehat{\bD}^{-1}-\bD^{-1}\|_{\max}\|\widehat{\bB} - \bB\|_1^{(q_1,q_2)} 
\\& & + \|\widehat{\bD}^{-1}-\bD^{-1}\|_{\max}\|\bD\|_{\max}\|\bPsi\|_1^{(q_1,q_2)}.
\end{eqnarray*}
This, together with (\ref{eq_delta_1^q}) implies that,
\begin{equation} \label{eq_B_1_block}
\begin{split}
\|\bB\|_{1}^{(q_1,q_2)} = O(\sqrt{\omega_0^X}s),
\end{split}
\end{equation}
and
\begin{equation}\label{psiwidehat-psi}
\begin{split}
\|\widehat{\bPsi} - \bPsi\|_1^{(q_1,q_2)} \leq  \frac{96\alpha_1^{1/2}q_1^{\alpha_1/2} s \lambda_n}{c_0^{1/2}\tau_2}\left\{1+o(1)\right\}.
\end{split}
\end{equation} 
Combining (\ref{eq_psi}) and (\ref{psiwidehat-psi}), we have
\begin{equation}\nonumber
\begin{split}
\|\widehat{\bPsi}\|_1^{(q_1,q_2)} =  O(s).
\end{split}
\end{equation} 
To bound the fourth term of (\ref{beta1}), we have that $\|R_{hj}\|_\cS = O(\|\sum_{l=1}^{q_1} \sum_{m=q_2+1}^\infty a_{hjlm}\psi_{jl}\phi_m \|_\cS \vee \| \sum_{l=q_1+1}^{\infty} \sum_{m=1}^{q_2} a_{hjlm}\psi_{jl}\phi_m\|_\cS )
= O(\mu_{hj}\min (q_1 , q_2)^{-\kappa+1/2}),$ for each $(h,j) \in S.$ For $(h,j) \in S^c,$ $\|R_{hj}\|_\cS = 0.$ Hence, $\sum\limits_{h,j}\|R_{hj}\|_\cS = O(s \min(q_1 , q_2)^{-\kappa+1/2}).$

Combining all the results with Proposition~\ref{pr_fof_eigen}, we obtain 
\begin{equation} \nonumber 
\begin{split}
\|\widehat{\bbeta}-\bbeta\|_{1} &\leq \|\widehat{\bPsi}\|_1^{(q_1,q_2)}\left\{q_1^{1/2}\max\limits_{j,l}\|\widehat{\psi}_{jl}-\psi_{jl}\|_\cS +q_2^{1/2}\max\limits_m\|\widehat{\phi}_{m}-\phi_{m}\|_\cS\right\} 
\\&\quad+\|\widehat{\bPsi} - \bPsi\|_1^{(q_1,q_2)} +
\sum\limits_{h,j}\|R_{hj}\|_\cS
\\&\leq  \frac{96\alpha_1^{1/2}q_1^{\alpha_1/2} s \lambda_n}{c_0^{1/2}\tau_2}\left\{1+o(1)\right\},
\end{split}
\end{equation}
which completes the proof.
$\square$

\subsection{Proofs of propositions}
\label{ap.pro.secB2}

\paragraph{Proof of Proposition~\ref{pr_fof_RE}}
Define $\bGamma = (n-L)^{-1}{\bD}^{-1}\mathbb{E}\{{\bZ}^{\T}\bZ\}\bD^{-1}.$ Note that $\btheta^{\T} \widehat \bGamma \btheta  = \btheta^{\T} \bGamma \btheta + \btheta^{\T} (\widehat \bGamma - \bGamma) \btheta.$ Hence we have
\begin{eqnarray*} \nonumber
	\btheta^{\T} \widehat \bGamma \btheta	 \ge \btheta^{\T} \bGamma \btheta -  \|\widehat \bGamma - \bGamma\|_{\max} \|\btheta\|_1^2.
\end{eqnarray*}
By Condition~\ref{con_fof_eigen_min}, $\omega_{\min}(\bGamma) \ge \underline{\mu},$ where $\omega_{\min}(\bGamma)$ denotes the minimum eigenvalue of $\bGamma.$ This, together with Lemma~\ref{lm_gamma}, completes our proof. 
$\square$
\\
\paragraph{Proof of Proposition~\ref{pr_fof_eigen}}
This proposition can be proved in similar way to Proposition~3 of \cite{guo2019} and hence the proof is omitted here. $\square$
\\

\paragraph{Proof of Proposition~\ref{pr_fof_max.error}}
Notice that $\widehat{\bU} =  \bZ  \bD^{-1}\tilde\bB + \widehat \bR + \widehat{\bE},$ where $ \widetilde\bB =  \bD\widetilde \bPsi$ and $\{(h+1)j\}$-th row block of $\widetilde \bPsi,$  $\widetilde\bPsi_{hj} = \int_\cV\int_\cU \bpsi_j(u) \beta_{hj}(u,v) \widehat \bphi(v)^{\T} du dv.$ The matrix $\widehat \bR$ and $\widehat \bE$ are both $(n-L)\times q_2$ matrices whose row vectors are formed by $\{\widehat \br_{t} = (\widehat r_{t1},\dots,\widehat r_{tq_2})^{\T}\}_{L+1}^n$ and $\{\widehat\beps_{t}= (\widehat \epsilon_{t1},\dots,\widehat \epsilon_{tq_2})^{\T}\}_{L+1}^n$ respectively, where $\widehat r_{tm} = \sum_{h=0}^{L} \sum_{j=1}^p\sum_{l=q_1+1}^\infty\langle\langle  \psi_{jl},\beta_{hj}\rangle,\widehat \phi_m \rangle  \zeta_{tjl}$ and $\widehat \epsilon_{tm} = \langle\epsilon_t,\widehat \phi_m\rangle.$ Then we rewrite
\begin{equation}\nonumber
\begin{split}
    & \quad \frac{1}{n-L} \widehat{\bD}^{-1}\widehat{\bZ}^{\T}(\widehat{\bU}-\widehat{\bZ}\widehat{\bD}^{-1}\bB)
     \\&=  \frac{1}{n-L} \widehat{\bD}^{-1}\widehat{\bZ}^{\T} ( \bZ  \bD^{-1}\tilde\bB - \widehat \bZ \widehat \bD^{-1}\bB) +  \frac{1}{n-L} \widehat{\bD}^{-1}\widehat{\bZ}^{\T} \widehat \bR  +\frac{1}{n-L} \widehat{\bD}^{-1}\widehat{\bZ}^{\T} \widehat{\bE}
     \\&= I_1 + I_2 +I_3.
     \end{split}
\end{equation}
Next, we show the deviation bounds of the above three parts.
\begin{equation} \nonumber
\begin{split}
   & \|I_1\|_{\max} ^{(q_1,q_2)} 
   \\=& \|\frac{1}{n-L}
   \widehat{\bD}^{-1}\widehat{\bZ}^{\T}(\bZ  \bD^{-1} - \widehat \bZ \widehat \bD^{-1})\bB\|_{\max}^{(q_1,q_2)} + \|\frac{1}{n-L} \widehat{\bD}^{-1}\widehat{\bZ}^{\T} \bZ  \bD^{-1}(\tilde\bB -\bB)\|_{\max} ^{(q_1,q_2)} 
    \\\leq &\|\frac{1}{n-L} \widehat{\bD}^{-1}\widehat{\bZ}^{\T}(\bZ  \bD^{-1} - \widehat \bZ \widehat \bD^{-1})\|_{\max}^{(q_1)}\|\bB\|_{1}^{(q_1,q_2)} \\ &+ \|\frac{1}{n-L} \widehat{\bD}^{-1}\widehat{\bZ}^{\T} \bZ  \bD^{-1}\|_{\max}^{(q_1)}\|\tilde\bB -\bB\|_{1} ^{(q_1,q_2)}
    \\ \leq &\|\frac{1}{n-L} \widehat{\bD}^{-1}\widehat{\bZ}^{\T}(\bZ  \bD^{-1} - \widehat \bZ \widehat \bD^{-1})\|_{\max}^{(q_1)}\|\bB\|_{1}^{(q_1,q_2)} + \|\widehat \bGamma\|_{\max}^{(q_1)}\|\tilde\bB -\bB\|_{1} ^{(q_1,q_2)} 
    \\ \quad &+  \|\frac{1}{n-L} \widehat{\bD}^{-1}\widehat{\bZ}^{\T}(\bZ  \bD^{-1} - \widehat \bZ \widehat \bD^{-1})\|_{\max}^{(q_1)}\|\tilde\bB -\bB\|_{1}^{(q_1,q_2)},
\end{split}
\end{equation}
where $\widehat{\bGamma} = (n-L)^{-1}\widehat{\bD}^{-1}\widehat{\bZ}^{\T}\widehat{\bZ}\widehat{\bD}^{-1}.$ By Lemmas~\ref{lm_blockinequality}, \ref{lm_gamma_2}, \ref{lm_B tilde} and (\ref{eq_B_1_block}) in Appendix~\ref{ap.th.secB1}, there exist some positive constants $C_1^*, c_1^*$ and $c_2^*$ such that 
\begin{equation}  \label{eq_pr_fof.max.error_1}
    \|I_1\|_{\max} ^{(q_1,q_2)} \leq C_1^* s  q_1^{1/2}(\cM_1^X  q_1^{\alpha_1+3/2} \vee \cM_1^Y  q_2^{\alpha_2+3/2})\sqrt{\frac{\log(pq_1\vee q_2)}{n}}
\end{equation}
with probability greater than $1 - c_1^*(pq_1\vee q_2)^{-c_2^*}.$

By Lemma~\ref{lm_Residual}, we obtain that there exist some positive constants $C_2^*, c_1^*$ and $c_2^*$ such that
\begin{equation} 
\begin{split}
     \|I_2\|_{\max}^{(q_1,q_2)}  \leq C_2^*  s q_1^{-\kappa+1} 
     \end{split}
\end{equation}
with probability greater than $1 - c_1^* (pq_1q_2)^{-c_2^*}.$

Let $\bQ = ((n-L)^{-1}\bU^{\T}\bU)^{1/2} = \text{diag}(\{\omega_{1}^Y\}^{1/2}, \dots, \{\omega_{q}^Y\}^{1/2}).$ It follows from Proposition~\ref{pr_FPCscores_XE}  and $\|\bQ\|_\tF \leq \sqrt{\omega_0^Y}$ that there exist some positive constants $C_3^*, c_1^*$ and $c_2^*$ such that
\begin{equation}\label{eq_pr_fof.max.error_2}
\begin{split}
     \|I_3\|_{\max}^{(q_1,q_2)}  &\leq q_1^{1/2}\|\widehat \bD^{-1} \bD\|_{\max}\|(n-L)^{-1} \bD^{-1}\widehat{\bZ}^{\T} \widehat \bE \bQ^{-1}\|_{\max} \|\bQ\|_\tF
     \\& \leq C_3^* q_1^{1/2}(\cM_1^X + \cM^\epsilon)(q_1^{\alpha_1}\vee q_2^{\alpha_2})\sqrt{\frac{\log(pq_1q_2)}{n}}
     \end{split}
\end{equation}
with probability greater than $1 - c_1^*(pq_1q_2)^{-c_2^*}.$

It follows from (\ref{eq_pr_fof.max.error_1})--(\ref{eq_pr_fof.max.error_2}) that there exist some positive constants $C_0, c_1^*$ and $c_2^*$ such that
\begin{eqnarray*}
    & &\frac{1}{n-L} \|\widehat{\bD}^{-1}\widehat{\bZ}^{\T}(\widehat{\bU}-\widehat{\bZ}\widehat{\bD}^{-1}\bB)\|_{\max}^{(q_1,q_2)} 
    \\&\leq& C_0  s  q_1^{1/2}((\cM_1^X+ \cM^\epsilon) \vee \cM_1^Y)\{(  q_1^{\alpha_1+3/2} \vee   q_2^{\alpha_2+3/2})\sqrt{\frac{\log(pq_1q_2)}{n}}+ q_1^{-\kappa+1/2}\}
\end{eqnarray*}
with probability greater than $1 - c_1^*(pq_1q_2)^{-c_2^*},$ which completes the proof.
$\square$

\subsection{Technical lemmas and their proofs}
\label{ap.lemma.secB3}
\begin{lemma} \label{lm_blockinequality}
$\| \widehat \bGamma\|_{\max}^{(q_1)} = O(q_1^{1/2}).$
\end{lemma}
{\bf Proof.} For a semi-positive definite block matrix
\begin{equation} \nonumber
    \bA = \left( 
    \begin{matrix}
    \boldsymbol{L} & \boldsymbol{X}  \\
    \boldsymbol{X^{\T}} & \boldsymbol{M}
    \end{matrix}
    \right),
\end{equation}
we have that $\|\bX\|_\tF^2 \leq \|\boldsymbol{L}\|_\tF\|\bM\|_\tF.$ This can be seen as a special case of $p = 1$ in Theorem~4.2 of \cite{horn1990}.
Without loss of generality, we take $L = 0$ as an example. Let $\widehat \bGamma_{jk} = (\widehat \Gamma_{jl,km})_{1\leq l,m \leq q_1}.$ Then for $j = k,$ by the diagonal structure of $\widehat \bGamma_{jj},$ we have $\|\widehat \bGamma_{jj}\|_\tF = O(q_1^{1/2}).$ Applying the above inequality, we obtain  $\|\widehat \bGamma_{jk}\|_\tF \leq \sqrt{\|\widehat\bGamma_{jj}\|_\tF\|\widehat\bGamma_{kk}\|_\tF} = O(q_1^{1/2}).$
$\square$
\\

\begin{lemma} \label{lm_gamma}
Suppose that Conditions~\ref{con_stability}--\ref{con_eigen} hold. Then there exist some positive constants $C_{\Gamma},$ $c_1^*$ and $c_2^*$ such that
\begin{equation} \nonumber
    	\big\|\widehat \bGamma - \bGamma\big\|_{\max} \leq C_\Gamma \cM_1^X q_1^{\alpha_1 +1}\sqrt{\frac{\log (pq_1)}{n}}
\end{equation}
with probability greater than $1 - c_1^*(pq_1)^{-c_2^*}.$
\end{lemma}
{\bf Proof.} The proof follows from Lemma~5 in \cite{guo2019}. $\square$
\\

\begin{lemma} \label{lm_gamma_2}
Suppose that Conditions~\ref{con_stability}--\ref{con_eigen} hold.  Then there exist some positive constants $\tilde C_{\Gamma},$ $c_1^*$ and $c_2^*$ such that
\begin{equation} \nonumber
    	\|\frac{1}{n-L} \widehat{\bD}^{-1}\widehat{\bZ}^{\T}(\bZ  \bD^{-1} - \widehat \bZ \widehat \bD^{-1})\|_{\max} \leq \tilde C_\Gamma \cM_1^X q_1^{\alpha_1 +1}\sqrt{\frac{\log (pq_1)}{n}}
\end{equation}
with probability greater than $1 - c_1^*(pq_1)^{-c_2^*}.$
\end{lemma}
{\bf Proof.} We first consider $\|\frac{1}{n-L} \widehat{\bD}^{-1}\widehat{\bZ}^{\T}\bZ  \bD^{-1} - \bGamma\|_{\max}$. By Lemma~\ref{lm_guo2019_phi_g}, Proposition~\ref{pr_fof_eigen} and following the similar argument in the proof of Lemma~\ref{lm_guo2019_score}, we obtain that
\begin{equation} \nonumber
\begin{split}
    &\max_{j,k,l,m}\frac{(n-L)^{-1} \sum_{t=L+1}^{n} \widehat\zeta_{(t-h)jl} \zeta_{tkm}}{\sqrt{\widehat \omega_{jl}^X\omega_{km}^X}} -\frac{\mathbb{E}(\zeta_{(t-h)jl} \zeta_{tkm})}{\sqrt{\omega_{jl}^X \omega_{km}^X}}
    \\&\lesssim  \max_{j,k,l,m}\frac{\langle\widehat{\psi}_{jl}, \langle \widehat \Sigma_{h,jk}^{X}, \phi_{km} \rrangle - \langle{\psi}_{jl}, \langle \Sigma_{h,jk}^{X}, \phi_{km} \rrangle  }{\sqrt{\omega_{jl}^X \omega_{km}^X}}
    \\& \lesssim  \max_{j,k,l,m}\frac{\langle\widehat{\psi}_{jl} - {\psi}_{jl} , \langle  \Sigma_{h,jk}^{X}, \phi_{km} \rrangle + \langle\widehat{\psi}_{jl}, \langle \widehat \Sigma_{h,jk}^{X}- \Sigma_{h,jk}^{X}, \phi_{km} \rrangle  }{\sqrt{\omega_{jl}^X \omega_{km}^X}} 
    \\&
    \lesssim \cM_1^X q_1^{\alpha_1 +1}\sqrt{\frac{\log (pq_1)}{n}}
\end{split}
\end{equation}
holds with probability greater than $1 - c_1^*(pq_1)^{-c_2^*}.$ This, together with Lemma~\ref{lm_gamma}, shows that
\begin{equation} \nonumber
\begin{split}
    &\quad \|\frac{1}{n-L} \widehat{\bD}^{-1}\widehat{\bZ}^{\T}(\bZ  \bD^{-1} - \widehat \bZ \widehat \bD^{-1})\|_{\max}
    \\& \leq \|\frac{1}{n-L} \widehat{\bD}^{-1}\widehat{\bZ}^{\T}\bZ  \bD^{-1} - \bGamma\|_{\max} +\|\widehat \bGamma - \bGamma\|_{\max}
    \\&=O_P\{ \cM_1^X q_1^{\alpha_1 +1}\sqrt{\frac{\log (pq_1)}{n}}\}.
\end{split}
\end{equation}
 $\square$
\\

\begin{lemma} \label{lm_B tilde}
Suppose that Conditions~\ref{con_stability}--\ref{con_fof_beta_a} hold. Then there exist some positive constants $C_{B},$ $c_1^*$ and $c_2^*$ such that
\begin{equation} \nonumber
    \|\widetilde\bB - \bB\|_1^{(q_1,q_2)} \leq  C_B s \cM_1^Y  q_2^{\alpha_2+3/2}\sqrt{\frac{\log(q_2)}{n}}
\end{equation}
with probability greater than $1 - c_1^*( q_2)^{-c_2^*}.$
\end{lemma}
{\bf Proof.} We start with the convergence rate of $\| \widetilde \bPsi - \bPsi\|_{1}^{(q_1,q_2)}.$ Elementwisely, for fixed $h$, $j$ and $l = 1, \dots, q_1, m = 1, \dots, q_2,$ we have that
\begin{equation}\nonumber
\begin{split}
    &\langle\langle  \psi_{jl},\beta_{hj}\rangle,\widehat \phi_m \rangle - \langle\langle  \psi_{jl},\beta_{hj}\rangle, \phi_m\rangle =   \langle\langle  \psi_{jl},\beta_{hj}\rangle, \widehat\phi_m  - \phi_m \rangle =I_1
\end{split}
\end{equation}
Recall that $\beta_{hj} = \sum_{l,m=1}^\infty a_{hjlm} \psi_{jl}(u)\phi_m(v)$ and $|a_{hjlm}| \leq u_{hj}(l+m)^{-\kappa-1/2}.$ 
\begin{equation}\nonumber
    \begin{split}
        I_1 &= \langle\langle \psi_{jl},\sum_{l',m'=1}^\infty a_{hjl'm'} \psi_{jl'}\phi_{m'}\rangle, \widehat \phi_m  - \phi_m \rangle = \sum_{m'=1}^\infty a_{hjlm'} \langle\phi_{m'},\widehat \phi_{m} -\phi_{m}\rangle
        \\& \lesssim \|\widehat \phi_{m} -\phi_{m}\| u_{hj} l^{- \kappa + 1/2}.
    \end{split}
\end{equation}
It follows from Lemma~\ref{lm_guo2019_eigen}, for $(h,j) \in S,$
\begin{equation} \nonumber
\begin{split}
    \| \widetilde \bPsi_{hj} - \bPsi_{hj}\|_\tF & =   \sqrt{\sum_{l=1}^{q_1} \sum_{m=1}^{q_2}I_1^2} \lesssim u_{hj} q_2^{1/2}\max_{1 \leq m \leq q_2} \|\widehat \phi_m - \phi_m\|
    \\& =O_P\{ u_{hj}  \cM_1^Y  q_2^{\alpha_2+3/2}\sqrt{\frac{\log(q_2)}{n}}\}.
    \end{split}
\end{equation}
Then $\| \widetilde \bPsi - \bPsi\|_{1}^{(q_1,q_2)} = \sum_{h = 0}^ L\sum_{j = 1}^ p\| \widetilde \bPsi_{hj} - \bPsi_{hj}\|_\tF  =O_P\{ s  \cM_1^Y  q_2^{\alpha_2+3/2}\sqrt{\frac{\log(q_2)}{n}}\}.$ This result, together with  $\|\bD\|_{\max} \leq \{\omega_0^X\}^{1/2}$, implies that there exists $C_B$ such that
\begin{eqnarray*}
    \|\tilde\bB - \bB\|_1^{(q_1,q_2)} &=&\|\bD (\widetilde \bPsi - \bPsi)\|_{1}^{(q_1,q_2)} \leq \|\bD\|_{\max}\| \widetilde \bPsi - \bPsi\|_{1}^{(q_1,q_2)}
    \\&\leq&  C_B s  \cM_1^Y  q_2^{\alpha_2+3/2}\sqrt{\frac{\log(q_2)}{n}},
\end{eqnarray*}
with probability greater than $1 - c_1^*(q_2)^{-c_2^*}.$  $\square$
\\

\begin{lemma} \label{lm_Residual}
Suppose that Conditions~\ref{con_stability}--\ref{con_fof_beta_a} hold. Then there exist some positive constants $C_{R},$ $c_1^*$ and $c_2^*$ such that
\begin{equation} \nonumber
\begin{split}
     \|(n-L)^{-1} \widehat{\bD}^{-1}\widehat{\bZ}^{\T} \widehat \bR\|_{\max}^{(q_1,q_2)}  \leq C_R s q_1^{-\kappa+1} 
     \end{split}
\end{equation}
with probability greater than $1 - c_1^* (pq_1q_2)^{-c_2^*}.$
\end{lemma}
{\bf Proof.} Recall $\widehat r_{tm} = \sum_{h = 0}^L\sum_{j=1}^p\sum_{l=q_1+1}^\infty\langle\langle  \psi_{jl},\beta_{hj}\rangle,\widehat \phi_m \rangle  \zeta_{tjl}= \sum_{h = 0}^L\sum_{j=1}^p\tilde r_{tmhj}.$ The matrix $\widehat \bR$ are $(n-L)\times q_2$ matrices whose row vectors are formed by $\{\widehat \br_{t} = (\widehat r_{t1},\dots,\widehat r_{tq_2})^{\T},t= L+1,\dots,n\}.$ By Cauchy-Schwarz inequality and the definition of $\widehat \omega_{jl}^X$, we obtain
\begin{equation} \nonumber
\begin{split}
    &\frac{(n-L)^{-1} \sum_{t=L+1}^{n} \widehat\zeta_{(t-h)jl} \sum_{h = 0}^L\sum_{j'=1}^p\tilde r_{tmhj'}}{\{\widehat \omega_{jl}^X\}^{1/2}}
    \leq  \sum_{h = 0}^L\sum_{j'=1}^p\sqrt{(n-L)^{-1}\sum_{t=L+1}^{n}  \tilde r_{tmhj'}^2}
    \\&=
    \sum_{h = 0}^L\sum_{j'=1}^p\sqrt{\mathbb{E}(\tilde r_{tmhj'}^2)+ (n-L)^{-1}\sum_{t=L+1}^{n}  \{\tilde r_{tmhj'}^2 - \mathbb{E}(\tilde r_{tmhj'}^2)\}} \\& =\sum_{h = 0}^L\sum_{j'=1}^p \sqrt{  I_{1, tmhj'} + I_{2, tmhj'}}.
\end{split}
\end{equation}
Recall that $\cov(\zeta_{tjl}, \zeta_{tjl'}) = \omega_{jl}^X I(l=l'),$ $\beta_{hj}(u,v) = \sum_{l,m=1}^\infty a_{hjlm} \psi_{jl}(u)\phi_m(v)$ and $|a_{hjlm}| \leq u_{hj}(l+m)^{-\kappa-1/2}.$ Then for $(h,j') \in S$,
\begin{equation} \nonumber
\begin{split}
 I_{1, tmhj'} &= \mathbb{E}[(\sum_{l'=q_1+1}^\infty\langle \psi_{j'l'},\langle\beta_{hj'},\widehat \phi_m \rrangle \zeta_{tj'l'}) ^2] = \sum_{l'=q_1+1}^\infty\langle \psi_{j'l'},\langle\beta_{hj'},\widehat \phi_m \rrangle^2 \omega_{j'l'}
    \\& \lesssim \sum_{l'=q_1+1}^\infty\langle  \psi_{j'l'},\langle\sum_{l'',m''=1}^\infty a_{hj'l''m''} \psi_{j'l''}\phi_{m''}\rangle,\phi_m + (\widehat \phi_m -\phi_m)\rrangle^2
        \\& \lesssim \sum_{l'=q_1+1}^\infty  a_{hj'l'm}^2 + \|\widehat \phi_m - \phi_m\|^2\sum_{l'=q_1+1}^\infty(\sum_{m''=1}^\infty a_{hj'l'm''})^2
        \\&\lesssim u_{hj'}^2 (q_1+m)^{-2\kappa} + u_{hj'}^2 \|\widehat \phi_m - \phi_m\|^2 q_1^{-2\kappa+2}.
\end{split}
\end{equation}
To provide the upper bound of $I_{2, tmhj'}$, we start with
\begin{equation} \nonumber
\begin{split}
& \frac{\sum_{t=L+1}^n [\zeta_{tj'l_1}\zeta_{tj'l_2} - \mathbb{E}(\zeta_{tj'l_1}\zeta_{tj'l_2})]}{n-L}
\\&= \langle \psi_{j'l_1}, \langle \widehat \Sigma_{0,j'j'}^X -\Sigma_{0,j'j'}^X, \psi_{j'l_2}\rrangle \leq  \|\widehat \Sigma_{0,j'j'}^X -\Sigma_{0,j'j'}^X\|_\cS = O_P\{\cM_1^X n^{-1/2} \}. 
\end{split}
\end{equation}
Combining this result with Lemmas~\ref{lm_guo2019_elementwise} and \ref{lm_guo2019_eigen} and following the similar argument in the proof of the upper bound of $I_{1, tmhj'}$, we obtain that, for $(h,j') \in S$,
\begin{equation} \nonumber
\begin{split}
 &\quad I_{2, tmhj'} \\&=  \sum_{l_1,l_2 =q_1+1}^\infty\langle \psi_{j'l_1},\langle\beta_{hj'},\widehat \phi_m \rrangle\langle \psi_{j'l_2},\langle\beta_{hj'},\widehat \phi_m \rrangle \frac{\sum_{t=L+1}^n [\zeta_{tj'l_1}\zeta_{tj'l_2} - \mathbb{E}(\zeta_{tj'l_1}\zeta_{tj'l_2})]}{n-L} 
 \\& \leq   \|\widehat \Sigma_{0,j'j'}^X -\Sigma_{0,j'j'}^X\|_\cS \{\sum_{l'=q_1+1}^\infty\langle \psi_{j'l'},\langle\beta_{hj'},\widehat \phi_m \rrangle\}^2= o_P(  I_{1, tmhj'}).
\end{split}
\end{equation}
Then
\begin{equation}
\begin{split}
     \|\frac{1}{n} \widehat{\bD}^{-1}\widehat{\bZ}^{\T} \widehat \bR\|_{\max}^{(q_1,q_2)} & \lesssim s \max_{1\leq j \leq p} \sqrt{\sum_{l=1}^{q_1}\sum_{m = 1 }^{q_2}\{(q_1+m)^{-2\kappa} +  \|\widehat \phi_m - \phi_m\|^2 q_1^{-2\kappa+2}\}}
     \\& \lesssim s \max_{1\leq j \leq p} \sqrt{q_1^{-2\kappa+2} + q_1^{-2\kappa+3}q_2\max_{1 \leq m \leq q_2} \|\widehat \phi_m - \phi_m\|^2 }
    \\& =O_P\{ s q_1^{-\kappa+1}\}.
     \end{split}
\end{equation}
$\square$

\section{Proofs of theoretical results in Section~\ref{sec_partial}}
\label{ap.thm.sec4}
This section is organized in the same manner as Appendix~\ref{ap.thm.sec3}. The proofs of Theorem~\ref{th_beta_partial} and Propositions~\ref{pr_partial_RE}--\ref{pr_partial_max.error_z} are presented in Appendices~\ref{ap.th.secC1} and \ref{ap.pro.secC2}, respectively, followed by supporting technical lemmas and their proofs in Appendix~\ref{ap.lemma.secC3}.

\subsection{Proof of Theorem~\ref{th_beta_partial} }
\label{ap.th.secC1}
Here $\widehat{B}\in \mathbb{R}^{pq}$ and $\widehat{\gamma} \in \mathbb{R}^d$ are  the minimizer of (\ref{target_partial}). Then
\begin{equation}
\nonumber
\begin{split}
&\frac{1}{2n}\|\mathcal{Y} -  \widehat{\mathcal{X}}\widehat D^{-1}\widehat B - \mathcal{Z}\widehat \gamma\|^2 + \lambda_{n1}\|\widehat B\|_1^{(q)}  + \lambda_{n2}\|\widehat \gamma\|_1\\\leq& \frac{1}{2n}\|\mathcal{Y} -  \widehat{\mathcal{X}}\widehat D^{-1}B - \mathcal{Z}\gamma\|^2  + \lambda_{n1}\|B\|_1^{(q)}+ \lambda_{n2}\|\gamma\|_1.
\end{split}
\end{equation}
Letting $\Delta = \widehat{B} - B,$ $\delta = \widehat \gamma - \gamma,$ $S_1^c$ be the complement of $S_1$ in the set $\{1,\dots,p\}$ and $S_2^c$ be the complement of $S_2$ in the set $\{1,\dots,d\}$ , we have
\begin{eqnarray*}
& &\frac{1}{2n} \{\Delta^{\T} \widehat \Omega^{\T} \widehat \Omega\Delta  + 2\Delta^{\T} \widehat\Omega^{\T} \mathcal{Z}\delta+  \delta^{\T}\mathcal{Z}^{\T} \mathcal{Z} \delta \}
\\ &\leq& \frac{1}{n} (\Delta^{\T} \widehat \Omega^{\T}+\delta^{\T} \mathcal{Z}^{\T})(\mathcal{Y}  - \widehat \Omega B - \mathcal{Z}\gamma)  
+ \lambda_{n1}(\|B\|_1^{(q)} - \|B+\Delta\|_1^{(q)})
\\ & &+ \lambda_{n2}(\|\gamma\|_1 - \|\gamma + \delta\|_1)
\\&\leq& \frac{1}{n} \Delta^{\T} \widehat \Omega^{\T}(\mathcal{Y}  - \widehat \Omega B - \mathcal{Z}\gamma) + \frac{1}{n} \delta^{\T} \mathcal{Z}^{\T}(\mathcal{Y} - \widehat \Omega B - \mathcal{Z}\gamma) 
 + \lambda_{n1}(\|\Delta_{S_1}\|_1^{(q)} - \|\Delta_{S_1^c}\|_1^{(q)}) 
\\& & + \lambda_{n2}(\|\delta_{S_2}\|_1 - \|\delta_{S_2^c}\|_1),
\end{eqnarray*}
where $\widehat{\Omega} = 
\widehat{\mathcal{X}}\widehat{D}^{-1}.$ By Propositions~\ref{pr_partial_max.error_X}, \ref{pr_partial_max.error_z} and the choice of $\lambda_n \asymp \lambda_{n1} \asymp \lambda_{n2} \geq 2 C_{0}^* s_1 (\cM_{X,Z} + \cM^\epsilon)[q^{\alpha+2}\{\log (pq+d)/n\}^{1/2}+ q^{-\kappa+1}],$ we obtain that
\begin{equation} \nonumber
\begin{split}
     \frac{1}{n} \Delta^{\T} \widehat \Omega^{\T}(\mathcal{Y}  - \widehat \Omega B - \mathcal{Z}\gamma) &\leq  \frac{1}{n}\|\Delta\|_1^{(q)}\|\widehat \Omega^{\T}(\mathcal{Y}  - \widehat \Omega B - \mathcal{Z}\gamma)\|_{\max}^{(q)} \\&\leq \frac{\lambda_{n}}{2}(\|\Delta_{S_1}\|_1^{(q)} + \|\Delta_{S_1^c}\|_1^{(q)}) ,
 \\     \frac{1}{n} \delta^{\T} \mathcal{Z}^{\T}(\mathcal{Y} - \widehat \Omega B - \mathcal{Z}\gamma) &\leq  \frac{1}{n} \|\delta\|_1 \| \mathcal{Z}^{\T}(\mathcal{Y} - \widehat \Omega B - \mathcal{Z}\gamma)\|_{\max} \\&\leq \frac{\lambda_{n}}{2}(\|\delta_{S_2}\|_1 + \|\delta_{S_2^c}\|_1).
\end{split}
\end{equation}
Combining the above results, we have 
\begin{equation} \nonumber
    0 \leq \frac{3}{2}( \|\Delta_{S_1}\|_1^{(q)} + \|\delta_{S_2}\|_1)  - \frac{1}{2}( \ \|\Delta_{S_1^c}\|_1^{(q)} + \|\delta_{S_2^c}\|_1).
\end{equation}
This ensures $ \|\Delta_{S_1^c}\|_1^{(q)}+ \|\delta_{S_2^c}\|_1 \leq 3( \|\Delta_{S_1}\|_1^{(q)}+\|\delta_{S_2}\|_1 ).$ Then we have that
\begin{equation} \nonumber
    \|\Delta\|_1^{(q)} + \|\delta\|_1  \leq 4(\|\Delta_{S_1}\|_1^{(q)}+\|\delta_{S_2}\|_1  ) \leq 4(  \sqrt{s_1} \|\Delta\| + \sqrt{s_2}\|\delta\| )  \leq 4\sqrt{s_1+ s_2} (\|\Delta\| + \|\delta\| ).
\end{equation}
This, together with Proposition~\ref{pr_partial_RE}, $\|\Delta\|_1 \leq \sqrt{q}\|\Delta\|_1^{(q)}$ and $\tau_2^* \geq 64\tau_1^*q (s_1+s_2)$ implies
\begin{equation}
\begin{split} \nonumber
    \frac{1}{n} \{\Delta^{\T} \widehat \Omega^{\T} \widehat \Omega\Delta  + 2\Delta^{\T} \widehat\Omega^{\T} \mathcal{Z}\delta+  \delta^{\T}\mathcal{Z}^{\T} \mathcal{Z} \delta \} &\geq \tau_2^*( \|\Delta\|^2 + \|\delta\|^2) -\tau_1^*(\sqrt{q}\|\Delta\|_1^{(q)}+\|\delta\|_1  )^2
    \\& \geq \frac{\tau_2^*}{2}(\|\Delta\|+\|\delta\|)^2 -\tau_1^*q( \|\Delta\|_1^{(q)}+ \|\delta\|_1  )^2
    \\& \geq \{\frac{\tau_2^*}{2}- 16\tau_1^*q (s_1+s_2) \}(\|\Delta\|+ \|\delta\| )^2  
    \\& \geq \frac{\tau_2^*}{4}(\|\Delta\|+ \|\delta\| )^2.
\end{split}
\end{equation}
This implies 
\begin{equation}\nonumber
    \frac{\tau_2^*}{8}(\|\Delta\|+ \|\delta\| )^2 \leq \frac{3\lambda_{n}}{2}( \|\Delta\|_1^{(q)}+ \|\delta\|_1  ) \leq 6 \lambda_{n} \sqrt{s_1 + s_2} (\|\Delta\|+ \|\delta\| ).
\end{equation}
Therefore, we obtain that
\begin{equation}\nonumber
    \|\Delta\|+ \|\delta\|  \lesssim \frac{\lambda_n\sqrt{s_1+s_2}}{\tau_2^*},
\end{equation}
\begin{equation}\nonumber
    \|\Delta\|_1^{(q)} +  \|\delta\|_1 \lesssim \frac{\lambda_n(s_1+s_2)}{\tau_2^*}.
\end{equation}
Provided that  $\|D^{-1}\|_{\max} \leq \alpha^{1/2}c_0^{-1/2}q^{\alpha/2}$, the rest can be proved in a similar way to the proof of Theorem~\ref{th_beta}, which shows
\begin{equation} \nonumber 
\begin{split}
\|\widehat{\eulB}-\eulB\|_{1} + q^{\alpha/2} \|\widehat \gamma - \gamma\|_1 &\leq 
\|\widehat \Psi- \Psi\|_1^{(q)} + q^{\alpha/2} \|\widehat \gamma - \gamma\|_1 + o(1)
\\& \leq \|D^{-1}\|_{\max} \|\widehat{B}-B\|_{1} + q^{\alpha/2} \|\widehat \gamma - \gamma\|_1 + o(1)
\\&\lesssim \frac{q^{\alpha/2} \lambda_n(s_1+s_2) }{\tau_2^*}\left\{1+o(1)\right\}.
\end{split}
\end{equation}
$\square$

\subsection{Proofs of propositions}

\paragraph{Proof of Proposition~\ref{pr_partial_RE}}
By Lemmas~\ref{lm_gamma}, \ref{lm_gamma_partial} and \ref{lm_basu_cov}, we obtain
\begin{equation} \nonumber
    \begin{split}
        &\|\frac{1}{n}\mathcal{S}^{\T}\mathcal{S}- \frac{1}{n}\mathbb{E}\{\mathcal{S}^{\T}\mathcal{S}\}\|_{\max}
        \\&= \max(\big\|\frac{1}{n} \mathcal{Z}^{\T} \mathcal{Z} - \frac{1}{n} \mathbb{E}\{\mathcal{Z}^{\T}   \mathcal{Z}\}\big\|_{\max},\big\|\frac{1}{n} \mathcal{Z}^{\T} \widehat \Omega - \frac{1}{n} \mathbb{E}\{\mathcal{Z}^{\T}  \Omega\}\big\|_{\max},\|\widehat \bGamma - \bGamma\|_{\max})
        \\& = O_P\{\max( \cM_1^Z \sqrt{\frac{\log (d)}{n}}, \cM_1^X q^{\alpha +1}\sqrt{\frac{\log (pq)}{n}}, \cM_{X,Z} q^{\alpha +1}\sqrt{\frac{\log (pqd)}{n}}) \}
        \\&= O_P\{ \cM_{X,Z} q^{\alpha +1}\sqrt{\frac{\log (pq+d)}{n}}\}.
    \end{split}
\end{equation}
 Combining this with Condition~\ref{con_partial_eigen_min_xz} and following the similar argument in the proof of Proposition~\ref{pr_fof_RE} implies Proposition~\ref{pr_partial_RE}. $\square$
\\
\paragraph{Proof of Proposition~\ref{pr_partial_max.error_X}}
\label{ap.pro.secC2}
Notice that
\begin{equation} \nonumber
    \frac{1}{n}  \widehat \Omega^{\T}(\mathcal{Y}  - \widehat \Omega B - \mathcal{Z}\gamma) =  \frac{1}{n}  \widehat \Omega^{\T}(   (\Omega - \widehat \Omega) B  +  R + E)
\end{equation}
where $\widehat{\Omega} = \widehat{\mathcal{X}}\widehat{D}^{-1},$ $  B =  D \Psi$ and $j$-th row of $ \Psi$ takes the form $\Psi_{j} = \int_\cU \bpsi_j(u) \beta_{j}(u) du.$  Recall that  $ r_{t} = \sum_{j=1}^p \sum_{l=q+1}^{\infty} \zeta_{tjl} \langle  \psi_{jl},\beta_{j}\rangle.$ 
Then it follows from Lemma~\ref{lm_gamma_2} when $L = 0$  that there exist some positive constants $C_{11}^*, c_1^*$ and $c_2^*$ such that 
\begin{equation} \nonumber
\begin{split}
    \|\frac{1}{n}\widehat \Omega^{\T} (\Omega - \widehat \Omega) B \|_{\max}^{(q)}  &\leq  \|\frac{1}{n}\widehat \Omega^{\T} (\Omega - \widehat \Omega)\|_{\max}^{(q)} \| B \|_1^{(q)} 
    \\ & \leq C_{11}^* s_1  \cM_1^X  q^{\alpha+2}\sqrt{\frac{\log(pq)}{n}},
\end{split}
\end{equation}
with probability greater than $1 - c_1^*(pq)^{-c_2^*}.$

Second, it follows from Lemma~$\ref{lm_Residual_partial_X}$ that there exist some positive constants $C_{12}^*, c_1^*$ and $c_2^*$ such that 
\begin{equation} \nonumber
    \|\frac{1}{n} \widehat \Omega^{\T}  R\|_{\max}^{(q)} \leq  C_{12}^*s_1q^{-\kappa+1} ,
\end{equation}
with probability greater than $1 - c_1^*(pq)^{-c_2^*}.$

Third, it follows from Proposition~\ref{pr_FPCscores_XE_partial} that there exist some positive constants $C_{13}^*, c_1^*$ and $c_2^*$ such that 
\begin{equation} \nonumber
    \|\frac{1}{n} \widehat \Omega^{\T} E\|_{\max}^{(q)} = \| \frac{1}{n}\widehat {D}^{-1}{D}{D}^{-1}\widehat{\mathcal{X}}^{\T} E\|_{\max}^{(q)} \leq C_{13}^* \{\mathcal{M}_1^{X} +\mathcal{M}^{\epsilon}\}q^{1/2}\sqrt{\frac{\log(pq)}{n}},
\end{equation}
with probability greater than $1 - c_1^*(pq)^{-c_2^*}.$

Combining the above results, we obtain that there exist some positive constants $C_{01}, c_1^*$ and $c_2^*$ such that
\begin{equation} \nonumber
    \frac{1}{n} \| \widehat \Omega^{\T}(\mathcal{Y}  - \widehat \Omega B - \mathcal{Z}\gamma) \|_{\max}^{(q)} \leq C_{01} s_1 (\cM_1^X + \cM^\epsilon)\{q^{\alpha+2}\sqrt{\frac{\log(pq)}{n}}+ q^{-\kappa+1}\}
\end{equation}
with probability greater than $1 - c_1^*(pq)^{-c_2^*}.$  
$\square$
\\

\paragraph{Proof of Proposition~\ref{pr_partial_max.error_z}}
Notice that
\begin{equation} \nonumber
 \frac{1}{n}  \mathcal{Z}^{\T}(\mathcal{Y} - \widehat \Omega B - \mathcal{Z}\gamma) = \frac{1}{n}  \mathcal{Z}^{\T}(   (\Omega - \widehat \Omega) B  +  R + E).
\end{equation}
First, we show the deviation bound of $\frac{1}{n}  \mathcal{Z}^{\T} (\Omega - \widehat \Omega) B.$ It follows from Lemma~\ref{lm_gamma_partial_2} and the fact that $\|\Psi_{j}\|_1 = \sum_{j = 1}^q u_j l^{-\kappa} = O(u_j)$, for $j \in S_1,$ that there exist some positive constants $C_{21}^*, c_1^*$ and $c_2^*$ such that 
\begin{equation} \nonumber
\begin{split}
    \|\frac{1}{n} \mathcal{Z}^{\T} (\Omega - \widehat \Omega) B \|_{\max}  &  \leq \|\frac{1}{n} \mathcal{Z}^{\T} (\Omega - \widehat \Omega)\|_{\max}\| B\|_1\
    \\& \leq \|\frac{1}{n} \mathcal{Z}^{\T} (\Omega - \widehat \Omega)\|_{\max}\|D\|_{\max}\| \Psi\|_1\
     \\ & \leq C_{21}^* s_1  \cM_{X,Z}  q^{\alpha+1}\sqrt{\frac{\log(pqd)}{n}},
\end{split}
\end{equation}
with probability greater than $1 - c_1^*({pqd})^{-c_2^*}.$

Second, it follows from Lemma~$\ref{lm_Residual_partial_Z}$ that there exist some positive constants $C_{22}^*, c_1^*$ and $c_2^*$ such that 
\begin{equation} \nonumber
    \|\frac{1}{n}  \mathcal{Z}^{\T}  R\|_{\max} \leq  C_{22}^* s_1q^{-\kappa+1/2},
\end{equation}
with probability greater than $1 - c_1^*(pqd)^{-c_2^*}.$

Third, it follows from Lemma~\ref{lm_basu_cov} that there exist some positive constants $C_{23}^*, c_1^*$ and $c_2^*$ such that 
\begin{equation} \nonumber
    \|\frac{1}{n}  \mathcal{Z}^{\T} E\|_{\max} \leq C_{23}^* \{\mathcal{M}_1^{Z} +\mathcal{M}^{\epsilon}\}\sqrt{\frac{\log(d)}{n}},
\end{equation}
with probability greater than $1 - c_1^*(d)^{-c_2^*}.$

Combining the above results, we obtain that there exist some positive constants $C_{02}, c_1^*$ and $c_2^*$ such that
\begin{equation} \nonumber
    \frac{1}{n} \| \mathcal{Z}^{\T}(\mathcal{Y} - \widehat \Omega B - \mathcal{Z}\gamma) \|_{\max} \leq C_{02} s_1 \{\cM_{X,Z}+ \cM^\epsilon)\}\{q^{\alpha+1}\sqrt{\frac{\log(pq+d)}{n}}+ q^{-\kappa+1/2}\}
\end{equation}
with probability greater than $1 - c_1^*(pq+d)^{-c_2^*}.$ $\square$

\subsection{Technical lemmas and their proofs}
\label{ap.lemma.secC3}
\begin{lemma} \label{lm_gamma_partial}
Suppose that Conditions~\ref{con_stability}--\ref{con_eigen} hold.  Then there exist some positive constants $\widetilde C_{1,Z\Gamma},$ $c_1^*$ and $c_2^*$ such that
\begin{equation} \nonumber
    	\big\|\frac{1}{n} \mathcal{Z}^{\T} \widehat \Omega - \frac{1}{n} \mathbb{E}\{\mathcal{Z}^{\T}  \Omega\}\big\|_{\max} \leq \widetilde C_{1,Z\Gamma} \cM_{X,Z} q^{\alpha +1}\sqrt{\frac{\log (pqd)}{n}}
\end{equation}
with probability greater than $1 - c_1^*(pqd)^{-c_2^*}.$
\end{lemma}
{\bf Proof}. Note that  
$$
\big\|\frac{1}{n} \mathcal{Z}^{\T} \widehat \Omega - \frac{1}{n} \mathbb{E}\{\mathcal{Z}^{\T}  \Omega\}\big\|_{\max} = \underset{\underset{1 \le l \le q}{1\le j\le p, 1\le k\le d}}{\max} {\left|\{\widehat\omega_{jl}^X \}^{-1/2}\widehat \varrho_{h,jkl}^{X,Z} - \{\omega_{jl}^X \}^{-1/2}\varrho_{h,jkl}^{X,Z}\right|}.
$$
Let $\widehat s_{jkl} = \left\{{\omega_{jl}^X}/{\widehat\omega_{jl}^X}\right\}^{1/2},$ then we obtain that
\begin{equation} \nonumber
    \{\widehat\omega_{jl}^X \}^{-1/2}\widehat \varrho_{h,jkl}^{X,Z} - \{\omega_{jl}^X \}^{-1/2}\varrho_{h,jkl}^{X,Z} = \widehat s_{jkl} \frac{\widehat\varrho_{h,jkl}^{X,Z}-\varrho_{h,jkl}^{X,Z}}{\{\omega_{jl}^X\}^{1/2}} + \frac{\{ \omega_{jl}^X\}^{1/2} -  \{\widehat\omega_{jl}^X\}^{1/2}}{\{ \widehat \omega_{jl}^X\}^{1/2}} \frac{\varrho_{h,jkl}^{X,Z}}{\{\omega_{jl}^X\}^{1/2}}.
\end{equation}
It follows Propositions~\ref{pr_FPCscores_XZ}, \ref{pr_fof_eigen} and the fact $\mathbb{E}({\zeta}_{tjl}Z_{tk}) \leq \sigma_{k}^Z \{\omega_{jl}^X \}^{1/2}$ that there exist some positive constants $\widetilde C_{1,Z\Gamma}, c_1^*$ and $c_2^*$ such that 
\begin{equation} \nonumber
\begin{split}
    	\big\|\frac{1}{n} \mathcal{Z}^{\T} \widehat \Omega - \frac{1}{n} \mathbb{E}\{\mathcal{Z}^{\T}  \Omega\}\big\|_{\max} &\leq \widetilde C_{1,Z\Gamma} \cM_{X,Z} q^{\alpha +1}\sqrt{\frac{\log (pqd)}{n}}
\end{split}
\end{equation}
with probability greater than $1 - c_1^*(pqd)^{-c_2^*}.$  $\square$
\\

\begin{lemma} \label{lm_gamma_partial_2}
Suppose that Conditions~\ref{con_stability}--\ref{con_eigen} hold.  Then there exist some positive constants $\widetilde C_{2,Z\Gamma},$ $c_1^*$ and $c_2^*$ such that
\begin{equation} \nonumber
    \|\frac{1}{n} \mathcal{Z}^{\T} (\Omega - \widehat \Omega)\|_{\max} \leq \widetilde C_{2,Z\Gamma} \cM_{X,Z} q^{\alpha +1}\sqrt{\frac{\log (pqd)}{n}}
\end{equation}
with probability greater than $1 - c_1^*(pqd)^{-c_2^*}.$
\end{lemma}
{\bf Proof}. We first consider $\|\frac{1}{n} \mathcal{Z}^{\T} \Omega - \frac{1}{n} \mathbb{E}\{\mathcal{Z}^{\T}  \Omega\}\|_{\max}$. By (\ref{eq_th_XZ}) in Appendix~\ref{ap.pro.secA2}, we obtain that
\begin{equation} \nonumber
\begin{split}
    &\max_{j,k,m}\frac{(n-L)^{-1} \sum_{t=L+1}^{n} Z_{(t-h)j} \zeta_{tkm}}{\sqrt{\omega_{km}^X}} -\frac{\mathbb{E}(Z_{(t-h)j} \zeta_{tkm})}{\sqrt{ \omega_{km}^X}}
    \\&= \max_{j,k,m}\frac{\langle \widehat \Sigma_{h,jk}^{Z,X}, \psi_{km} \rangle -  \langle \Sigma_{h,jk}^{Z,X}, \psi_{km} \rangle  }{\sqrt{ \omega_{km}^X}} =O_P\{\cM_{X,Z} \sqrt{\frac{\log (pqd)}{n}}\}.
\end{split}
\end{equation}
This, together with Lemma~\ref{lm_gamma_partial}, implies that 
\begin{equation}\nonumber
\begin{split}
      \|\frac{1}{n} \mathcal{Z}^{\T} (\Omega - \widehat \Omega)\|_{\max} 
      &
      \leq   \|\frac{1}{n} \mathcal{Z}^{\T} \Omega - \frac{1}{n} \mathbb{E}\{\mathcal{Z}^{\T}  \Omega\}\|_{\max}+ \|\frac{1}{n} \mathcal{Z}^{\T} \widehat \Omega - \frac{1}{n} \mathbb{E}\{\mathcal{Z}^{\T}  \Omega\}\|_{\max}  
      \\&= O_P\{\cM_{X,Z} q^{\alpha+1 }\sqrt{\frac{\log (pqd)}{n}}\}.
\end{split}
\end{equation}
$\square$
\\
\begin{lemma} \label{lm_Residual_partial_X}
Suppose that Conditions~\ref{con_stability}--\ref{con_eigen} and \ref{con_partial_beta_b} hold. Then there exist some positive constants $C_{R1},$ $c_1^*$ and $c_2^*$ such that
\begin{equation} \nonumber
\begin{split}
     \|n^{-1} \widehat \Omega^{\T}  R\|_{\max}^{(q)}  \leq C_{R1}s_1q^{-\kappa+1}
     \end{split}
\end{equation}
with probability greater than $1 - c_1^*(pq)^{-c_2^*}.$
\end{lemma}
{\bf Proof.} This lemma can be proved in a similar way to Lemma~\ref{lm_Residual} and hence the proof is omitted here.  $\square$
\\
\begin{lemma} \label{lm_Residual_partial_Z}
Suppose that Conditions~\ref{con_stability}--\ref{con_eigen} and \ref{con_partial_beta_b} hold. Then there exist some positive constants $C_{R2},$ $c_1^*$ and $c_2^*$ such that
\begin{equation} \nonumber
\begin{split}
     \|n^{-1} \mathcal{Z}^{\T}  R\|_{\max}  \leq C_{R2}s_1q^{-\kappa+1/2}
     \end{split}
\end{equation}
with probability greater than $1 - c_1^*(pqd)^{-c_2^*}.$
\end{lemma}
{\bf Proof.} This lemma can be proved in a similar way to Lemma~\ref{lm_Residual} and hence the proof is omitted here.  $\square$

\section{Existing results for sub-Gaussian (functional) linear processes}
\label{ap.subgaussian}
For ease of reference, we present some useful existing results in \cite{guo2019}, including non-asymptotic error bounds on estimated covariance matrix function, estimated eigenpairs and estimated (auto)covariance between estimated FPC scores. By Theorem~\ref{th_XY}, we can easily extend these results from Gaussian functional time series to accommodate sub-Gaussian functional linear processes in Lemmas~\ref{lm_guo2019_elementwise}--\ref{lm_guo2019_score}. Moreover, we also present non-asymptotic error bounds on estimated (cross-)covariance matrix in \cite{basu2015a} to accommodate sub-Gaussian linear processes in Lemma~\ref{lm_basu_cov}.
\begin{lemma} \label{lm_guo2019_elementwise}
Suppose that Conditions~\ref{con_stability}, \ref{con_sub_coefficient} and \ref{con_e} hold for sub-Gaussian linear process $\{\bX_t(\cdot)\}_{t \in \mathbb{Z}}.$ Then there exists some universal constant $\tilde c_1>0$ such that for any $\eta>0$ and each $j,k = 1,\dots,p,$
\begin{equation} \nonumber
P\left\{\| \widehat{\Sigma}_{0,jk}^{X}-\Sigma_{0,jk}^{X}\|_\cS > 2\omega_0^X\mathcal{M}_1^{X} \eta\right\} \leq 4\exp\{-\tilde c_1 n\min(\eta^2,\eta)\}.
\end{equation}
\end{lemma}
\textbf{Proof}. This lemma follows directly from Theorem~\ref{th_XY} and Theorem~2 of \cite{guo2019} and hence the proof is omitted here.
$\square$
\\
\begin{lemma} \label{lm_guo2019_eigen}
Suppose that Conditions~\ref{con_stability}, \ref{con_sub_coefficient}, \ref{con_e} and \ref{con_eigen} hold for sub-Gaussian linear process $\{\bX_t(\cdot)\}_{t \in \mathbb{Z}}.$ Let $M$ be a positive integer possibly depending on $(n,p).$ If $n \gtrsim \log(pM)M^{4\alpha+2}(M_1^X)^2,$ then there exist some constants $\tilde c_2, \tilde c_3>0$ such that, with probability greater than $1 - \tilde c_2 (pM)^{-\tilde c_3},$ the estimates $\{\widehat \omega_{jl}^X\}$ and $\{\widehat \psi_{jl}\}$ satisfy
	\begin{equation} \label{eq_guo_eigen}
	\begin{split}
	\max_{1\leq j\leq p, 1 \leq l \leq M}\left\{\Big|\frac{\widehat \omega_{jl}^X - \omega_{jl}^X}{\omega_{jl}^X}\Big| + \Big\|\frac{\widehat \psi_{jl} - \psi_{jl}}{l^{\alpha+1}}\Big\|\right\} \lesssim \cM_1^X \sqrt{\frac{\log(pM)}{n}}.
	\end{split}
	\end{equation}
\end{lemma}
\textbf{Proof}. This lemma follows directly from Theorem~\ref{th_XY} and Theorem~3 of \cite{guo2019} and hence the proof is omitted here.
$\square$
\\

\begin{lemma} \label{lm_guo2019_phi_g}
Suppose that conditions in Lemma~\ref{lm_guo2019_eigen} hold. Then there exists some universal constant $\tilde c_4>0$ such that for each $j=1,\dots,p, l=1, \dots,d_j,$ any given function $g \in \cH$ and $\eta > 0,$
	\begin{equation} \nonumber
	\begin{split}
	&P\left\{\left|\big \langle \widehat \psi_{jl} - \psi_{jl}, g\big \rangle \right| \ge
	\tilde \rho_1\|g^{-jl}\|_{\omega}  \cM_1^{X} \{\omega_{jl}^X\}^{1/2} l^{\alpha + 1}  \eta
	+ \tilde \rho_2 \|g\|\{\cM_1^X\}^2  l^{2(\alpha+1)}\eta^2 \right\} \\
	&\le 8 \exp\Big\{ - \tilde c_4 n\min(\eta^2, \eta)\Big\} + 4\exp\Big\{ - \tilde c_4  \{\cM_1^X\}^{-2} n l^{-2(\alpha + 1)}\Big\},
	\end{split}
	\end{equation}
where $g(\cdot)=\sum_{l=1}^{\infty}g_{jl}\psi_{jl}(\cdot),$ $\|g^{-jl}\|_{\omega} = \big({\sum}_{l': l'\neq l}\omega_{jl'} g_{jl'}^2 \big)^{1/2},$
$\tilde \rho_{1} = 2 c_0^{-1}\omega_0^X
$
and
$\tilde \rho_{2} = 4(6+2\sqrt{2}) c_0^{-2} \{\omega_0^X\}^2$ with $c_0 \leq 4\cM_1^X\omega_0^X l^{\alpha+1}.$

\end{lemma}
\textbf{Proof}. This lemma follows directly from Theorem~\ref{th_XY} and Lemma~3 of \cite{guo2019} and hence the proof is omitted here.
$\square$
\\

\begin{lemma} \label{lm_guo2019_score}
Suppose that conditions in Lemma~\ref{lm_guo2019_eigen} hold. Let $M$ be a positive integer possibly depending on $(n,p).$ If $n \gtrsim \log(pM)M^{4\alpha+2}(\cM_1^X)^2,$ then there exist some constants $\tilde c_5, \tilde c_6>0$ such that, with probability greater than $1 - \tilde c_5 (pM)^{-\tilde c_6},$ the estimates $\{\widehat \sigma_{h,jklm}^X\}$ satisfies
\begin{equation} \label{eq_guo_scores}
\underset{\underset{1 \le l,m \le M}{1\le j, k \le p}}{\max} \frac{\left|\widehat \sigma_{h,jklm}^{X} - \sigma_{h,jklm}^{X}\right|}{(l \vee m)^{\alpha+1}{\sqrt{\omega_{jl}^X \omega_{km}^X}}}\lesssim \mathcal{M}_1^{X}\sqrt{\frac{\log(pM)}{n}}.
\end{equation}
\end{lemma}
\textbf{Proof}. This lemma follows directly from Theorem~\ref{th_XY} and Theorem~4 of \cite{guo2019} and hence the proof is omitted here.
$\square$
\\
\begin{lemma} \label{lm_basu_cov}
(i)Suppose $\{\bZ_t\}$ is from $d$-dimensional sub-Gaussian linear process with absolute summable coefficients
and bounded $\cM^Z.$ For any given vector $\bnu\in \mathbb{R}_0^d$ with $\|\bnu\|_0 \leq k$ $(k =1,\dots,d),$ denote $\cM(\bbf_Z,\bnu) = 2\pi \cdot\text{ess} \sup_{\theta \in [-\pi,\pi]} \bnu^{\T}\bbf_Z\bnu $. Then there exists some constants $ c, \tilde c_{16}, \tilde c_{17} > 0$ such that for any $\eta > 0,$
\begin{equation} \nonumber
     P\left\{\left|{\bnu^{\T}(\widehat{\bSigma}_{0}^Z - \bSigma_{0}^Z)\bnu}\right| > \cM(\bbf_Z,\bnu)\eta\right\} \leq 2 \exp \left\{ - cn \min \left(\eta^2,\eta\right)\right\},
\end{equation}
and
\begin{equation} \nonumber
     P\left\{\left|\frac{\bnu^{\T}(\widehat{\bSigma}_{0}^Z - \bSigma_{0}^Z)\bnu}{\bnu^{\T}\bSigma_{0}^Z\bnu}\ \right| > \cM_{k}^Z\eta\right\} \leq 2 \exp \left\{ -cn \min \left(\eta^2,\eta\right)\right\}.
\end{equation}
In particular, with probability greater than $1 - \tilde c_{16}(d)^{-\tilde c_{17}},$
\begin{equation} \nonumber
\max_{1\leq j,k\leq d}| \widehat{\Sigma}_{0,jk}^{Z}-\Sigma_{0,jk}^{Z}|   \lesssim \mathcal{M}_1^{Z}\sqrt{\frac{\log(d)}{n}}.
\end{equation}
(ii)Suppose $\{\epsilon_t\}$ is from sub-Gaussian linear process with absolute summable coefficients, bounded $\cM^\epsilon$ and independent of $\{\bZ_t\}.$ Then there exist some positive constants $\tilde c_{18}, \tilde c_{19}$ such that with probability greater than $1 - \tilde c_{18}(d)^{-\tilde c_{19}},$
\begin{equation} \nonumber
    \max_{1\leq j \leq d} \left|\sum_{t= 1}^n Z_{tj} \epsilon_t/n\right| \lesssim  (\mathcal{M}_1^{Z} +\mathcal{M}^{\epsilon})\sqrt{\frac{\log(d)}{n}}.
\end{equation}
\end{lemma}
\textbf{Proof}. This lemma can be proved in similar way to Proposition~2.4 of \cite{basu2015a} and be extended to sub-Gaussian linear process setting following the similar techniques used in the proof of Theorem~\ref{th_XY}.
$\square$

\section{Matrix representation of model~(\ref{model_fully})}
\label{sec_sm_fully_matrix}
It follows from the Karhunen-Lo\`eve expansion that model~(\ref{model_fully}) can be rewritten as 
\begin{equation}\nonumber
    \sum_{m=1}^\infty\xi_{tm}\phi_m(v) = \sum_{h=0}^L\sum_{j=1}^p\sum_{l=1}^\infty\langle\psi_{jl}(u),\beta_{hj}(u,v)\rangle\zeta_{(t-h)jl} + \epsilon_t(v),
\end{equation}
This, together with orthonormality of \{$\phi_m(\cdot)\}_{m\geq 1},$ implies that
\begin{equation}\nonumber
    \xi_{tm}= \sum_{h=0}^L\sum_{j=1}^p\sum_{l=1}^{q_{1j}}\langle\langle\psi_{jl}(u),\beta_{hj}(u,v)\rangle,\phi_m(v) \rangle\zeta_{(t-h)jl} + r_{tm} + \epsilon_{tm},
\end{equation}
where $r_{tm} = \sum_{h=0}^L\sum_{j=1}^p\sum_{l=q_{1j}+1}^\infty\langle\langle\psi_{jl}(u),\beta_{hj}(u,v)\rangle,\phi_m(v) \rangle\zeta_{(t-h)jl}$ and $\epsilon_{tm} = \langle\phi_m,\epsilon_t\rangle$ for $m = 1,\dots,q_2,$ represent the approximation and random errors, respectively. Let $\br_{t} = (r_{t1},\dots,r_{tq_2})^{\T}$ and $\beps_{t} = (\epsilon_{t1},\dots,\epsilon_{tq_2})^{\T}.$ Let $\bR$ and $\bE$ be $(n-L)\times q_2$ matrices whose row vectors are formed by $\{\br_{t},t=L+1,\dots,n\}$ and $\{\beps_{t},t=L+1,\dots,n\}$ respectively. Then (\ref{model_fully}) can be represented in the matrix form of (\ref{model_fully_matrix}).

\linespread{1.05}\selectfont
\bibliography{paperbib}
\bibliographystyle{dcu}

\end{document}